\documentclass[11pt]{amsart}
\usepackage{mathrsfs}
\usepackage{amssymb}
\usepackage{titletoc}
\usepackage{amscd}
\usepackage{amsmath, amssymb}
\usepackage{amsfonts}
\usepackage[colorlinks,linkcolor=blue,citecolor=blue, pdfstartview=FitH]{hyperref}
\usepackage[all,cmtip]{xy}

\usepackage{stmaryrd}
\usepackage{upgreek}
\usepackage{tikz-cd}
\usepackage{bm}

\numberwithin{equation}{section}

\theoremstyle{plain}
\newtheorem{thm}{Theorem}[section]

\newtheorem{lem}[thm]{Lemma}
\newtheorem{corollary}[thm]{Corollary}
\newtheorem{prop}[thm]{Proposition}

\theoremstyle{definition}

\newtheorem{rmk}[thm]{Remark}

\newtheorem{claim}[thm]{Claim}

\newtheorem{defn-thm}[thm]{Definition-Theorem}
\newtheorem{defn-pro}[thm]{Definition-Proposition}

\usepackage{fancyhdr}

\newcommand{\sA}{{\mathcal A}}

\newcommand{\D}{{\mathsf D}}

\newcommand{\F}{{\mathbb F}}

\newcommand{\sF}{{\mathcal F}}

\newcommand{\sG}{{\mathcal G}}

\newcommand{\sH}{{\mathcal H}}
\newcommand{\I}{{\mathrm I}}
\newcommand{\sI}{{\mathcal I}}

\renewcommand{\k}{\Bbbk}

\newcommand{\sL}{{\mathcal L}}

\newcommand{\sM}{{\mathcal M}}

\renewcommand{\P}{{\mathsf P}}

\newcommand{\sP}{{\mathcal P}}
\newcommand{\Q}{{\mathbb Q}}

\newcommand{\fS}{{\mathfrak{S}}}

\newcommand{\sW}{{\mathcal W}}
\newcommand{\X}{{\mathbb X}}
\newcommand{\Xd}{\check{\X}}

\newcommand{\Z}{{\mathbb Z}}

\newcommand{\id}{\mathrm{id}}

\newcommand{\xs}{{\ \xrightarrow{\sim}\ }}

\newcommand{\Hom}{\mathrm{Hom}}

\newcommand{\cha}{{\mathrm{char}}}

\newcommand{\Mod}{\text{-}\mathrm{Mod}}

\newcommand{\Rep}{\mathrm{Rep}}

\newcommand{\Vect}{\mathsf{Vect}}

\newcommand{\Spec}{\mathrm{Spec}}

\newcommand{\rk}{\mathrm{rk}} 
 
\newcommand{\pt}{\mathrm{pt}}

\newcommand{\Gm}{{\mathbb{G}_{\mathrm{m}}}}
\newcommand{\Ga}{{\mathbb{G}_{\mathrm{a}}}}

\newcommand{\Ad}{\mathrm{Ad}}

\newcommand{\For}{\mathsf{For}}

\newcommand{\semi}{{{\frac{\infty}{2}}}}
\newcommand{\quar}{{{\frac{\infty}{4}}}}
\newcommand{\oblv}{\mathrm{oblv}}
\newcommand{\Av}{\mathrm{Av}}

\newcommand{\IW}{\mathcal{IW}}

\newcommand{\simto}{\xrightarrow{\sim}}
\newcommand{\Phis}{\Phi^{\mathrm{s}}}
\newcommand{\ex}{\mathrm{ext}}
\newcommand{\af}{\mathrm{aff}}
\newcommand{\Shv}{\mathsf{Shv}}
\newcommand{\Shvc}{\mathsf{Shv}_{\mathrm{c}}}
\newcommand{\Shvfg}{\mathsf{Shv}_{\mathrm{fg}}}
\newcommand{\Ras}{\mathrm{Ras}}
\newcommand{\Fl}{\mathrm{Fl}}
\newcommand{\Gr}{\mathrm{Gr}}
\newcommand{\tFl}{\widetilde{\Fl}}

\newcommand{\bs}{\mathbf{s}}
\newcommand{\Iu}{\mathrm{I}_{\mathrm{u}}}

\newcommand{\comp}{\mathrm{c}}
\newcommand{\rS}{\mathfrak{S}}
\newcommand{\uFl}{\underline{\Fl}}
\newcommand{\fd}{\mathrm{fd}}
\newcommand{\ZCen}{\mathsf{Z}}
\newcommand{\Gait}{\mathrm{Ga}}
\newcommand{\Sat}{\mathsf{Sat}}

\DeclareMathOperator*{\colim}{colim}

\begin{document}

\title{Parabolic semiinfinite sheaves on affine flag varieties}

\makeatletter
\let\MakeUppercase\relax
\makeatother

\author{Simon Riche}
\address{Universit\'e Clermont Auvergne, CNRS, LMBP, F-63000 Clermont-Ferrand, France, Institut Universitaire de France.}
\email{simon.riche@uca.fr}

\author{Quan Situ}
\address{Institute of Science and Technology Austria (ISTA), Am Campus 1, 3400 Klosterneuburg, Austria.}
\email{quan.situ@ist.ac.at}

\thanks{
This project has received
funding from the European Research Council (ERC) under the European Union's Horizon 2020
research and innovation programme (grant agreement No.~101002592).
}

\date{}

\begin{abstract}
We study sheaves on affine flag varieties of reductive algebraic groups in parabolic semiinfinite settings, generalizing some constructions and results from \cite{Gaits18,adr} (for the case of the Borel subgroup) and~\cite{DL23} (for characteristic-$0$ coefficients and a particular case of parabolic semiinfinite condition). In particular we define perverse t-structures and Raskin equivalences in this context, and study an analogue of the Gaitsgory sheaf. In future work, these $\infty$-categories will be related to representation-theoretic $\infty$-categories associated with the Langlands dual group.
\end{abstract}

\maketitle

{\small{\setcounter{tocdepth}{1} \tableofcontents}}

\section{Introduction} 

\subsection{Semiinfinite sheaves on affine flag varieties}

This paper is a continuation of a joint project with Pramod Achar and Gurbir Dhillon aiming at developing a framework of \emph{semiinfinite sheaves} on affine flag varieties of connected reductive algebraic groups and exploring some applications to combinatorial, geometric and representation-theoretic questions involving the reductive group under consideration or its Langlands dual.
This study has its roots in predictions in the 1990s (by Lusztig and Finkelberg--Mirkovi{\'c}) that constructible sheaves on \emph{semiinfinite affine flag varieties} should provide geome\-tric models for some categories of representations of the Langlands dual Lie algebra; see~\cite[\S 1.1]{adr} for precise references. 
The concrete realizations of these predictions have been hindered by the difficulty of making sense of such categories of sheaves, but a new construction has been proposed more recently (with different applications in mind) by Gaitsgory~\cite{Gaits18}, in Lurie's framework of (presentable stable) $\infty$-categories, and we expect that this construction will allow us to realize these geometric models and use them more easily.

The first step of this project has been accomplished in~\cite{adr}, where Pramod Achar, Gurbir Dhillon and the first author have explained how to adapt Gaitsgory's construction, initially written for sheaves on affine Grassmannians with coefficients in characteristic $0$, to sheaves on affine flag varieties with coefficients of arbitrary characteristic, and explored some basic properties of these $\infty$-categories. An application in geometry has been explained in~\cite{adr-Drinfeld}, and first applications to representation theory will be developed in~\cite{adr-Gr,adrs}.

\subsection{Parabolic versions}

The goal of the present paper is to explain how to construct some variants of the main players of~\cite{adr} associated with parabolic subgroups, and generalize many of the results of~\cite{adr} to that setting.
We also take this opportunity to work out in detail the version of these constructions on general \emph{partial} affine flag varieties. (In the setting of~\cite{adr}, the possibility of this extension was mentioned, but the concrete details were left to the reader except in the important special case of affine Grassmannians.) Parabolic semiinfinite sheaves (on affine Grassmannians, and for coefficients of characteristic $0$) have already been investigated by Gurbir Dhillon and Sergey Lysenko in~\cite{DL23}, which provided some suggestions of constructions and statements; but even in their setting our versions are more general, and many of our proofs are different.

Let us now explain what we mean more precisely. Let $G$ be a connected reductive algebraic group over an algebraically closed field $\mathbb{F}$, and let $T \subset B \subset G$ be a maximal torus and a Borel subgroup. The choice of $B$ determines an Iwahori subgroup $\I$ in the loop group $\sL G$ of $G$, and we can consider the affine flag variety $\Fl=\sL G / \I$.
Let $W_\ex$ be the extended affine Weyl group of $G$, and $S_\af \subset W_\ex$ be the subset of simple reflections determined by $B$. Then for any finitary subset $K \subset S_\af$ we have a parahoric subgroup $\I_K \subset \sL G$ containing $\I$, and the associated partial affine flag variety $\Fl_K=\sL G/\I_K$. Given now a subset $J$ of the set $S \subset S_\af$ of simple reflections in the \emph{finite} Weyl group $W$, we have a standard parabolic subgroup $P_J \subset G$, with a Levi decomposition $P_J = M_J \ltimes U_J$. Any further subset $L \subset J$ is a finitary subset of the set of affine simple reflections of $M_J$, so we have an associated parahoric subgroup $\I_{L \subset J} \subset \sL M_J$, and the semiinfinite parabolic subgroup of $\sL G$ defined by $L \subset J$ is
\[
\I^\semi_{L \subset J} = \I_{L \subset J} \ltimes \sL U_J.
\]

This construction recovers several well-known and well-studied objects:
\begin{itemize}
\item
if $L=J=\varnothing$ one obtains the semiinfinite Iwahori subgroup appearing in particular in~\cite{Gaits18,adr};
\item
if $J=S$ we have $\I_{L \subset S} = \I_L$;
\item
if $L=J$ one obtains the parabolic semiinfinite subgroups appearing in~\cite{DL23};
\item
if $L=\varnothing$ we obtain the ``quarterinfinite'' parahoric subgroups that appear in the discussion in~\cite[\S 8.11]{bez-losev}.
\end{itemize}

\subsection{Contents}

Our main objects of study in this paper are the $\infty$-categories
\begin{equation}
\label{eqn:intro-cat}
\Shv \left( \I^\semi_{L \subset J} \backslash \Fl_K \right)
\end{equation}
of \'etale $\I^\semi_{L \subset J}$-equivariant sheaves on $\Fl_K$, with coefficients in a field $\k$ (of characteristic either $0$ or a prime invertible in $\F$), for the various possible choices of $L$, $J$ and $K$, and the natural functors that relate them. In particular:
\begin{enumerate}
\item
We study the relevant combinatorics of $W_\ex$, defining in particular a quarterinfinite length and a quarterinfinite order associated with $J$; see Section~\ref{sec:Wext}.
\item
We define the standard and costandard objects in the $\infty$-categories~\eqref{eqn:intro-cat}, and study their behaviour under various natural functors; see Section~\ref{sec:quar-semi-sheaves}.
\item
We define and study a ``Raskin equivalence'' relating~\eqref{eqn:intro-cat} to the $\infty$-category $\Shv \left( \I_{L} \backslash \Fl_K \right)$; see Section~\ref{sec:Raskin}. (This construction is a direct generalization of constructions in~\cite{Gaits18,adr,DL23}, which have their roots in earlier work of Raskin.)
\item
We define a ``perverse'' t-structure on~\eqref{eqn:intro-cat}, prove t-exactness of various natural functors, prove that standard and costandard objects are perverse in some cases, and describe the t-structure on $\Shv \left( \I \backslash \Fl_K \right)$ obtained by transporting this t-structure via the Raskin equivalence (for $L=\varnothing$) in terms similar to the description of the ``new t-structure'' of Bezrukavnikov--Lin; see Section~\ref{sec:perv}.
\item
We study a natural parabolic analogue of the ``Gaitsgory sheaf'' (or semiinfinite intersection cohomology sheaf) considered in~\cite{Gaits18,adr}; see Section~\ref{sec:Gaitsgory}.
\end{enumerate}
All of these results are natural extensions of statements appearing in~\cite{Gaits18,adr,DL23}, and their proofs are mostly based on similar arguments (with sometimes additional technical difficulties).

The only section containing a result which might be considered more original is Section~\ref{sec:t-exactness}. Here we prove t-exactness of a certain functor
\[
\Shv \left( \I^\semi_{J \subset J} \backslash \Fl_K \right) \to \Shv \left( \I^\semi_{\varnothing \subset \varnothing} \backslash \Fl_K \right)
\]
defined by combining two different choices of Raskin equivalences. The validity of this statement is suggested by expected relations with representation theory described in~\S\ref{ss:into-RT} below; but in fact we are planning to use it to prove some of these conjectural equivalences. This proof requires the consideration of analogues of the $\infty$-categories~\eqref{eqn:intro-cat} where the equivariance condition on the right is replaced by a Whittaker-type condition; for this reason we also study such $\infty$-categories systematically.

\subsection{Relation to representation theory}
\label{ss:into-RT}

Our main motivation for the systematic study of the $\infty$-categories $\Shv \left( \I^\semi_{L \subset J} \backslash \Fl_K \right)$ is that we expect them to provide useful ``geometric models'' for several interesting $\infty$-categories of representations of objects related to the Langlands dual reductive group over $\k$. Namely, recall that the \emph{geometric Satake equivalence}~\cite{mv} provides an equivalence of abelian monoidal categories
\begin{equation}
\label{eqn:intro-geom-Satake}
(\Shvfg(\sL^+ G \backslash \Gr)^\heartsuit, \star^{\sL^+ G}) \simeq (\Rep^{\heartsuit,\fd}(G^\vee_\k), \otimes)
\end{equation}
between (finitely generated) $\sL^+G$-equivariant perverse sheaves on the affine Grassmannian $\Gr=\Fl_S = \sL G / \sL^+ G$ (where $\I_S = \sL^+ G$ is the positive loop group) and finite-dimensional representations of the split reductive Langlands dual group $G^\vee_\k$ over $\k$. 

In the particular case $K=S$ and $L=J$ we will construct in~\cite{adrs} an equivalence of abelian monoidal categories
\begin{equation}
\label{eqn:intro-Satake-parabolic}
 \Shv_\omega \left( \I^\semi_{J \subset J} \backslash \Gr \right)^\heartsuit \simeq \Rep^{\heartsuit,\fd}(P_{J,\k}^\vee)
\end{equation}
where $P^\vee_{J,\k} \subset G^\vee_\k$ is the standard parabolic subgroup associated with $J$ (now seen as a subset of the simple reflections of the Weyl group of $G^\vee_\k$, which canonically identifies with $W$), and the symbol ``$\omega$'' indicates a full subcategory introduced formally in the body of the paper. Here the left-hand side is equipped with the fusion product. In the special case $J=S$ this equivalence will coincide with~\eqref{eqn:intro-geom-Satake}. It will generalize an equivalence conjectured in~\cite[\S 1.3.3]{DL23}\footnote{The parabolic subgroup appearing in~\cite[\S 1.3.3]{DL23} is the \emph{negative} parabolic subgroup attached to $J$. The precise relation between the two equivalences should be given by conjugation by the longest element $w_0$ in the Weyl group. For instance, for $J=\varnothing$, in our version the simple object $\overline{\Delta}^\semi_\lambda$ of~\cite[Proposition~6.20]{adr} will correspond to the $1$-dimensional representation of the (positive) Borel subgroup $P^\vee_{\varnothing,\k}$ corresponding to $\lambda$, while the version considered in~\cite{DL23} or~\cite[Proposition~1.5.7]{Gaits18} involves the action of $w_0$.} and the special case $J=\varnothing$ studied (under some assumptions on $\mathrm{char}(\k)$, and without the monoidal structure) in~\cite{adr-Gr}.

Now we turn to the $\infty$-category $\Shv \left( \I^\semi_{J \subset J} \backslash \tFl \right)$, where $\tFl$ is the quotient of $\sL G$ by the pro-unipotent radical $\Iu$ of $\I$. (Sheaves on $\tFl$ are also studied in the body of the paper, in parallel with the $\infty$-categories~\eqref{eqn:intro-cat}.) We assume that $\k$ has positive characteristic $\ell \geq h$ where $h$ is the Coxeter number of $G$, and consider the Lie algebra $\mathfrak{g}^\vee_\k$ of $G^\vee_\k$. 
Following~\cite{adr-cat} we consider the $\infty$-category
\[
(\mathfrak{g}^\vee_\k,P^\vee_{J,\k})\Mod
\]
of $(\mathfrak{g}^\vee_\k,P^\vee_{J,\k})$-modules, and its block $(\mathfrak{g}^\vee_\k,P^\vee_{J,\k})\Mod_0$ containing the trivial module. (Here $(\mathfrak{g}^\vee_\k,P^\vee_{J,\k})\Mod$ is a presentable stable $\infty$-category equipped with a t-structure whose heart is the abelian category of $P^\vee_{J,\k}$-equivariant $U(\mathfrak{g}^\vee_\k)$-modules such that the differential of the action of $P^\vee_{J,\k}$ coincides with the restriction of the action of $\mathfrak{g}^\vee_\k$.) Generalizing the picture sketched in~\cite[\S 1.4]{adr} (corresponding to $J=\varnothing$), we expect that for any $J$ there exists a t-exact equivalence of $\infty$-categories
\begin{equation}\label{equ:parabolic-cat-O}
\Shv \left( \I^\semi_{J \subset J} \backslash \tFl \right) \simeq (\mathfrak{g}^\vee_\k,P^\vee_{J,\k})\Mod_0,
\end{equation}
and that these equivalences satisfy some natural compatibility properties with each other and with some other constructions like~\eqref{eqn:intro-Satake-parabolic}. We expect to prove at least some particular cases of such equivalences in future work. 
An abelian version of the equivalence \eqref{equ:parabolic-cat-O} will be shown in an upcoming work of the second author with Peter Fiebig and Jonathan Gruber \cite{FGS}. 

Finally, in~\cite[\S 8.11]{bez-losev}, Bezrukavnikov and Losev conjecture an equivalence relating a (characteristic-$0$ analogue of a) certain derived category of $U(\mathfrak{g}^\vee_\k)$-modules with a fixed nontrivial Frobenius character and a category of $\Iu$-equivariant sheaves on $\sL G / \I_{\varnothing \subset J}^\semi$, mentioning that the definition of the latter category is problematic since ``the space $\sL G / \I_{\varnothing \subset J}^\semi$ behaves pretty badly.'' We expect that a more suitable replacement of this category is the $\infty$-category $\Shv \left( \I^\semi_{\varnothing \subset J} \backslash \tFl \right)$, and are planning to come back to this question in future work.

\subsection{Acknowledgements}

As noted above, this paper is part of a joint project with Pramod Achar and Gurbir Dhillon, and although they were not directly involved in the study of the specific aspects considered here their insights on these questions were very helpful. 
We thank Ivan Losev for explaining to us the idea of \cite{bez-losev} during the early stages of this work. 
The second author thanks Xinyu Li for his interest and helpful discussions. 

We acknowledge the use of Copilot and Claude for editing purposes and to check correctness of our proofs. No mathematical content was generated or suggested by AI.

\section{Affine Weyl group and affine flag varieties} 
\label{sec:Wext}

\subsection{Roots and weights} 
\label{ss:roots-weights}

Let $\F$ be an algebraically closed field, whose characteristic will be denoted $p$, and let $G$ be a connected reductive algebraic group over $\F$. 
We fix a Borel subgroup $B\subset G$ and a maximal torus $T\subset B$, and denote by $U\subset B$ the unipotent radical of $B$. 

Let $\X=X^*(T)$ be the character lattice, and $\Xd=X_*(T)$ be the cocharacter lattice. 
Let $\Phi \subset \X$ be the root system of $G$ with respect to $T$, let $\Phi^+\subset \Phi$ be the subset of positive roots determined by $B$ (i.e.~the set of nonzero $T$-weights in $\mathrm{Lie}(B)$), and let $\Phis \subset \Phi^+$ be the associated system of simple roots.
Let $\Xd^{++}\subset \Xd^+ \subset \Xd$ be the subsets of strictly dominant and dominant coweights respectively, and $\check{\Phi} \subset \Xd$ be the set of coroots. 
For $\alpha \in \Phi$, we will denote as usual by $\check{\alpha} \in \check{\Phi}$ the associated coroot. 
We set $2\rho=\sum_{\alpha\in \Phi^+} \alpha \in \X$.
For $\alpha\in \Phi$, we denote by $U_\alpha$ the associated one-parameter subgroup in $G$, and fix an isomorphism $u_\alpha:\Ga\xs U_\alpha$ such that $t u_\alpha(a) t^{-1} = u_\alpha(\alpha(t) a)$ for any $a \in \Ga$ and $t \in T$. 

Let $W$ be the Weyl group of $(G,T)$, endowed with the Coxeter group structure with set of Coxeter generators $S=\{s_\alpha : \alpha\in \Phis\}$, where $s_\beta$ is the reflection associated with $\beta\in \Phi$. 
We will denote by $w_0\in W$ the longest element. 
For every $w \in W$, we fix a lift $\dot{w}\in \mathrm{N}_G(T)$ of $w$. 

If $J$ is a subset of $S$, we set $J^\comp:=S\smallsetminus J$. 
For such $J$, we will denote by $P_J$ the standard parabolic subgroup of $G$ associated with $J$, with Levi decomposition $P_J=M_J\ltimes U_{J}$ where $M_J$ is the Levi subgroup of $P_J$ containing $T$ and $U_{J}$ is the unipotent radical of $P_J$. 
Let
\[
\Phis_J=\{\alpha\in \Phis \mid s_\alpha\in J\}, \quad \Phi^+_J=\Phi^+ \cap \Z\Phis_J, \quad \Phi_J=\Phi \cap \Z\Phis_J
\]
be the subsets of simple roots, positive roots, and roots of the group $M_J$, with respect to the maximal torus $T \subset M_J$ and the Borel subgroup $B_J:=B \cap M_J$. 
We set $2\rho_J=\sum_{\alpha\in \Phi^+_J} \alpha$. 
We will consider the following subsets of $\Xd$ attached to $J$:
\begin{align*}
\Xd_J^0 &:=\{\lambda\in \Xd \mid \forall \alpha\in \Phis_J, \, \langle \lambda, \alpha\rangle = 0\}, \\
\Xd^+_J &:=\{\lambda\in \Xd \mid \forall \alpha\in \Phis_J, \, \langle \lambda, \alpha\rangle \geq 0\}, \\
\Xd^{++}_J &:=\{\lambda\in \Xd \mid \forall \alpha\in \Phis_J, \, \langle \lambda, \alpha\rangle > 0\}. 
\end{align*}
(To compare with the notation of~\cite{DL23}, our $\Xd_J^0$ corresponds to their $\Lambda_{M,ab}$ where $M=M_J$, and the subset $\Xd_J^0 \cap \Xd^+=\Xd_J^0 \cap \Xd_{J^\comp}^+$ that will appear frequently below corresponds to their $\Lambda_{M,ab}^+$.) We will also denote by $W_J \subset W$ the parabolic subgroup generated by $J$, i.e.~the Weyl group of $(M_J,T)$.

\subsection{Affine Weyl group} 

\subsubsection{Affine roots and affine Weyl group} 

For the statements in this subsection, we refer to~\cite[\S\S 1--2]{macdonald}.
We consider the affine real roots 
\[
\Phi_\af=\Phi+\Z\hbar= \{\alpha+n\hbar : \alpha\in \Phi,n\in \Z\} \subset \X \oplus \Z \hbar,
\] 
together with the collection of positive roots
\[
\Phi^+_\af=(\Phi^++\Z_{\geq 0}\hbar) \sqcup (-\Phi^++\Z_{> 0}\hbar).
\] 
Following usual conventions, for $\alpha \in \Phi_\af$, we will write $\alpha >0$ for the condition that $\alpha \in \Phi^+_\af$, and $\alpha<0$ for the condition that $\alpha \in -\Phi^+_\af$.

We define the extended affine Weyl group $W_\ex$ and the affine Weyl group $W_\af$ as
\[
W_\ex=W\ltimes \Xd \quad \text{and} \quad W_\af= W\ltimes \Z\check{\Phi}. 
\]
For $\lambda\in \Xd$ we will denote by $t_\lambda := (1 \ltimes \lambda) \in W_\ex$ the translation element corresponding to $\lambda$. 
For $\alpha \in \Phi$ and $n \in \Z$,
we set $s_{\alpha+n\hbar} = t_{n\check{\alpha}} s_\alpha =s_\alpha t_{-n\check{\alpha}} \in W_\af$. 
The group $W_\af$ admits a structure of Coxeter group, with set of Coxeter generators $S_\af=S \cup \{s_{-\gamma+\hbar}\}_\gamma$ where $\gamma$ runs over the maximal roots in $\Phi^+$. 
Let $\ell: W_\af\rightarrow \Z_{\geq 0}$ be the length function associated with this Coxeter structure. 
It is well known that we have
\begin{equation}
\label{eq:length}
\ell(t_\lambda w) = \sum_{\substack{\alpha \in \Phi^+ \\ w^{-1}(\alpha) >0}} |\langle \lambda, \alpha \rangle | + \sum_{\substack{\alpha \in \Phi^+ \\ w^{-1}(\alpha) <0}} |1 + \langle \lambda,\alpha\rangle |
\end{equation}
for $w \in W$ and $\lambda \in \Z\check{\Phi}$. 
This formula makes sense more generally when $\lambda \in \Xd$, which allows to extend the function $\ell$ to $W_\ex$. 

Setting $\Omega := \{w \in W_\ex \mid \ell(w)=0\}$, it is a standard fact that $\Omega$ is a finitely generated abelian group, that conjugation by $\Omega$ on $W_\ex$ preserves $S_\af$, that multiplication induces a group isomorphism
\[
\Omega \ltimes W_\af \simto W_\ex,
\]
and finally that $\ell(w\omega) = \ell(w) = \ell(\omega w)$ for any $\omega \in \Omega$ and $w \in W_\ex$.

\begin{rmk}
\label{rmk:Coxeter-structures}
The Coxeter structure on $W_\af$ considered above is the same as the one used in~\cite{adr}. It is \emph{not} the same as the one used in~\cite{AR18}; more precisely these two structures are twists of each other by the automorphism of $W_\af$ which is the identity on $W$ and sends each $t_\lambda$ to $t_{-\lambda}$ ($\lambda \in \Z\check{\Phi}$), or alternatively are twists of each other by the automorphism given by conjugation by $w_0$.
\end{rmk}

Let $\leq$ be the Bruhat order on $W_\ex$, which is generated by 
\[
x< xt \quad \text{for any reflection $t\in W_\af$ such that $\ell(x)<\ell(xt)$},
\]
where we call $t\in W_\af$ a reflection if it is contained in
\[
\bigcup_{y\in W_\ex} yS_\af y^{-1} = \{s_\beta : \beta \in \Phi_\af\}.
\] 

The group $W_\ex$ acts on $\Phi_\af$ by the formula
\begin{equation}
\label{eq:act-aff}
(t_\lambda w)(\alpha+n\hbar)=w(\alpha)+(n+\langle \lambda,w(\alpha)\rangle)\hbar 
\end{equation}
for $\lambda\in \Xd$, $w\in W$, $\alpha \in \Phi$ and $n \in \Z$.
In terms of this action, we have
\[
\ell(w) = \#\{\alpha \in \Phi_\af^+ \mid w^{-1}(\alpha)<0\}.
\]
For any $x\in W_\ex$ and $\alpha\in \Phi^+_\af$ we also have
\begin{equation}
\label{eq:Bruhat}
x < xs_\alpha \quad \text{if and only if}\quad x(\alpha)\in \Phi^+_\af. 
\end{equation}


\subsubsection{Quarterinfinite length and order} 
\label{sss:quar-length}

Let $J\subset S$.
We define the \emph{quarterinfinite length}\footnote{The term ``quarterinfinite'' used here and below is taken from the discussion in~\cite[\S 8.11]{bez-losev}, which was one of the starting points for our investigations (see also~\S\ref{ss:into-RT}).} associated with $J$ on $W_\ex$ by 
\begin{multline*}
\ell^J(t_\lambda w) = \langle \lambda, 2\rho-2\rho_J \rangle 
+ \#\{ \alpha \in \Phi^+\smallsetminus \Phi^+_J \mid w^{-1}(\alpha) <0 \} \\ 
+ \sum_{\substack{\alpha \in \Phi_J^+ \\ w^{-1}(\alpha) >0}} | \langle \lambda, \alpha \rangle | 
+ \sum_{\substack{\alpha \in \Phi_J^+ \\ w^{-1}(\alpha) <0}} | 1+ \langle \lambda, \alpha \rangle |
\end{multline*}
for $\lambda \in \Xd$ and $w \in W$.
Note that for a fixed $x \in W_\ex$, for $N \gg 0$ we have
\[
\ell^J(x)=\ell(t_{\theta} x)-\ell(t_\theta)
\]
for any $\theta \in \Xd_J^0$ such that $\langle \theta,\alpha \rangle \geq N$ for all $\alpha \in \Phis_{J^\comp}$. 
If $J=S$ then $\ell^S$ is the length function $\ell$ from~\eqref{eq:length}; if $J=\varnothing$ then $\ell^{\varnothing}$ is the semiinfinite length considered e.g.~in~\cite[Remark~3.5]{adr}. 

\begin{rmk}
\phantomsection
\label{rmk:lJ}
\begin{enumerate}
\item
\label{it:lJ-translation}
It is clear from the definition that for $\theta \in \Xd_J^0$ and $x \in W_\ex$ we have $\ell^J(t_\theta x) = \ell^J(x) + \langle \theta, 2\rho \rangle$.
\item
\label{it:lJ-WJext}
Consider the extended affine Weyl group $W_{J,\ex}$ of the reductive group $M_J$ (with maximal torus $T$ and Borel subgroup $B_J$), and denote its length function by $\ell_J$. Then $W_{J,\ex}$ identifies with a subgroup of $W_\ex$, and any $w \in W_\ex$ can be written uniquely as $w=yx$ with $y \in W_{J,\ex}$ and $x \in W$ the unique minimal-length representative of its coset $W_J x$. Writing $y=t_\lambda z$ with $\lambda \in \Xd$ and $z \in W_J$, using the facts that $z^{-1}$ permutes $\Phi^+ \smallsetminus \Phi_J^+$ and that $x^{-1}(\Phi_J^+) \subset \Phi^+$ (which is classical, and follows from~\eqref{eq:Bruhat}) we see that $\ell^J(w) = \langle \lambda, 2\rho-2\rho_J \rangle + \ell(x) + \ell_J(y)$.
\end{enumerate}
\end{rmk}

We also consider the \emph{quarterinfinite order} $\preceq^J$ on $W_\ex$ defined as follows. Consider the subset $W_\ex^{+,J} \subset W_\ex$ consisting of elements $x \in W_\ex$ which satisfy $\ell(t_\mu x)=\ell(t_\mu) + \ell(x)$ for any $\mu \in \Xd_J^0 \cap \Xd^+_{J^\comp}$. In particular, this subset contains all the elements of the form $t_\lambda w$ with $w \in W$ and $\lambda \in \Xd$ which satisfies $\langle \lambda, \beta \rangle \geq 0$ for all $\beta \in \Phi^+ \smallsetminus \Phi_J^+$, and is stable under left multiplication by any $t_\theta$ with $\theta \in \Xd_J^0 \cap \Xd^+_{J^\comp}$. As in~\cite[\S 3.1]{adr} one checks that for $x,y \in W_\ex$ the following conditions are equivalent:
\begin{enumerate}
\item
there exists $\mu \in \Xd^0_J$ such that $t_\mu x$ and $t_\mu y$ belong to $W_\ex^{+,J}$ and $t_\mu x \leq t_\mu y$;
\item
for any $\mu \in \Xd^0_J$ such that $t_\mu x$ and $t_\mu y$ belong to $W_\ex^{+,J}$, we have $t_\mu x \leq t_\mu y$.
\end{enumerate}
When these conditions are satisfied, we write $x \preceq^J y$.

This order is ``$\Xd^0_J$-periodic'' in the sense that for $x,y \in W_\ex$ and $\lambda \in \Xd^0_J$ we have $x \preceq^J y$ if and only if $t_\lambda x \preceq^J t_\lambda y$.
When $J=S$ we have $W_\ex^{+,S} = W_\ex$, and $\preceq^S$ is the Bruhat order. When $J=\varnothing$, $W_\ex^{+,\varnothing}$ is the subset denoted $W_\ex^+$ in~\cite[\S 3.1]{adr}, and $\preceq^{\varnothing}$ is the ``semiinfinite order,'' which is denoted $\leq^\semi$ in~\cite{adr}.
The order $\preceq^{J}$ can also be characterized using the quarterinfinite length: it is the order generated by 
\[
x \preceq^{J} xt \quad \text{for any reflection $t\in W_\af$ such that}\ \ell^J(x)<\ell^J(xt). 
\]
In particular, $\ell^J$ is increasing with respect to $\preceq^J$ (and the standard order on $\Z$).

We set 
\[
{}^J\Phi^+_\af:= (\Phi_J^++\Z_{\geq 0}\hbar) \sqcup (-\Phi_J^++\Z_{> 0}\hbar) \sqcup ((\Phi^+\smallsetminus \Phi_J^+) +\Z \hbar). 
\]
Note that for $\alpha \in \Phi_\af$, we have $\alpha \in {}^J\Phi^+_\af$ if and only if there exists $\theta \in \Xd_J^0$ such that $t_{\theta+\theta'} (\alpha) \in \Phi_\af^+$ for any $\theta' \in \Xd_J^0 \cap \Xd_{J^\comp}^+$.
The following lemma is an analogue of~\eqref{eq:Bruhat}.

\begin{lem}
\label{lem:quar-order}
For any $x\in W_\ex$ and $\alpha\in \Phi^+_\af$, we have $x\prec^{J} xs_\alpha$ if and only if $x(\alpha)\in {}^J\Phi_\af^+$. 
\end{lem}

\begin{proof}
Assume that $x\prec^{J} xs_\alpha$, and choose $\theta \in \Xd_J^0$ such that $t_\theta x$ and $t_\theta x s_\alpha$ belong to $W_\ex^{+,J}$. Then for any $\theta' \in \Xd_J^0 \cap \Xd_{J^\comp}^+$ we have  $t_{\theta+\theta'} x \leq t_{\theta + \theta'} x s_\alpha$, so $t_{\theta+\theta'} x(\alpha) \in \Phi_\af^+$ by~\eqref{eq:Bruhat}. This forces $x(\alpha)\in {}^J\Phi_\af^+$.

Conversely, assume that $x(\alpha)\in {}^J\Phi_\af^+$, and choose $\theta \in \Xd^0_J$ such that $t_\theta x$ and $t_\theta x s_\alpha$ belong to $W_\ex^{+,J}$ and $t_{\theta} x(\alpha) \in \Phi_\af^+$. (Such a $\theta$ indeed exists.) Then we have $t_{\theta} x < t_{\theta} x s_\alpha$ by~\eqref{eq:Bruhat}, hence $x\prec^{J} xs_\alpha$.
\end{proof}

\subsubsection{Minimal coset representatives} 

For any $J\subset S_\af$, we denote by $W_J$ the subgroup of $W_\ex$ generated by $J$. 
Recall that $J$ is called \emph{finitary} if $W_J$ is finite; in this case
we will denote by $w_J$ the longest element in $W_J$. (Note in particular that any subset of $S$ is finitary, and that we have $w_S=w_0$.)

For $J\subset S$ and $K\subset S_\af$ finitary we set 
\begin{gather*}
{}^JW_\ex:=\{x\in W_\ex \mid 
\forall w\in W_J, \, \ell(x)\leq \ell(wx)\}, 
\\
W_\ex^{(K,J)}:=\{x\in W_\ex \mid 
\forall w\in W_K, \, \ell^{J}(x)\leq \ell^{J}(xw)\}.
\end{gather*}
In particular, ${}^J W_\ex$ is the set of minimal-length representatives for the cosets in $W_J \backslash W_\ex$ (with respect to $\ell$), and for $x \in W_\ex$ we have $x \in {}^J W_\ex$ if and only if $sx>x$ for any $s \in J$. When $J=S$, we will write $W_\ex^K=W_\ex^{(K,S)}$; this set is the set of minimal-length representatives for the cosets in $W_\ex / W_K$. (We will use similar conventions for the variants of this notation introduced below.)


\begin{lem}
\label{lem:min-length}
Let $J\subset S$ and $K\subset S_\af$ be a finitary subset. Let $x\in W_\ex$.
\begin{enumerate}
\item 
\label{it:min-length-0}
If $\theta \in \Xd_J^0$, we have
\[
x \in {}^J W_\ex \Leftrightarrow t_\theta x \in {}^J W_\ex, \qquad
x \in W_\ex^{(K,J)} \Leftrightarrow t_\theta x \in W_\ex^{(K,J)}.
\]
\item 
\label{it:min-length-1}
We have $x\in {}^JW_\ex$ if and only if $x$ is minimal in $W_Jx$ for the quarterinfinite length $\ell^J$. 
In this case, for any $w\in W_J$ we have $\ell^{J}(wx)=\ell^{J}(x)+\ell(w)$. 
\item 
\label{it:min-length-2}
If $x\in W_\ex^{(K,J)}$, for any $w\in W_K$ we have $\ell^J(xw)=\ell^J(x)+\ell(w)$. 
\end{enumerate}
\end{lem}

\begin{proof}
\eqref{it:min-length-0}
Let $\theta \in \Xd_J^0$.
By~\eqref{eq:Bruhat}, we have $x\in {}^JW_\ex$ if and only if $x^{-1}(\alpha)\in \Phi^+_\af$ for any $\alpha\in \Phis_J$. 
Since $t_{-\theta}$ fixes $\Phis_J$ (under the action \eqref{eq:act-aff}), 
this is equivalent to the condition that $x^{-1}t_{-\theta}(\alpha)\in \Phi^+_\af$ for any $\alpha\in \Phis_J$, i.e.~$t_\theta x \in {}^J W_\ex$.

Now, choose $\theta' \in \Xd_J^0$ such that for any $w \in W_K$ we have
\[
\ell^J(xw) = \ell(t_{\theta+\theta'} xw) - \ell(t_{\theta+\theta'}), \quad \ell^J( t_{\theta}xw) = \ell(t_{\theta'} t_{\theta}xw) - \ell (t_{\theta'}).
\]
We have $x \in W_\ex^{(K,J)}$ if and only if $\ell(t_{\theta+\theta'} xw) \geq \ell(t_{\theta+\theta'} x)$ for any $w \in W_K$, which is equivalent to $t_\theta x \in W_\ex^{(K,J)}$.

\eqref{it:min-length-1}
Let $x \in W_\ex$, and fix $\theta\in \Xd_J^0$ 
such that for any $w \in W_J$ we have
\[
\ell^J(wx) = \ell(t_\theta wx) - \ell(t_\theta).
\]
If $x \in {}^J W_\ex$ then $t_\theta x \in {}^J W_\ex$ by~\eqref{it:min-length-0}; since each $w$ as above commutes with $t_\theta$ we deduce that $\ell(t_\theta wx)=\ell(w) + \ell(t_\theta x)$, hence $\ell^J(wx)=\ell^J(x)+\ell(w)$. In particular, $x$ is minimal in $W_J x$ for $\ell^J$. Conversely, if $\ell^J(wx) \geq \ell^J(x)$ for any $w \in W_J$, then $t_\theta x \in {}^J W_\ex$, hence $x \in {}^J W_\ex$ by~\eqref{it:min-length-0}.

\eqref{it:min-length-2}
Choose $\theta \in \Xd_J^0$ such that 
for any $w \in W_K$ we have
\[
\ell^J(xw) = \ell(t_\theta xw)-\ell(t_\theta).
\]
We have $t_\theta x\in W_\ex^{K}$, and then $\ell(t_\theta xw) = \ell(t_\theta x)+\ell(w)$ for any $w \in W_K$, which implies the desired formula.
\end{proof}

It follows in particular from Lemma~\ref{lem:min-length}\eqref{it:min-length-2} that $W_\ex^{(K,J)}$ is a subset of representatives for the cosets in $W_\ex/W_K$. For $L \subset S$ we set 
\[
{}^{L}W_\ex^{(K,J)}:={}^{L}W_\ex\cap W_\ex^{(K,{J})}.
\]
When $J=S$, ${}^{L}W_\ex^{K}$ is the set of minimal-length representatives for the cosets in $W_{L} \backslash W_\ex / W_K$.
For any subset $D \subset S_\af$ and $x\in W_\ex$, we set ${}^x D=xDx^{-1}$, $D^x=x^{-1}Dx$. For $D_2 \subset D_1$ finitary, we also set $W_{D_1}^{D_2} := W_{D_1} \cap W_\ex^{D_2}$, ${}^{D_2} W_{D_1} := W_{D_1} \cap {}^{D_2} W_\ex$.

\begin{lem}
\label{lem:bi-min}
Let $L \subset J\subset S$ and $K\subset S_\af$ be finitary. 
\begin{enumerate}
\item 
\label{it:bi-min-1}
The composition
\[
{}^{L} W_\ex^{(K,J)} \hookrightarrow W_\ex \twoheadrightarrow W_{L}\backslash W_\ex /W_K
\]
is a bijection. 
\item 
\label{it:bi-min-2}
Let $x\in {}^{L}W_\ex^{(K,J)}$. 
We have $W_{L}\cap xW_Kx^{-1} = W_{L \cap ({}^xK)}$. Moreover, for any $w \in W_{L}^{L\cap ({}^xK)}$ we have $wx\in W_\ex^{(K,J)}$, and all the representatives in $W_\ex^{(K,J)}$ of the cosets of the form $wxW_K$ with $w \in W_L$ are of this form.
\item 
\label{it:bi-min-3}
Let $x\in {}^{L} W_\ex^{(K,J)}$. 
We have $x^{-1}W_{L} x\cap W_K = W_{(L^x)\cap K}$. Moreover, for any $w \in {}^{(L^x)\cap K}W_K$ we have $xw\in {}^{L}W_\ex$, and the minimal representatives in all cosets $W_L xw$ with $w \in W_K$ are of this form.
\end{enumerate}
\end{lem}

\begin{proof}
All the statements reduce to well-known facts about the subsets ${}^{L} W_\ex^K$, using the facts that for $\theta \in \Xd_J^0$, $t_\theta$ commutes with $W_{L}$ and 
left multiplication by $t_\theta$ stabilizes 
${}^{L}W_\ex$ and
$W_\ex^{(K,J)}$ (see Lemma~\ref{lem:min-length}\eqref{it:min-length-0}). More specifically, one proceeds as follows.

For~\eqref{it:bi-min-1}, we start with injectivity.
Let $x,y\in {}^{L}W_\ex^{(K,J)}$ be two elements which have the same image, i.e.~such that $W_{L} x W_K = W_{L} y W_K$. 
As in the proof of Lemma~\ref{lem:min-length}\eqref{it:min-length-2}, there exists $\theta \in \Xd_J^0$ 
such that $t_\theta x$ and $t_\theta y$ belong to $W_\ex^{K}$, and then these elements belong to ${}^{L} W_\ex^{K}$. Since $W_{L} t_\theta x W_K = t_\theta W_{L} x W_K = t_\theta W_{L} y W_K = W_{L} t_\theta y W_K$, we have $t_\theta x=t_\theta y$, and then $x=y$.

To show surjectivity, we consider a double coset $X \in W_{L} \backslash W_\ex /W_K$.
Since $X$ is finite, there exists $\theta \in \Xd_J^0$ such that for all $w \in X$ we have
\[
\ell^J( w) = \ell(t_\theta w) - \ell(t_\theta).
\]
If $z$ is the unique element in $t_\theta X$ (which is a coset in $W_{L}\backslash W_\ex /W_K$) which belongs to ${}^{L} W_\ex^{K}$, then $t_{-\theta} z$ is an element of $X$ which belongs to ${}^{L} W_\ex^{(K,J)}$.


For~\eqref{it:bi-min-2} and \eqref{it:bi-min-3}, assume first that $J=S$. In this case,~\eqref{it:bi-min-3} follows from~\cite[Lemma~2.1.9 and Theorem~2.1.12]{geck-pfeiffer} and~\eqref{it:bi-min-2} is obtained by applying these statements to $x^{-1}$, with the roles of $K$ and $L$ switched. (These statements are proved for Coxeter groups in~\cite{geck-pfeiffer}, but they follow for an ``almost Coxeter group'' such as $W_\ex$.)
The general case follows,
using multiplication by an appropriate $t_\theta$ with $\theta \in \Xd_J^0$.
\end{proof}

For $J$, $K$ as above, we now set
\[
W_\ex^{(K,J),\max} := W_\ex^{(K,J)} \cdot w_K,
\]
so that by Lemma~\ref{lem:min-length}\eqref{it:min-length-2}
we have
\[
W_\ex^{(K,J),\max} = \{x \in W_\ex \mid \forall w\in W_K, \, \ell^{J}(x)\geq \ell^{J}(xw) \}.
\]
For $L\subset S$ we set 
\[
{}^{L}W_\ex^{(K,J),\max}:={}^{L}W_\ex\cap W_\ex^{(K,{J}),\max}. 
\]
Here, one can check (using e.g.~Lemma~\ref{lem:quar-min-max} below) that the natural map
\[
{}^{L}W_\ex^{(K,J),\max} \to W_{L} \backslash W_\ex / W_K
\]
is injective, but it is \emph{not} surjective in general.

\begin{lem}
\label{lem:quar-min-max}
Let $L \subset J\subset S$ and $K\subset S_\af$ be finitary.
For 
any $x\in W_\ex$, the following properties are equivalent:
\begin{enumerate}
\item 
$x\in {}^{L} W_\ex^{(K,J),\max}$; 
\item 
$x\in {}^{L} W_\ex$, and for any $w\in W_{L}$ we have $wx\in W_\ex^{(K,J),\max}$; 
\item 
$x\in W_\ex^{(K,J),\max}$, and for any $u\in W_K$ we have $xu\in {}^{L} W_\ex$; 
\item 
for any $w\in W_{L}$ and $u\in W_K$, we have $\ell^J(wxu)=\ell(w)+\ell^J(x)-\ell(u)$. 
\end{enumerate}
\end{lem}

\begin{proof}
As for Lemma~\ref{lem:bi-min}, the assertion follows from the analogous statement for the usual length function, which follows (using right multiplication by $w_K$) from~\cite[Lemma~2.4]{AR23}.
\end{proof}


In the following lemma we consider the special case $K=S$ and $L=J$.

\begin{lem}
\label{lem:JSJw}
For any $J \subset S$ we have 
\[
{}^JW_\ex^{(S,J),\max}=\{w_Jt_\lambda w_0 : \lambda\in \Xd_J^{++}\}. 
\]
\end{lem}

\begin{proof}
This statement will be deduced from~\cite[Lemma~10.2]{AR18}. (This reference has some technical assumptions on the reductive group under consideration, but the general case can be easily reduced to that of~\cite{AR18}.) The latter statement describes the elements $x \in W_\ex$ which are minimal in $Wx$ and maximal in $xW_J$, with respect to the ``other'' Coxeter structure mentioned in Remark~\ref{rmk:Coxeter-structures}. Noting that the assignment $x \mapsto w_J x^{-1} w_0$ induces a bijection between this set and 
the set of elements $y$ which are maximal in $yW$ and minimal in $W_J y$, and taking into account the difference of Coxeter group structures, we obtain that with our present conventions the set of elements $y$ which are maximal in $yW$ and minimal in $W_J y$ consists of the elements $(w_J t_\nu v_\nu w_0 : \nu \in \Xd_J^{++})$ where $v_\nu \in W$ is the element of minimal length such that $v_\nu(\nu) \in \Xd^+$.

Let now $x \in {}^JW_\ex^{(S,J),\max}$, and choose $\theta \in \Xd_J^0$ such that
\[
\ell^J(w) = \ell(t_\theta w)-\ell(t_\theta)
\]
for any $w \in x W$. Writing $x=t_\mu y$ with $\mu \in \Xd$ and $y \in W$, we can (and shall) also assume that 
$\langle \theta+\mu, w_J(\alpha) \rangle \geq 0$ for any $\alpha \in \Phis_{J^\comp}$.
Then $t_\theta x$ is minimal in $W_J t_\theta x$ by Lemma~\ref{lem:min-length}\eqref{it:min-length-0} and maximal in $t_\theta x W$ by our choice of $\theta$, hence $t_\theta x = w_J t_\nu v_\nu w_0$ for some $\nu \in \Xd_J^{++}$. We then have
\[
t_\theta x = t_{\theta+\mu} y = t_{w_J(\nu)} w_J v_\nu w_0.
\]
It follows that $\theta+\mu = w_J(\nu)$, hence $\langle \nu, \alpha \rangle \geq 0$ for $\alpha \in \Phis_{J^\comp}$, so that $\nu \in \Xd^+$ and therefore $v_\nu=1$. We deduce that
\[
x = w_J t_{\nu-\theta} w_0,
\]
hence $x$ has the required form with $\lambda=\nu-\theta \in \Xd_J^{++}$.

Conversely, let $\lambda \in \Xd_J^{++}$, and set $x = w_J t_\lambda w_0$. Choose $\theta \in \Xd_J^0$ such that 
\[
\ell^J(w) = \ell(t_\theta w)-\ell(t_\theta)
\]
for any $w \in x W$, and moreover
$\lambda+\theta \in \Xd^+$. Then the element $t_\theta x = w_J t_{\lambda + \theta} w_0$ is minimal in $W_J t_\theta x$ and maximal in $t_\theta x W$. It follows that $x \in W_\ex^{(S,J),\max}$ by our choice of $\theta$, and that $x \in {}^J W_\ex$ by Lemma~\ref{lem:min-length}\eqref{it:min-length-0}.
Hence $x \in {}^JW_\ex^{(S,J),\max}$, as desired.
\end{proof}

\begin{rmk}
\label{rmk:JSJw-length}
A direct computation shows that if $\lambda\in \Xd_J^{++}$ we have $\ell^J(w_Jt_\lambda w_0)= \langle \lambda, 2\rho \rangle + \ell(w_0) - \ell(w_J)$.
\end{rmk}

\subsection{Affine flag varieties} 

\subsubsection{Loop groups, parahoric subgroups, and affine flag varieties} 
\label{sss:def-Fl}

For an affine scheme $X$ over $\F$, we denote by 
\[
\sL{X},\quad \text{resp.}\ \sL^+X, \quad \text{resp.}\ \sL^-X,
\]
the ind-affine ind-scheme, resp. affine scheme, resp. ind-affine ind-scheme, represen\-ting the functor sending an $\F$-algebra $R$ to $X(R( \hspace{-1pt} (z) \hspace{-1pt} ))$, resp.~$X(R[ \hspace{-1pt} [z] \hspace{-1pt} ])$, resp.~$X(R[z^{-1}])$, where $z$ is an indeterminate. 
There is a natural morphism
\[
\mathrm{ev}^X_+:\sL^+X\rightarrow X, \quad \text{resp.} \quad \mathrm{ev}^X_-:\sL^-X\rightarrow X,
\]
induced by evaluating $z$, resp. $z^{-1}$, to zero. 
When $X$ is an affine group scheme, $\sL{X}$, resp. $\sL^+X$, resp. $\sL^-X$, admits a structure of group ind-scheme, resp. group scheme, resp. group ind-scheme. 
In this case we set
\[
\sL^{++}X=\ker(\mathrm{ev}^X_+), \quad \sL^{--}X=\ker(\mathrm{ev}^X_-).
\] 
We call $\sL X$ the loop group of $X$, and $\sL^+X$ the arc group of $X$. 

We will consider in particular these constructions for our groups $T$, $U$ and $G$.
For $\lambda\in \Xd$, we denote by $z^\lambda$ the $\F$-point of $\sL{G}$ given by the image of the $\F$-point $z\in \F(\!(z)\!)^\times = \sL\Gm(\F)$ under the composition $\sL\Gm\xrightarrow{\sL\lambda} \sL{T} \rightarrow \sL{G}$. 

We define the Iwahori subgroup $\I$ of $\sL^+ G$, resp.~its pro-unipotent radical $\Iu$, as
\[
\I=(\mathrm{ev}^G_+)^{-1}(B), \quad \text{resp.} \quad \Iu=(\mathrm{ev}_+^G)^{-1}(U).
\] 
For any finitary subset $J \subset S_\af$, we have an associated parahoric subgroup $\I_J \subset \sL G$ containing $\I$, see e.g.~\cite[\S 3.4]{AR23}. If we set $\Iu^J := \dot{w}_J \Iu (\dot{w}_J)^{-1}$, the intersection $\Iu \cap \Iu^J$ is the pro-unipotent radical of $\I_J$, and we have a canonical Levi factor $M_J \subset \I_J$ such that multiplication induces an isomorphism
\[
M_J \ltimes (\Iu \cap \Iu^J) \simto \I_J.
\]
If $J \subset S$ 
we have
\[
\I_{J}=(\mathrm{ev}^G_+)^{-1}(P_J), 
\]
and $M_J$ is the subgroup of $G$ defined in~\S\ref{ss:roots-weights} (now seen as a subgroup of $\sL G$).

If $L \subset J \subset S$, then $J$ is the set of simple reflections of the Weyl group $W_J$ of $(M_J,T)$, and we will denote by
$\I_{L \subset J}\subset \sL^+M_J$ the parahoric subgroup of $\sL^+ M_J$ associated with the subset $L \subset J$. (Here, the fixed Borel subgroup of $M_J$ used for this construction is $B_J$.)
We also set
\[
\I^\semi_{L \subset J} := \I_{L \subset J} \ltimes \sL U_J.
\]
The following special cases are particularly relevant.
\begin{itemize}
\item
When $L=J$, we will simplify the notation and set $\I_{J}^{\semi} := \I_{J \subset J}^{\semi}$. This group ind-scheme
is considered in particular in~\cite{DL23}.
\item
When $L=\varnothing$, we will set $\I_{J}^{\quar}:= \I_{\varnothing \subset J}^\semi$, and call this group ind-scheme the quarterinfinite Iwahori subgroup associated with $J$. This object appears in the discussion in~\cite[\S 8.11]{bez-losev}.
\item
When $J=\varnothing$ (so that $L=\varnothing$) we have $\I_{\varnothing}^{\quar}=\I_{\varnothing}^{\semi}$, and this group coincides with the \textit{semiinfinite Iwahori subgroup} $\I^\semi$ considered in~\cite{adr} (and earlier, using different notation, in various references including~\cite{Gaits18}).
\item
When $J=S$, we have $\I_{L \subset S}^\semi = \I_L$; in particular $\I_{S}^\quar = \I$ and $\I_S^\semi = \sL^+ G$.
\end{itemize}


We define the affine flag variety and the enhanced affine flag variety  by 
\[
\Fl=\sL G/\I, \quad \tFl=\sL G/\Iu.
\]
For a finitary subset $K\subset S_\af$, we also have the partial affine flag variety
\[
\Fl_K=\sL G/\I_K. 
\]
We will denote by $q:\tFl\rightarrow \Fl$ and $p_K: \Fl\rightarrow \Fl_K$ the natural projections. In case $K=S$, $\Fl_S$ is the familiar affine Grassmannian $\Gr = \sL G / \sL^+ G$.

\subsubsection{Orbits} 

For $x=t_\lambda w$ ($\lambda\in \Xd$, $w\in W$), we set $\dot{x}=z^\lambda \dot{w}\in \sL{G}(\F)$. 
Let $L \subset J \subset S$, and let $K \subset S_\af$ be a finitary subset.
We have decompositions of $\Fl_K$ into $\I_{L \subset J}^\semi$-orbits
\begin{equation}
\label{eq:decomp}
\Fl_K= \bigsqcup_{x\in {}^L W_\ex^{(K,J)}} \rS^{L \subset J}_{K,x},
\end{equation}
where 
\[
\rS^{L \subset J}_{K,x}:= \I_{L \subset J}^{\semi} \dot{x} \I_K/\I_K.
\]
When $K=S$ and $L=J=\varnothing$, these orbits are the familiar ``semiinfinite orbits'' which play a crucial role in many constructions related to the \emph{geometric Satake equivalence}, see~\cite{mv}. When $K=\varnothing$ and $L=J=\varnothing$, these orbits play a role in some constructions in~\cite{ab} (see also~\cite[\S 4.5]{ARbook}).
When $J=S$ (so that $\I_{L \subset S}^\semi = \I_L$),
$\rS^{L \subset S}_{K,x}$ is a (smooth) scheme of finite type over $\F$, and we have
\begin{equation}
\label{eqn:dim-par-orbit}
\dim(\rS^{L \subset S}_{K,x}) = \ell(x) + \ell(w_L w_{L \cap ({}^xK)}) = \ell(x) + \ell(w_L) -\ell(w_{L \cap ({}^xK)}).
\end{equation}
(This formula follows e.g.~from standard facts about the Bruhat decomposition, together with Lemma~\ref{lem:bi-min}\eqref{it:bi-min-2}.)

We denote by $\bs^{L \subset J}_{K,x}$ the locally closed embedding of $\rS^{L \subset J}_{K,x}$ in $\Fl_K$. 
We will omit $L$, $J$ or $K$, in the notations above if it is the empty set. (In particular, $\rS^J_{K,x}$ will mean $\rS^{\varnothing \subset J}_{K,x}$, the $\I_J^\quar$-orbit of $\dot{x} \I_K/\I_K$.)
We also have a decomposition 
\begin{equation*}
\tFl = \bigsqcup_{x\in {}^L W_\ex} q^{-1}(\rS^{L \subset J}_x)
\end{equation*}
into $\I_{L \subset J}^{\semi}$-orbits.

For $\alpha\in \Phi$ and $n\in \Z$, we denote by $U_{\alpha+n\hbar}$ the one-parameter subgroup in $\sL G$ given by the image of $a\in\Ga \mapsto u_\alpha(az^n)$. The following lemma, which can be checked by standard computations, describes the $\I^\quar_J$-orbits on $\Fl$.
Here and below, a product $\prod_{n=0}^k$ with $k<0$ or a product $\prod_{n=1}^k$ with $k \leq 0$ should be interpreted as $\Spec(\F)$.

\begin{lem}
Let $x=t_\lambda w$ ($\lambda\in \Xd$, $w\in W$).
Multiplication on $\dot{x}\I/\I$ induces an isomorphism 
\begin{equation}
\label{eq:quar-orbit}
\begin{aligned}
&\left(\prod_{\substack{\alpha \in \Phi_J^+ \\ w^{-1}(\alpha)>0}} \prod_{n=0}^{\langle\lambda,\alpha\rangle-1}U_{\alpha+n\hbar} \right) \times 
\left(\prod_{\substack{\alpha \in \Phi_J^+ \\ w^{-1}(\alpha)<0}} \prod_{n=0}^{\langle\lambda,\alpha\rangle}U_{\alpha+n\hbar} \right) \times \\ 
&\qquad \left(\prod_{\substack{\alpha \in -\Phi_J^+ \\ w^{-1}(\alpha)>0}} \prod_{n=1}^{\langle\lambda,\alpha\rangle-1}U_{\alpha+n\hbar} \right) \times 
\left(\prod_{\substack{\alpha \in -\Phi_J^+ \\ w^{-1}(\alpha)<0}} \prod_{n=1}^{\langle\lambda,\alpha\rangle}U_{\alpha+n\hbar} \right) \times \\ 
&\qquad\qquad \left(\prod_{\substack{\alpha \in \Phi^+\smallsetminus \Phi_J^+ \\ w^{-1}(\alpha)>0}} \Ad_{z^{\lambda}}(\sL^{--} U_{\alpha}) \right) \times 
\left(\prod_{\substack{\alpha \in \Phi^+\smallsetminus \Phi_J^+ \\ w^{-1}(\alpha)<0}} \Ad_{z^{\lambda}}(\sL^{-} U_{\alpha}) \right) 
\xs \rS^{J}_{x}. 
\end{aligned}
\end{equation}
\end{lem}

Note that the affine root subgroups appearing in the left-hand side of~\eqref{eq:quar-orbit} are parametrized by
\[
\{ \alpha \in {}^J \Phi_\af^+ \mid x^{-1}(\alpha)<0\}.
\]
In the special case $J=S$, we obtain an isomorphism
\begin{equation*}
\begin{aligned}
&\left(\prod_{\substack{\alpha \in \Phi^+ \\ w^{-1}(\alpha)>0}} \prod_{n=0}^{\langle\lambda,\alpha\rangle-1}U_{\alpha+n\hbar} \right) \times 
\left(\prod_{\substack{\alpha \in \Phi^+ \\ w^{-1}(\alpha)<0}} \prod_{n=0}^{\langle\lambda,\alpha\rangle}U_{\alpha+n\hbar} \right) \times \\
&\qquad\qquad \left(\prod_{\substack{\alpha \in -\Phi^+ \\ w^{-1}(\alpha)>0}} \prod_{n=1}^{\langle\lambda,\alpha\rangle-1}U_{\alpha+n\hbar} \right) \times 
\left(\prod_{\substack{\alpha \in -\Phi^+ \\ w^{-1}(\alpha)<0}} \prod_{n=1}^{\langle\lambda,\alpha\rangle}U_{\alpha+n\hbar} \right)
\xs \I \dot{x} \I/\I, 
\end{aligned}
\end{equation*}
which generalizes~\cite[Lemma~3.2]{adr}.


Given locally closed sub-ind-schemes $X,Y \subset \Fl$, we will write $X \times^\I Y$ for the twisted product $\widetilde{X} \times^\I Y$ where $\widetilde{X}$ is the preimage of $X$ in $\sL G$.

\begin{lem}
\label{lem:min-push}
Let $J \subset S$, and let $K \subset S_\af$ be finitary.
Let $x\in W_\ex^{(K,J)}$.

\begin{enumerate}
\item
\label{it:min-push-1}
We have 
\[
(p_K)^{-1}(\rS^J_{K,x})=\bigsqcup_{{w\in W_K}} \rS^J_{xw},
\]
and $p_K$ restricts to an isomorphism $\rS^J_{x}\xs \rS^J_{K,x}$.
\item
\label{it:min-push-2}
For any $w\in W_K$, multiplication induces an isomorphism
\[
\rS^J_{x}\times^{\I}\I \dot{w}\I/\I \xs 
\rS^J_{xw}. 
\]
\end{enumerate}
\end{lem}

\begin{proof}
This follows from~\eqref{eq:quar-orbit} 
and the techniques from Lemma~\ref{lem:min-length}. For instance, \eqref{it:min-push-2} follows from the equality
\begin{multline*}
\{ \alpha \in {}^J \Phi_\af^+ \mid (xw)^{-1}(\alpha)<0\} = \\
\{ \alpha \in {}^J \Phi_\af^+ \mid x^{-1}(\alpha)<0\} \sqcup x(\{\alpha \in \Phi_\af^+ \mid w^{-1}(\alpha)<0\}),
\end{multline*}
which can be deduced from the corresponding fact for the standard Bruhat order (corresponding to $J=S$) using the fact that if $\alpha \in {}^J \Phi_\af^+$, we have $t_\theta(\alpha) \in \Phi_\af^+$ if $\theta \in \Xd^0_J$ satisfies $\langle \theta, \beta \rangle \gg 0$ for any $\beta \in \Phis_{J^\comp}$.
\end{proof}

\begin{lem}
\label{lem:bi-min-av}
Let $L \subset J \subset S$, and let $K \subset S_\af$ be finitary.
Let $x\in {}^L W_\ex^{(K,J)}$. 
\begin{enumerate}
\item 
\label{it:bi-min-av-1}
For any $w\in W_L^{L\cap ({}^xK)}$, we have canonical isomorphisms 
\[
B\dot{w}B\times^{B} (\rS_{K,x}^{J}) \xs \I_J^{\quar}\dot{w}\I_J^{\quar}\times^{\I_J^{\quar}} (\rS_{K,x}^{J})
 \xs \rS_{K,wx}^J,
\]
where the second map is induced by multiplication. 
\item 
\label{it:bi-min-av-2}
The morphism
\[
\I^\semi_{L \subset J} \times^{\I^\quar_J} \rS_{K,x}^J  \rightarrow \rS_{K,x}^{L \subset J}
\]
induced by multiplication
is a locally trivial fibration with fiber isomorphic to the flag variety of $M_{L \cap ({}^xK)}$. 
\item 
\label{it:bi-min-av-3}
We have 
\[
\rS_{K,x}^{L \subset J}
=\bigsqcup_{w\in W_L^{L\cap ({}^xK)}} \rS^J_{K,wx}. 
\]
\end{enumerate}
\end{lem}

\begin{proof}
\eqref{it:bi-min-av-1}
By Lemma~\ref{lem:bi-min}\eqref{it:bi-min-2} we have $wx \in W_\ex^{(K,J)}$. Then by Lemma~\ref{lem:min-push}\eqref{it:min-push-1} we can assume that $K=\varnothing$. In this case, 
the proof is similar to that of Lemma~\ref{lem:min-push}\eqref{it:min-push-2}, the corresponding equality being
\begin{multline*}
\{\alpha \in {}^J \Phi_\af^+ \mid (wx)^{-1} (\alpha) <0\} = \\
\{\alpha \in \Phi_\af^+ \mid w^{-1}(\alpha)<0\} \sqcup w(\{\alpha \in {}^J \Phi_\af^+ \mid x^{-1}(\alpha)<0\}).
\end{multline*}

\eqref{it:bi-min-av-2}
Since $\rS_{K,x}^{L \subset J}$ is an $\I^\semi_{L \subset J}$-orbit,
it suffices to analyze the fiber of this map over the image of $\dot{x}$.
This fiber is isomorphic to 
\[
\I^\semi_{L \subset J} \cap \Ad_{\dot{x}}(\I_K)/\I^\quar_J\cap \Ad_{\dot{x}}(\I_K)
\simeq \bigl( M_L \cap \Ad_{\dot{x}}(\I_K) \bigr) / \bigl( B_L \cap \Ad_{\dot{x}}(\I_K) \bigr). 
\]
Here, since $x \in {}^L W_\ex$, by~\eqref{eq:Bruhat} we have $x^{-1}(\Phi_L^+) \subset \Phi_\af^+$, i.e.~$\Phi_L^+ \subset x(\Phi_\af^+)$, hence $B_L \subset \Ad_{\dot{x}}(\I_K)$. On the other hand we have $W_L \cap x W_K x^{-1} = W_{L \cap ({}^xK)}$ by Lem\-ma~\ref{lem:bi-min}\eqref{it:bi-min-2}, hence $(M_L \cap \Ad_{\dot{x}}(\I_K)) / (B_L \cap \Ad_{\dot{x}}(\I_K))$ identifies with the flag variety of the Levi factor $M_{L\cap ({}^xK)}$ of $M_L$.


\eqref{it:bi-min-av-3}
It is clear that $\rS_{K,x}^{L \subset J}$ is the union of the $\I^\quar_J$-orbits associated with the representatives in $W_\ex^{(K,J)}$ of the cosets $wxW_K$ for $w \in W_L$. The claim therefore follows from
Lemma~\ref{lem:bi-min}\eqref{it:bi-min-2}.
%
\end{proof}

\begin{lem}
\label{lem:bi-min-push}
Let $L \subset J \subset S$, and let $K \subset S_\af$ be finitary.
Let $x\in {}^L W_\ex^{(K,J)}$. 
\begin{enumerate}
\item 
\label{it:bi-min-push-1}
For any $w\in {}^{(L^x)\cap K}W_K$, multiplication induces an isomorphism
\[
\rS^{L \subset J}_{x}\times^{\I}\I \dot{w}\I/\I \xs \rS^{L \subset J}_{xw}.
\]
\item
\label{it:bi-min-push-2}
The morphism $p_K: \rS^{L \subset J}_{x}\rightarrow \rS^{L \subset J}_{K,x}$ is a locally trivial fibration with fiber isomorphic to the flag variety of $M_{(L^x) \cap K}$. 
\item
\label{it:bi-min-push-3}
We have 
\[
(p_K)^{-1}(\rS^{L \subset J}_{K,x})=\bigsqcup_{w\in {}^{(L^x) \cap K}W_K} \rS^{L \subset J}_{xw}. 
\] 
\end{enumerate}
\end{lem}

\begin{proof}
The proof is similar to that of
Lemma~\ref{lem:bi-min-av}, using Lemma~\ref{lem:bi-min}\eqref{it:bi-min-3}.
\end{proof}

\section{Parabolic semiinfinite sheaves} 
\label{sec:quar-semi-sheaves}

\subsection{Definition} 
\label{ss:qs-sheaves}

We will work in the sheaf-theoretic context explained in \cite[\S\S 2.1--2.2]{adr}, in the case of \'etale sheaves.
So let $\ell$ be a prime number invertible in $\F$, and let $\k$ be either a finite extension or an algebraic closure of $\F_\ell$ or $\Q_\ell$.
For a prestack $X$ over $\F$, we have a $\k$-linear presentable stable $\infty$-category $\Shv(X)$ of $\k$-sheaves on $X$. The functorialities that these $\infty$-categories satisfy are briefly discussed in~\cite[\S 2.2]{adr} (following in particular some constructions from~\cite{zhu}); in particular, for any map $f : X \to Y$ of $\F$-prestacks we have a $!$-pullback functor $f^! : \Shv(Y) \to \Shv(X)$, which sends the dualizing complex $\omega_Y$ of $Y$ to the dualizing complex $\omega_X$ of $X$. These functors, as well as all the functors between such $\infty$-categories that we will consider below, are continuous, i.e.~colimit-preserving.

\begin{rmk}
It is likely that all the results below have analogues for sheaves in the analytic topology. However, since a detailed treatment of the construction of such $\infty$-categories is not yet available in the literature, we will stick to \'etale sheaves. (For the statements involving Whittaker conditions, their counterparts should involve the ``Kirillov model.'')
\end{rmk}


Fix subsets $L \subset J \subset S$, and
let $X$ be either $\Fl_K$ (for some finitary subset $K \subset S_\af$) or $\tFl$. 
Below we will study the $\infty$-category 
\[
\Shv \left( \I^\semi_{L \subset J} \backslash X \right).
\]
These $\infty$-categories (for the various choices of $L$, $J$ and $X$) admit natural descriptions in terms of limits/colimits of $\infty$-cate\-gories of sheaves on algebraic stacks of finite type, which can be obtained by adapting (in the obvious way) the discussions in~\cite[\S\S 2.5--2.6]{adr} or~\cite[\S 1.2]{Gaits18}.


If $L_1 \subset L_2 \subset J$,
since $\I_{L_1 \subset J}^{\semi}$ is a subgroup of $\I_{L_2 \subset J}^{\semi}$, we
have a natural forgetful functor
\[
\oblv^{\I_{L_2 \subset J}^{\semi}|\I_{L_1 \subset J}^{\semi}}: \Shv \left( \I_{L_2 \subset J}^{\semi} \backslash X \right) \rightarrow \Shv \left( \I_{L_1 \subset J}^{\semi} \backslash X \right)
\]
induced by $!$-pullback along the quotient map $\I_{L_1 \subset J}^{\semi} \backslash X \to \I_{L_2 \subset J}^{\semi}\backslash X$,
which admits left and right adjoints 
\[
\Av_!^{\I_{L_2 \subset J}^{\semi}|\I_{L_1 \subset J}^{\semi}}, \Av_*^{\I_{L_2 \subset J}^{\semi}|\I_{L_1 \subset J}^{\semi}} : \Shv \left( \I_{L_1 \subset J}^{\semi} \backslash X \right) \rightarrow \Shv \left( \I_{L_2 \subset J}^{\semi} \backslash X \right). 
\]
Here $\Av_!^{\I_{L_2 \subset J}^{\semi}|\I_{L_1 \subset J}^{\semi}}$ is given by $!$-pushforward along the map $\I_{L_1 \subset J}^{\semi} \backslash X \to \I_{L_2 \subset J}^{\semi}\backslash X$, and $\Av_*^{\I_{L_2 \subset J}^{\semi}|\I_{L_1 \subset J}^{\semi}}$ is given by $*$-pushforward along this map, followed by $[-2\ell(w_{L_2}w_{L_1})]$. In fact, since the projection map is proper we have
\[
\Av_*^{\I_{L_2 \subset J}^{\semi}|\I_{L_1 \subset J}^{\semi}} = \Av_!^{\I_{L_2 \subset J}^{\semi}|\I_{L_1 \subset J}^{\semi}}[-2\ell(w_{L_2}w_{L_1})].
\]

We also have pullback and pushforward functors 
\begin{gather*}
(p_K)^*, (p_K)^!: \Shv \left( \I^\semi_{L \subset J} \backslash \Fl_K \right) \rightarrow \Shv \left(\I^\semi_{L \subset J}  \backslash \Fl \right), \\ 
(p_K)_* = (p_K)_!: \Shv \left( \I^\semi_{L \subset J} \backslash \Fl \right) \rightarrow \Shv \left( \I^\semi_{L \subset J}  \backslash \Fl_K \right), \\
q^*, q^!: \Shv \left( \I^\semi_{L \subset J} \backslash \Fl \right) \rightarrow \Shv \left( \I^\semi_{L \subset J} \backslash \tFl \right), \ 
q_*,q_!: \Shv \left( \I^\semi_{L \subset J}  \backslash \tFl \right) \rightarrow \Shv \left( \I^\semi_{L \subset J} \backslash \Fl \right). 
\end{gather*}
We have $(p_K)^! \simeq (p_K)^*[2 \ell(w_K)]$, and we set
\[
(p_K)^\dag = (p_K)^*[\ell(w_K)].
\]
Similarly we have $q^! \simeq q^*[2\dim(T)]$.

\subsection{Standard and costandard sheaves} 

For $L$, $J$, $K$ as above and $x \in {}^LW_\ex^{(K,J)}$,
we will also consider the $\infty$-category 
\begin{equation}
\label{eq:cat-orbt}
\Shv \left( \I^{\semi}_{L \subset J} \backslash \rS^{L \subset J}_{K,x} \right).
\end{equation}
As in~\cite[Proposition~1.3.5]{Gaits18} or \cite[Lemma~3.3.2]{DL23}, $!$-restriction to the point $\dot{x}\I_K/\I_K$ induces an equivalence 
\begin{equation}
\label{eq:semi-orbt}
\Shv \left( \I^{\semi}_{L \subset J} \backslash \rS^{L \subset J}_{K,x} \right) \xs \Shv \left( \mathrm{Stab}_{M_L}(\dot{x}\I_K/\I_K)\backslash \pt \right).
\end{equation}
To simplify notation we will write $\omega_{\rS^{L \subset J}_{K,x}}$ for $\omega_{\I^{\semi}_{L \subset J} \backslash \rS^{L \subset J}_{K,x}}$. (We will use similar conventions whenever convenient below.) Under the equivalence~\eqref{eq:semi-orbt}, this object corresponds to $\omega_{\mathrm{Stab}_{M_L}(\dot{x}\I_K/\I_K)\backslash \pt}$.
We define the \textit{perverse t-structure} on~\eqref{eq:cat-orbt} to be the cohomological shift of the transfer of the natural t-structure on the $\infty$-category $\Shv \left( \mathrm{Stab}_{M_L}(\dot{x}\I_K/\I_K)\backslash \pt \right)$ such that the object
\begin{equation}
\label{eqn:omega-orbit-perverse}
\omega_{\rS^{L \subset J}_{K,x}}[-\ell(w_Lw_{L\cap ({}^xK)})-\ell^J(x)]
\end{equation}
belongs to the heart. 
Note that when $J=S$, the forgetful functor $\Shv(\I_L \backslash \rS^{L \subset S}_{K,x}) \to \Shv(\rS^{L \subset S}_{K,x})$ is t-exact with respect to the usual perverse t-structure on $\Shv(\rS^{L \subset S}_{K,x})$, see~\eqref{eqn:dim-par-orbit}. In general we are dealing with schemes of infinite dimension, so that there is no natural choice for the shift appearing in~\eqref{eqn:omega-orbit-perverse}; as in~\cite{adr} the shifts are chosen so as to provide a nice theory. (In particular, when $L=J=\varnothing$, our conventions are consistent with those in~\cite{adr}.)

We set
\begin{align*}
\Delta^{L \subset J,\semi}_{K,x} &:= (\bs^{L \subset J}_{K,x})_! \omega_{\rS^{L \subset J}_{K,x}}[-\ell(w_L w_{L \cap ({}^xK)})-\ell^J(x)], \\ 
\nabla^{L \subset J,\semi}_{K,x} &:= (\bs^{L \subset J}_{K,x})_* \omega_{\rS^{L \subset J}_{K,x}}[-\ell(w_L w_{L \cap ({}^xK)})-\ell^J(x)].
\end{align*}
These objects will be called the standard and costandard sheaves in $\Shv \left( \I^{\semi}_{L \subset J} \backslash \Fl_{K} \right)$ respectively.
We will omit $K$ in the notation for these objects if it is the empty set, and will also omit ``$L \subset J$'' in case $L=J=\varnothing$ (which corresponds to the setting of~\cite{adr}). In case $L=\varnothing$, we will also replace  ``$\varnothing \subset J, \semi$'' by ``$J,\quar$,'' and when $L=J$ we will replace ``$J \subset J, \semi$'' by ``$J,\semi$.'' Finally, when $J=S$ (so that $\I_{L \subset S}^\semi = \I_L$) we will replace ``$L \subset S, \semi$'' by ``$L$'', and will omit $L$ if $L=\varnothing$ (so that we work with usual Iwahori-equivariant sheaves on $\Fl_K$). We will use similar conventions for the variants of these objects introduced below.

\begin{lem}
\label{lem:Av-quar-semi-D-N}
Let $L \subset J \subset S$, and let $K \subset S_\af$ be finitary.
Let $x\in {}^LW_\ex^{(K,J)}$. 
\begin{enumerate}
\item 
\label{it:Av-quar-semi-D-N-1}
Let $w\in W_L^{L\cap ({}^x K)}$. 
The object 
\[
\Av_!^{\I_{L \subset J}^{\semi}|\I_J^{\quar}}(\Delta^{J,\quar}_{K,wx}), \quad 
\text{resp.}\quad \Av_*^{\I_{L \subset J}^{\semi}|\I_J^{\quar}}(\nabla^{J,\quar}_{K,wx}),
\]
is an iterated extension of objects
\[
\Delta^{L \subset J,\semi}_{K,x}[\ell(w_{L\cap ({}^xK)}u)+\ell(w_Lwu)], \
\text{resp.}\ \nabla^{L \subset J,\semi}_{K,x}[-\ell(w_{L\cap ({}^xK)}u)-\ell(w_Lwu)],
\]
with $u\in W_{L\cap ({}^xK)}$, each of them appearing once. 
\item 
\label{it:Av-quar-semi-D-N-2}
The object 
\[
\oblv^{\I_{L \subset J}^{\semi}|\I_J^{\quar}}(\Delta^{L \subset J,\semi}_{K,x}), 
\quad \text{resp.}\quad 
\oblv^{\I_{L \subset J}^{\semi}|\I_J^{\quar}}(\nabla^{L \subset J,\semi}_{K,x}),
\]
is an iterated extension of objects
\[
\Delta^{J,\quar}_{K,wx}[\ell(w_L w_{L\cap ({}^xK)})-\ell(w)], 
\ \text{resp.}\ 
\nabla^{J,\quar}_{K,wx}[-\ell(w_L w_{L \cap ({}^xK)})+\ell(w)],
\]
with $w\in W_L^{L\cap ({}^xK)}$, each of them appearing once. 
\end{enumerate}
\end{lem}

\begin{proof}
\eqref{it:Av-quar-semi-D-N-1}
Lemma~\ref{lem:bi-min-av}\eqref{it:bi-min-av-1} allows to reduce the proof to the case $w=1$. Then we use Lemma~\ref{lem:bi-min-av}\eqref{it:bi-min-av-2}. For the first claim we have to compute the $!$-pushforward under that map of $\omega_{\I_{L \subset J}^\semi \times^{\I_J^\quar} \rS_{K,x}^{J}}$, and the statement follows from the fact that the cohomology of $\omega_{M_{L \cap ({}^x K)}/B_{L \cap ({}^x K)}}$ is $\bigoplus_{u \in W_{L \cap ({}^x K)}} \k[2\ell(w_{L \cap ({}^x K)}) - 2\ell(u)]$. For the second claim we use the fact that $\Av_*^{\I_{L \subset J}^{\semi}|\I_J^{\quar}}$ is the composition of the appropriate $*$-pushforward functor with 
$[-2\ell(w_L)]$ (see~\S\ref{ss:qs-sheaves}), and write now the cohomology above as $\bigoplus_{u \in W_{L \cap ({}^x K)}} \k[2\ell(u)]$.

\eqref{it:Av-quar-semi-D-N-2}
Using the decomposition from Lemma~\ref{lem:bi-min-av}\eqref{it:bi-min-av-3}, we obtain that $\omega_{\rS^{L \subset J}_{K,x}}$ is an iterated extension of objects $f_! f^* \omega_{\rS^{L \subset J}_{K,x}}$, resp.~$f_* f^! \omega_{\rS^{L \subset J}_{K,x}}$, where $f$ runs over the locally closed immersions $\rS^{J}_{K,wx} \hookrightarrow \rS^{L \subset J}_{K,x}$ ($w \in W_L^{L \cap ({}^xK)}$). The second assertion immediately follows. For the first one, we need to compute the complexes $f^* \omega_{\rS^{L \subset J}_{K,x}}$ for $f$ as above.
Using Lemma~\ref{lem:bi-min-av}\eqref{it:bi-min-av-1},
we write the embedding $\rS^{J}_{K,wx} \hookrightarrow \rS^{L \subset J}_{K,x}$ as the composition
\[
\rS^{J}_{K,wx} \hookrightarrow \I_{L \subset J}^\semi \times^{\I_J^\quar} \rS_{K,x}^J \to \rS^{L \subset J}_{K,x}.
\]
The second map is smooth of relative dimension $\ell(w_{L \cap ({}^x K)})$ by Lemma~\ref{lem:bi-min-av}\eqref{it:bi-min-av-2}, so the $*$-pullback under this map of $\omega_{\rS^{L \subset J}_{K,x}}$ is $\omega_{\I_{L \subset J}^\semi \times^{\I_J^\quar} \rS_{K,x}^J}[-2\ell(w_{L \cap ({}^x K)})]$. On the other hand, the $*$-pullback under the first map of $\omega_{\I_{L \subset J}^\semi \times^{\I_J^\quar} \rS_{K,x}^J}$ is $\omega_{\rS^{J}_{K,wx}}[2\ell(w_L)-2\ell(w)]$. The claim follows.
\end{proof}

\begin{lem}
\label{lem:push-pull-pK-D-N-semi}
Let $L \subset J \subset S$, and let $K \subset S_\af$ be finitary.
Let $x\in {}^L W_\ex^{(K,J)}$. 
\begin{enumerate}
\item 
\label{it:push-pull-pK-D-N-semi-1}
Let $w\in {}^{(L^x)\cap K}W_K$. 
The object 
\[
(p_K)_*(\Delta^{L \subset J,\semi}_{xw})[\ell(w_K)], 
\quad \text{resp.}\quad 
(p_K)_*(\nabla^{L \subset J, \semi}_{xw})[-\ell(w_K)], 
\]
is an iterated extension of objects
\[
\Delta^{L \subset J,\semi}_{K,x}[\ell(uw_{(L^x)\cap K})+\ell(uww_K)], \
\text{resp.}\ \nabla^{L \subset J,\semi}_{K,x}[-\ell(uw_{(L^x)\cap K})-\ell(uww_K)],
\]
with $u\in W_{(L^x)\cap K}$, each appearing once. 
\item
\label{it:push-pull-pK-D-N-semi-2}
The object 
\[
(p_K)^\dag(\Delta^{L \subset J,\semi}_{K,x}), 
\quad \text{resp.}\quad 
(p_K)^\dag(\nabla^{L \subset J,\semi}_{K,x}), 
\]
is an iterated extension of objects
\[
\Delta^{L \subset J, \semi}_{xw}[\ell(w_{(L^x)\cap K}w_K)-\ell(w)], \
\text{resp.}\ \nabla^{L \subset J, \semi}_{xw}[-\ell(w_{(L^x)\cap K}w_K)+\ell(w)],
\]
with $w\in {}^{(L^x)\cap K}W_K$, each appearing once. 
\end{enumerate}
\end{lem}

\begin{proof}
\eqref{it:push-pull-pK-D-N-semi-1}
Using Lemma~\ref{lem:bi-min-push}\eqref{it:bi-min-push-1}, we reduce the proof to the case $w=1$. In that case we use Lemma~\ref{lem:bi-min-push}\eqref{it:bi-min-push-2} to see that $(p_K)_*(\Delta^{L \subset J,\semi}_{x})$ is an iterated extension of objects $\Delta^{L \subset J,\semi}_{K,x}[-\ell(w_{L \cap ({}^x K)}) + 2\ell(w_{(L^x) \cap K})-2\ell(u)]$ with $u\in W_{(L^x)\cap K}$. 
The first assertion follows, using the fact that $\ell(w_{L \cap ({}^x K)}) = \ell(w_{(L^x) \cap K})$. 
The second assertion follows from the same considerations (and replacing $u$ by $u w_{(L^x) \cap K}$ in the formula above).

\eqref{it:push-pull-pK-D-N-semi-2}
By the base change theorem and smoothness of $p_K$, the complex $(p_K)^\dag(\Delta^{L \subset J,\semi}_{K,x})$, resp.~$(p_K)^\dag(\nabla^{L \subset J,\semi}_{K,x})$, is the $!$-pushforward, resp.~$*$-pushforward, of
\[
\omega_{(p_K)^{-1}(\rS^{L \subset J}_{K,x})}[-\ell(w_L)+\ell(w_{L \cap ({}^x K)})-\ell^J(x)-\ell(w_K)]
\]
along the immersion $(p_K)^{-1}(\rS^{L \subset J}_{K,x}) \hookrightarrow \Fl$. Now, by the decomposition in Lem\-ma~\ref{lem:bi-min-push}\eqref{it:bi-min-push-3}, $\omega_{(p_K)^{-1}(\rS^{L \subset J}_{K,x})}$ is an iterated extension of $*$-pushforwards of complexes $\omega_{\rS^{L \subset J}_{xw}}$ with $w \in {}^{(L^x) \cap K} W_K$, which implies the second assertion using the formula $\ell(w_{L \cap ({}^x K)}) = \ell(w_{(L^x) \cap K})$. Similarly, $\omega_{(p_K)^{-1}(\rS^{L \subset J}_{K,x})}$ is an iterated extension of objects $f_! f^* \omega_{(p_K)^{-1}(\rS^{L \subset J}_{K,x})}$ where $f$ runs over the immersions $\rS^{L \subset J}_{xw} \hookrightarrow (p_K)^{-1}(\rS^{L \subset J}_{K,x})$, and arguments similar to those encountered in the proof of Lemma~\ref{lem:Av-quar-semi-D-N}\eqref{it:Av-quar-semi-D-N-2}
imply the first assertion.
\end{proof}

\subsection{Whittaker sheaves} 
\label{ss:Whittaker}

In this subsection we assume that $p>0$, and that $\k$ contains a nontrivial $p$-th root of unity.

\subsubsection{Definition} 
\label{sss:Whittaker}

Let $A\subset S_\af$ be a finitary subset.
As in~\cite[\S 3.4]{AR23}, for any $s \in S_\af$ we have an associated (negative, affine) root subgroup $U_s^- \subset \sL G$ and a (noncanonical, but fixed) isomorphism $\Ga \simeq U_s^-$.
Recall from~\S\ref{sss:def-Fl} that we have a parahoric subgroup $\I_A \subset \sL G$ associated with $A$, with pro-unipotent radical $\Iu \cap \Iu^A$ and Levi factor $M_A$, which therefore identifies with the quotient $\I_A / (\Iu \cap \Iu^A)$. The image $U_A^-$ of $\Iu^A$ in this quotient is the unipotent radical of a Borel subgroup, and as in~\cite[\S 3.4]{AR23} we consider the composition
\[
\Iu^A \to U_A^- \to U_A^- / [U_A^-, U_A^-] \simeq \prod_{s \in A} U_s^- \simeq (\Ga)^{\# A} \xrightarrow{\mathrm{sum}} \Ga
\]
where the first two maps are the natural quotient morphisms.
Under our assumptions on $p$ and $\k$ the local system $\mathrm{AS}_* \underline{\k}_{\Ga}$, where $\mathrm{AS} : \Ga \to \Ga$ is the Artin--Schreier map sending $x$ to $x^p-x$, splits as a direct sum of $p$ rank-$1$ local systems, parametrized by morphisms from $(\mathbb{F}_p,+)$ to $\k^\times$. We choose one nontrivial irreducible direct summand,
and denote by $\chi_A$, resp.~$\chi'_A$, its pullback to $\Iu^A$, resp.~to $U_A^-$. 

Combining the constructions from~\S\ref{ss:qs-sheaves} with the consideration of Whittaker shea\-ves as in~\cite[End of~\S 4.2]{zhu},
one can define for any subsets $L \subset J \subset S$ the $\infty$-category
\begin{equation}
\label{eqn:Whit-cat}
\Shv \left( \I_{L \subset J}^{\semi}\backslash \sL G/(\Iu^A,\chi_A) \right). 
\end{equation}
More specifically, let $\widetilde{\Iu^A} := \Iu^A \times_{\Ga} \Ga$, where the map $\Ga \to \Ga$ is $\mathrm{AS}$. Then one can consider the stack $\I_{L \subset J}^\semi \backslash \sL G / \widetilde{\Iu^A}$. The $\infty$-category
\[
\Shv \left( \I_{L \subset J}^\semi \backslash \sL G / \widetilde{\Iu^A} \right)
\]
breaks as a direct sum indexed by the $p$ different irreducible direct summands of $\mathrm{AS}_* \underline{\k}_{\Ga}$, and $\Shv \left( \I_{L \subset J}^{\semi}\backslash \sL G/(\Iu^A,\chi_A) \right)$ is defined as the direct summand corresponding to our choice of Artin--Schreier local system. As in~\S\ref{ss:qs-sheaves} we have functors
\[
\oblv^{\I_{L_2 \subset J}^{\semi}|\I_{L_1 \subset J}^{\semi}}, \Av_!^{\I_{L_2 \subset J}^{\semi}|\I_{L_1 \subset J}^{\semi}}, \Av_*^{\I_{L_2 \subset J}^{\semi}|\I_{L_1 \subset J}^{\semi}}
\]
relating the $\infty$-categories~\eqref{eqn:Whit-cat} when $L$ varies.

It is convenient to introduce the ind-scheme 
\[
\uFl_A := \sL G / (\Iu \cap \Iu^A),
\]
which admits an $\sL G$-action from the left and an $M_A$-action from the right; then we have
\[
\Shv \left( \I_{L \subset J}^{\semi}\backslash \sL G/(\Iu^A,\chi_A) \right) \simeq \Shv \left( \I_{L \subset J}^{\semi} \backslash \uFl_A /(U_A^-,\chi'_A) \right)
\]
where the right-hand side is defined as above.
Denote by $r_A: \uFl_A \rightarrow \Fl$ the natural projection. 
By~\eqref{eq:decomp} we have a decomposition
\begin{equation}
\label{eq:decomp-tGr}
\uFl_A= \bigsqcup_{x\in {}^L W_\ex} (r_A)^{-1}(\rS^{L \subset J}_{x}) \cdot \dot{w}_A
\end{equation}
into $(\I_{L \subset J}^\semi \times (U_A^-)^\mathrm{op})$-orbits.

\begin{lem}
\label{lem:whit}
For $x\in {}^L W_\ex$, $!$-restriction to the point $\dot{x}\dot{w}_A$ induces an equivalence 
\begin{equation*}
\Shv\left( \I_{L \subset J}^{\semi} \backslash (r_A)^{-1}(\rS^{L \subset J}_x)\dot{w}_A / (U_A^-,\chi'_A) \right)
\simeq
\begin{cases}
\Vect_\k & \text{if $x\in {}^L W_\ex^{(A,J),\max}$;} 
\\ 0 & \text{otherwise.} 
\end{cases}
\end{equation*}
\end{lem}

\begin{proof}
This follows from standard arguments; namely, the same arguments as in~\cite[Proposition~1.3.5]{Gaits18} allow to reduce the proof to a finite-dimensional setting, and then one uses the same arguments as in~\cite[Lemma~4.2.1]{by}. In the end we use Lemma~\ref{lem:quar-order} and Lemma~\ref{lem:quar-min-max}, which ensure that $x \in {}^L W_\ex^{(A,J),\max}$ if and only if we have $\Ad_{\dot{x}\dot{w}_A}(U_s^-) \cap \I_{L \subset J}^{\semi} = \{1\}$ for any $s \in A$.
\end{proof}

If $x \in {}^L W_\ex^{(A,J),\max}$, we will denote by 
\[
\omega_{(r_A)^{-1}(\rS^{L \subset J}_x)\dot{w}_A}^{\chi_A} \in \Shv\left( \I_{L \subset J}^{\semi} \backslash (r_A)^{-1}(\rS^{L \subset J}_x)\dot{w}_A / (U_A^-,\chi'_A) \right)
\]
the object corresponding to $\k\in \Vect_\k$ under the equivalence of Lemma~\ref{lem:whit}.
Let $\widetilde{\bs}^{L \subset J}_{A,x}$ be the embedding of $(r_A)^{-1}(\rS^{L \subset J}_x)\dot{w}_A$ in $\uFl_A$. 
We define the standard and costandard sheaves in $\Shv \left( \I_{L \subset J}^{\semi}\backslash \sL G/(\Iu^A,\chi_A) \right)$ as 
\begin{align*}
\Delta_{\chi_A,x}^{L \subset J, \semi} &:= (\widetilde{\bs}^{L \subset J}_{A,x})_! \omega_{(r_A)^{-1}(\rS^{L \subset J}_x)\dot{w}_A}^{\chi_A} [-\ell^J(x)-\ell(w_L)],
\\
\nabla_{\chi_A,x}^{L \subset J, \semi} &:= (\widetilde{\bs}^{L \subset J}_{A,x})_* \omega_{(r_A)^{-1}(\rS^{L \subset J}_x)\dot{w}_A}^{\chi_A} [-\ell^J(x)-\ell(w_L)].
\end{align*}

It follows from Lemma~\ref{lem:quar-min-max} that if $x\in W_\ex^{(A,J),\max}$, there exists at most one element $w \in W_L$ such that $wx \in {}^L W_\ex^{(A,J),\max}$; on the other hand, if $x\in {}^L W_\ex^{(A,J),\max}$, then $wx \in W_\ex^{(A,J),\max}$ for any $w \in W_L$.

\begin{lem}
\label{lem:av-oblv-whit}
Let $L \subset J \subset S$, and let $A \subset S_\af$ be a finitary subset.
\begin{enumerate}
\item 
\label{it:av-oblv-whit-1}
Let $x\in W_\ex^{(A,J),\max}$. 
We have 
\[
\Av_!^{\I_{L \subset J}^{\semi}|\I_J^{\quar}}(\Delta^{J,\quar}_{\chi_A,x})\simeq
\begin{cases}
\Delta^{L \subset J,\semi}_{\chi_A, wx}[\ell(w_Lw)] & \text{if $w\in W_L$ satisfies $wx\in {}^LW_\ex^{(A,J),\max}$;} \\ 
0 & \text{if $W_L x \cap {}^LW_\ex^{(A,J),\max} = \varnothing$,}
\end{cases}
\]
and 
\[
\Av_*^{\I_{L \subset J}^{\semi}|\I_J^{\quar}}(\nabla^{J,\quar}_{\chi_A,x})\simeq
\begin{cases}
\nabla^{L \subset J,\semi}_{\chi_A, wx}[-\ell(w_Lw)] & \text{if $w\in W_L$ satisfies $wx\in {}^LW_\ex^{(A,J),\max}$;} \\ 
0 & \text{if $W_L x \cap {}^LW_\ex^{(A,J),\max} = \varnothing$.} 
\end{cases}
\]
\item 
\label{it:av-oblv-whit-2}
Let $x\in {}^L W_\ex^{(A,J),\max}$. 
Then $\oblv^{\I_{L \subset J}^{\semi}|\I_J^{\quar}}(\Delta^{L \subset J,\semi}_{\chi_A,x})$, resp. $\oblv^{\I_{L \subset J}^{\semi}|\I_J^{\quar}}(\nabla^{L \subset J,\semi}_{\chi_A,x})$, is an iterated extension of objects $\Delta^{J,\quar}_{\chi_A, wx}[\ell(w_Lw)]$, resp. $\nabla^{J,\quar}_{\chi_A,wx}[-\ell(w_Lw)]$, with $w\in W_L$, each appearing once. 
\end{enumerate}
\end{lem}

\begin{proof}
\eqref{it:av-oblv-whit-1}
We only consider the first assertion; the second one can be justified similarly.
The case when $W_L x \cap {}^L W_\ex^{(A,J),\max}=\varnothing$ directly follows from Lemma~\ref{lem:whit}. Now, let us assume that there exists $w\in W_L$ such that $wx\in {}^L W_\ex^{(A,J),\max}$. Using Lemma~\ref{lem:bi-min-av}\eqref{it:bi-min-av-1} we obtain that
\[
\Av_!^{\I_{L \subset J}^{\semi}|\I_J^{\quar}}(\Delta^{J,\quar}_{\chi_A,x}) \simeq \Av_!^{\I_{L \subset J}^{\semi}|\I_J^{\quar}}(\Delta^{J,\quar}_{\chi_A,wx})[-\ell(w)],
\]
so we can assume $w=1$. In this case, the claim follows from Lemma~\ref{lem:bi-min-av}\eqref{it:bi-min-av-2} (applied with $K=\varnothing$, so that the flag variety of the Levi subgroup under consideration is a point).

\eqref{it:av-oblv-whit-2}
The proof is the same as that of Lemma~\ref{lem:Av-quar-semi-D-N}\eqref{it:Av-quar-semi-D-N-2}, using
Lemma~\ref{lem:bi-min-av}\eqref{it:bi-min-av-3}.
\end{proof}

\subsubsection{Whittaker averaging functors} 
\label{sss:av-fun}

Let $L$, $J$, $A$ be as above.
We have a forgetful functor 
\[
\oblv^{(\Iu^A,\chi_A)\mid \Iu^A \cap \Iu}: 
\Shv \left( \I_{L \subset J}^{\semi}\backslash \sL G/(\Iu^A,\chi_A) \right) \rightarrow \Shv \left( \I_{L \subset J}^{\semi}\backslash \sL G/(\Iu^A \cap \Iu) \right)
\]
which admits left and right adjoints 
\begin{multline*}
\Av^{(\Iu^A,\chi_A)\mid \Iu^A \cap \Iu}_!, \Av^{(\Iu^A,\chi_A)\mid \Iu^A \cap \Iu}_* : \Shv \left( \I_{L \subset J}^{\semi}\backslash \sL G/(\Iu^A \cap \Iu) \right) \\
\to \Shv \left( \I_{L \subset J}^{\semi}\backslash \sL G/(\Iu^A,\chi_A) \right).
\end{multline*}
Similarly we have the forgetful functor 
\[
\oblv^{\Iu\mid \Iu^A \cap \Iu}: \Shv \left( \I_{L \subset J}^{\semi} \backslash \tFl \right) \rightarrow \Shv \left( \I_{L \subset J}^{\semi}\backslash \sL G/(\Iu^A \cap \Iu) \right),
\]
and its left and right adjoints 
\[
\Av^{\Iu\mid \Iu^A \cap \Iu}_!, \Av^{\Iu\mid \Iu^A \cap \Iu}_* : \Shv \left( \I_{L \subset J}^{\semi}\backslash \sL G/(\Iu^A \cap \Iu) \right) \to \Shv \left( \I_{L \subset J}^{\semi} \backslash \tFl \right).
\]
We set 
\begin{align*}
\Av^A_{\chi,!} &:= \Av^{(\Iu^A,\chi_A)\mid \Iu^A \cap \Iu}_! \circ \oblv^{\Iu\mid \Iu^A \cap \Iu}[-\ell(w_A)], \\ 
\Av^A_{\chi,*} &:= \Av^{(\Iu^A,\chi_A)\mid \Iu^A \cap \Iu}_* \circ \oblv^{\Iu\mid \Iu^A \cap \Iu}[\ell(w_A)], \\
\Av^A_{!} &:= \Av^{\Iu\mid \Iu^A \cap \Iu}_! \circ \oblv^{(\Iu^A,\chi_A)\mid \Iu^A \cap \Iu}[-\ell(w_A)], \\ 
\Av^A_{*} &:= \Av^{\Iu\mid \Iu^A \cap \Iu}_* \circ \oblv^{(\Iu^A,\chi_A)\mid \Iu^A \cap \Iu}[\ell(w_A)].
\end{align*}
Then we have adjoint pairs $(\Av^A_{\chi,!},\Av^A_{*})$ and $(\Av^A_{!},\Av^A_{\chi,*})$. 

The following lemma follows from standard arguments, as e.g.~in~\cite[Lemma~3.2]{AR23}.

\begin{lem}
There exists a canonical isomorphism of functors
\[
\Av^A_{\chi,!} \simeq \Av^A_{\chi,*}.
\]
\end{lem}

In view of this lemma, we will denote both $\Av^A_{\chi,!}$ and $\Av^A_{\chi,*}$ by $\Av^A_\chi$.
We will say more about these averaging functors in~\S\ref{sss:prop-av} below. 

When $J=S$, we are in the traditional setting of Whittaker sheaves on (partial) affine flag varieties as studied e.g.~in~\cite[\S 3]{AR23} (except that the roles of multiplication on the left and on the right are switched compared to this reference). In particular, it follows from~\cite[Lemma~3.3]{AR23} (or, more explicitly, an immediate generalization of this statement to all parabolic affine flag varieties) that for $x \in {}^L W_\ex$ we have
\begin{equation}
\label{eqn:Av-chi-D-classical}
\Av^A_\chi(q^! \Delta^L_{x}) \simeq
\begin{cases}
\Delta^L_{\chi_A,xw} & \text{if $w \in W_A$ satisfies $xw \in {}^L W_\ex^{A,\max}$;} \\
0 & \text{if $xW_A \cap {}^L W_\ex^{A,\max} = \varnothing$,}
\end{cases}
\end{equation}
and similarly for costandard objects.

\subsection{Convolution actions} 
\label{ss:convolution}

We now drop the assumptions on $p$ and $\k$ from~\S\ref{ss:Whittaker}.
For any $K \subset S_\af$ finitary,
the $\infty$-category $\Shv(\I_K \backslash \sL G/\I_K)$ admits a canonical structure of algebra object in $\k$-linear presentable stable $\infty$-categories, with binary product given by convolution and denoted $\star^{\I_K}$ (see the discussion and references in~\cite[\S 2.5]{adr}). For any $L \subset J \subset S$, we have a natural right action of $\Shv \left( \I_K \backslash \sL G/\I_K \right)$ on $\Shv \left( \I^\semi_{L \subset J} \backslash \Fl_K \right)$; the corresponding binary operation
\[
\Shv \left( \I^\semi_{L \subset J} \backslash \Fl_K \right) \times \Shv \left( \I_K \backslash \sL G/\I_K \right) \to \Shv \left( \I^\semi_{L \subset J} \backslash \Fl_K \right)
\]
will also be denoted $\star^{\I_K}$.

Let now $X$ be either $\Fl_K$ for some finitary $K \subset S_\af$, $\tFl$, or the symbol $\sL G/(\Iu^A,\chi_A)$ for some finitary $A \subset S_\af$ (in this case we assume that $p>0$ and that $\k$ contains a nontrivial $p$-th root of unity). 
For any $L \subset J \subset S$ we have a natural action by convolution from the left 
\[
(-) \star^{\I_{L \subset J}} (-) : \Shv \left( \I_{L \subset J}\backslash \sL{M_J}/\I_{L \subset J} \right) \times \Shv \left( \I_{L \subset J}^{\semi} \backslash X \right) \rightarrow \Shv \left( \I_{L \subset J}^{\semi}\backslash X \right).
\]

If $\theta\in \Xd_J^0$, then $z^\theta$ commutes with $\sL^+ M_J$ and normalizes $\sL U_J$. Hence pushforward under the action of this element on $X$ induces an autoequivalence of $\Shv \left( \I_{L \subset J}^{\semi} \backslash X \right)$, which we will denote $z^\theta \cdot (-)$. We will denote by
$\langle \theta \rangle$ the auto-equivalence given by $\sF\mapsto z^\theta\cdot \sF [-\langle \theta, 2\rho \rangle]$. It is clear (using Remark~\ref{rmk:lJ}\eqref{it:lJ-translation}) that for any $\theta \in \Xd_J^0$ and $x \in {}^L W_\ex^{(K,J)}$ we have
\begin{equation}
\label{eqn:D-N-shift}
\Delta^{L \subset J, \semi}_{K,x} \langle \theta \rangle \simeq \Delta^{L \subset J, \semi}_{K,t_\theta x}, \quad
\nabla^{L \subset J, \semi}_{K,x} \langle \theta \rangle \simeq \nabla^{L \subset J, \semi}_{K,t_\theta x}.
\end{equation}

\section{Raskin's equivalences} 
\label{sec:Raskin}

Let $X$ be either $\Fl_K$ for some finitary $K \subset S_\af$, $\tFl$, or the symbol $\sL G/(\Iu^A,\chi_A)$ for some finitary $A \subset S_\af$ (assuming in this case that $p>0$ and that $\k$ contains a nontrivial $p$-th root of unity).

\subsection{Construction and functorialities} 

The following proposition generalizes constructions
discussed in~\cite[\S 5.2]{Gaits18},~\cite[Proposition~3.1.7]{DL23} and~\cite[\S2.7]{adr}, which themselves generalize an earlier construction of Raskin. These equivalences sometimes allow one to reduce the proof of properties of semiinfinite sheaves to the more familiar case of sheaves equivariant with respect to a parahoric subgroup.

\begin{prop}
\label{prop:Raskin}
For any subsets $L \subset J \subset S$, there exists a canonical equivalence of $\infty$-categories
\[
\Ras^{\I_{L \subset J}^{\semi} \mid \I_L}: \Shv \left( \I_{L \subset J}^{\semi} \backslash X \right) \simto \Shv \left( \I_L \backslash X \right).
\]
If $L_1 \subset L_2 \subset J \subset S$ we have a commutative diagram
\begin{equation}
\label{eq:comm-diag-ras-oblv}
\begin{tikzcd}[column sep=3cm]
\Shv \left( \I_{L_2 \subset J}^{\semi}\backslash X \right) \ar[r, "\Ras^{\I_{L_2 \subset J}^{\semi} \mid \I_{L_2}}", "\sim"'] \ar[d, "\oblv^{\I_{L_2 \subset J}^{\semi}|\I_{L_1 \subset J}^{\semi}}"'] & \Shv(\I_{L_2} \backslash X) \ar[d, "\oblv^{\I_{L_2}|\I_{L_1}}"] \\
\Shv \left( \I_{L_1 \subset J}^{\semi}\backslash X \right) \ar[r, "\Ras^{\I_{L_1 \subset J}^{\semi} \mid \I_{L_1}}", "\sim"'] & \Shv(\I_{L_1} \backslash X).
\end{tikzcd}
\end{equation}
\end{prop}

For $L \subset J \subset S$,
we will denote by $\Ras^{\I_L \mid \I_{L \subset J}^{\semi}}$ the inverse of $\Ras^{\I_{L \subset J}^{\semi} \mid \I_L}$.
If $J_1, J_2 \subset S$ and $J_1 \supset L \subset J_2$, we also set
\[
\Ras^{\I_{L \subset J_1}^{\semi} \mid \I_{L \subset J_2}^{\semi}} := \Ras^{\I_L \mid \I_{L \subset J_2}^{\semi}} \circ \Ras^{\I_{L \subset J_1}^{\semi} \mid \I_L} : \Shv \left( \I_{L \subset J_1}^{\semi} \backslash X \right) \simto \Shv \left( \I_{L \subset J_2}^\semi \backslash X \right).
\] 

The proof of Proposition~\ref{prop:Raskin} is similar to that of its counterparts in~\cite{Gaits18,DL23,adr}, so we will not repeat it, but only explain some aspects that will be relevant below.
Consider the natural forgetful functors 
\[
\Shv \left( \I_{L \subset J}^{\semi} \backslash X \right) \xrightarrow{\oblv^{\I_{L \subset J}^{\semi} \mid \I_{L \subset J}}} 
\Shv \left( \I_{L \subset J} \backslash X \right) \xleftarrow{\oblv^{\I_L \mid \I_{L \subset J}}}
\Shv \left( \I_L \backslash X \right). 
\]
The functor $\oblv^{\I_L \mid \I_{L \subset J}}$ admits a continuous right adjoint 
\[
\Av^{\I_L \mid \I_{L \subset J}}_*: \Shv(\I_{L \subset J}\backslash X)\rightarrow \Shv(\I_L \backslash X),
\]
and the equivalence $\Ras^{\I_{L \subset J}^{\semi} \mid \I_L}$ is given by 
\[
\Ras^{\I_{L \subset J}^{\semi} \mid \I_L} := \Av^{\I_L \mid \I_{L \subset J}}_*  \circ  \oblv^{\I_{L \subset J}^{\semi} \mid \I_{L \subset J}}. 
\]
In the other direction, the functor $\oblv^{\I_{L \subset J}^{\semi} \mid \I_{L \subset J}}$ admits a (continuous) left adjoint 
\[
\Av_!^{\I_{L \subset J}^{\semi} | \I_{L \subset J}}: \Shv \left( \I_{L \subset J}\backslash X \right) \to \Shv \left( \I_{L \subset J}^{\semi} \backslash X \right),
\]
and we have 
\[
\Ras^{\I_L \mid \I_{L \subset J}^{\semi}} = 
\Av^{\I_{L \subset J}^{\semi} \mid \I_{L \subset J}}_! \circ \oblv^{\I_L \mid \I_{L \subset J}}. 
\]
More explicitly, 
let $\theta\in \Xd_J^0 \cap \Xd_{J^\comp}^{++}$, and for $n \geq 0$ set 
\[
\I_{L \subset J,(n)}^{\semi}=\I_{L \subset J}\ltimes \Ad_{z^{-n\theta}}(\sL^+U_{J});
\]
then $( \I_{L \subset J,(n)}^{\semi} : n \geq 0)$ forms an exhaustive filtration of $\I_{L \subset J}^{\semi}$ by subgroup schemes, so that as in~\cite[\S 2.6]{adr} we have
\[
\Shv \left( \I_{L \subset J}^{\semi} \backslash X \right) = \lim_n \Shv \left( \I_{L \subset J, (n)}^{\semi} \backslash X \right).
\]
For any $\sF\in \Shv(\I_L \backslash X)$, we have 
\begin{equation}
\label{eq:Raskin-formula}
\begin{aligned}
\Ras^{\I_L \mid \I_{L \subset J}^{\semi}}(\sF)
&\simeq \colim_{n} \, \Av^{\I_{L \subset J,(n)}^{\semi} \mid \I_{L \subset J}}_!\circ \oblv^{\I_L \mid \I_{L \subset J}}(\sF) \\
&\simeq \colim_{n} \, z^{-n\theta} \Delta^L_{L,t_{n\theta}} \star^{\I_L} \sF [n\cdot \ell(t_{\theta})]
\end{aligned}
\end{equation}
(where the second identification can be justified e.g.~using~\eqref{eqn:decomposition-IL-J} below).

\begin{rmk}
\phantomsection
\label{rmk:Raskin}
\begin{enumerate}
\item
\label{it:ttheta-minimal}
For any $\theta \in \Xd_J^0$, we have $t_\theta \in {}^J W_\ex^J$. Hence one can consider $\Delta^J_{J,t_\theta}$, and we have $(p_J)_* \Delta_{t_\theta} \simeq \Delta^J_{J,t_\theta}$. In particular, the object $\Delta^L_{L,t_{n\theta}}$ considered in~\eqref{eq:Raskin-formula} indeed makes sense, and for $\sF$ as above we have $\oblv^{\I_L \mid \I}(\Delta^L_{L,t_{n\theta}} \star^{\I_L} \sF) \simeq \Delta_{t_{n\theta}} \star^{\I} \oblv^{\I_L \mid \I}(\sF)$.
\item
\label{it:rmk-Raskin-2}
The functor $\Ras^{\I_{L \subset J}^{\semi} \mid \I_L}$ can also be described as the composition
\begin{multline*}
\Shv \left( \I_{L \subset J}^{\semi} \backslash X \right) \xrightarrow{\oblv^{\I_{L \subset J}^{\semi} \mid \I_{L \subset J} \ltimes \sL^+ U_J}} 
\Shv \left( \I_{L \subset J} \ltimes \sL^+ U_J \backslash X \right) \\
\xrightarrow{\Av_*^{\I_L \mid \I_{L \subset J} \ltimes \sL^+ U_J}}
\Shv \left( \I_L \backslash X \right). 
\end{multline*}
where the first functor is the forgetful functor and the second functor is the right adjoint to the corresponding forgetful functor. In fact, this follows from the fact that the forgetful functor $\Shv \left( \I_{L \subset J} \ltimes \sL^+ U_J \backslash X \right) \to \Shv \left( \I_{L \subset J} \backslash X \right)$ is fully faithful, since $\sL^+ U_J$ is unipotent. Similarly, $\Ras^{\I_L \mid \I_{L \subset J}^{\semi}}$ identifies with the composition
\begin{multline*}
\Shv \left( \I_L \backslash X \right) \xrightarrow{\oblv^{\I_L \mid \I_{L \subset J} \ltimes \sL^+ U_J}}
\Shv \left( \I_{L \subset J} \ltimes \sL^+ U_J \backslash X \right) \\
\xrightarrow{\Av_!^{\I_{L \subset J}^{\semi} \mid \I_{L \subset J} \ltimes \sL^+ U_J}} 
\Shv \left( \I_{L \subset J}^{\semi} \backslash X \right). 
\end{multline*}
\item
\label{it:rmk-Raskin-3}
Denote by $U_J^-$ the unipotent radical of the parabolic subgroup of $G$ opposite to $P_J$ with respect to $T$. Then
a version of the Iwahori decomposition gives that multiplication induces an isomorphism
\begin{equation}
\label{eqn:decomposition-IL-J}
\sL^{++} U_J^- \times \I_{L \subset J} \times \sL^+ U_J \simto \I_L.
\end{equation}
Combining this with the considerations in~\eqref{it:rmk-Raskin-2}, we obtain that for any $\sF$ in $\Shv \left( \I_{L \subset J}^{\semi} \backslash X \right)$ we have
\[
\Ras^{\I_{L \subset J}^{\semi} \mid \I_L}(\sF) \simeq \Av^{\sL^{++} U_J^-}_*(\sF)
\]
where we omit the forgetful functors, and as usual $\Av^{\sL^{++} U_J^-}_*$ is the right adjoint to the corresponding forgetful functor. Similarly, for $\sF$ in $\Shv \left( \I_L \backslash X \right)$ we have
\[
\Ras^{\I_L \mid \I_{L \subset J}^{\semi}}(\sF) \simeq \Av^{\sL^{--} U_J}_!(\sF).
\]
\item
Consider subsets $L \subset J_1 \subset J_2 \subset S$. Then by~\eqref{it:rmk-Raskin-3} $\Ras^{\I_{L \subset J_2}^{\semi} \mid \I_L}$ is the restriction of $\Av^{\sL^{++} U_{J_2}^-}_*$, and $\Ras^{\I_L \mid \I_{L \subset J_1}^{\semi}}$ is the restriction of $\Av^{\sL^{--} U_{J_1}}_!$; as a consequence $\Ras^{\I_{L \subset J_2}^{\semi} \mid \I_{L \subset J_1}^{\semi}}$ is the restriction of $\Av^{\sL^{--} U_{J_1}}_! \circ \Av^{\sL^{++} U_{J_2}^-}_*$. Now we have $U_{J_1} = (B_{J_2} \cap U_{J_1}) \ltimes U_{J_2}$, hence $\sL^{--}U_{J_1} = (\sL^{--} (B_{J_2} \cap U_{J_1})) \ltimes (\sL^{--} U_{J_2})$. Since $\Av_!^{\sL^{--} U_{J_2}} \circ \Av^{\sL^{++} U_{J_2}^-}_* \simeq \id$ on $\Shv \left( \I^{\semi}_{L \subset J_2} \backslash X \right)$, it follows that $\Ras^{\I_{L \subset J_2}^{\semi} \mid \I_{L \subset J_1}^{\semi}}$ is the restriction of $\Av_!^{\sL^{--} (B_{J_2} \cap U_{J_1})}$. Using this one sees that
\[
\Ras^{\I_{L \subset J_2}^{\semi} \mid \I_{L \subset J_1}^{\semi}} = \Av_!^{\I^{\semi}_{L \subset J_1} | \I_{L \subset J_1} \ltimes \sL^+(B_{J_2} \cap U_{J_1}) \ltimes \sL U_{J_2}} \circ \oblv^{\I^{\semi}_{L \subset J_2} | \I_{L \subset J_1} \ltimes \sL^+(B_{J_2} \cap U_{J_1}) \ltimes \sL U_{J_2}}.
\]
\end{enumerate}
\end{rmk}


The following lemma is clear. 

\begin{lem}
\label{lem:Av-Ras}
Let $L \subset J\subset S$, and let $K \subset S_\af$ and $A \subset S_\af$ be finitary. 
\begin{enumerate}
\item
\label{it:Av-Ras-1}
There are canonical isomorphisms
\[
\Ras^{\I_{L \subset J}^\semi \mid \I_L} \circ (p_K)^\dag \simeq (p_K)^\dag \circ \Ras^{\I_{L \subset J}^\semi \mid \I_L}, \quad
\Ras^{\I_{L \subset J}^\semi \mid \I_L} \circ (p_K)_* \simeq (p_K)_* \circ \Ras^{\I_{L \subset J}^\semi \mid \I_L},
\]
and similarly for $q^!$, $q_*$ and $q_!$.
\item 
\label{it:Av-Ras-2}
Assume that $p>0$ and that $\k$ contains a nontrivial $p$-th root of unity. There is a canonical isomorphism
\[
\Ras^{\I_{L \subset J}^{\semi} \mid \I_L} \circ \Av^A_{\chi} \simeq \Av^A_{\chi} \circ \Ras^{\I_{L \subset J}^{\semi} \mid \I_L}, 
\]
and similarly for $\Av^A_*$ and $\Av^A_!$.
\end{enumerate}
\end{lem}

\subsection{\texorpdfstring{$J$}{J}-Wakimoto sheaves}
\label{sss:Wakimoto}

Let $J \subset S$, and
recall the subset $W_\ex^{+,J} \subset W_\ex$ defined in~\S\ref{sss:quar-length}. For $x\in W_\ex$, 
we set
\[
\Xd_{J,x}^{0,+} = \{ \mu\in \Xd_J^0\cap \Xd_{J^{\comp}}^+ \mid t_\mu x\in W_\ex^{+,J} \}.
\]
Using~\cite[Lemma~4.1]{adr} one sees that the object
\[
\nabla_{t_{-\theta}} \star^\I \Delta_{t_\theta x} \quad \in \Shv(\I\backslash \Fl)
\]
does not depend on the choice of element $\theta \in \Xd_{J,x}^{0,+}$, up to canonical isomorphism. This complex will be denoted $\sW^J_x$, and called the \textit{$J$-Wakimoto sheaf} associated with $x$. If $J=\varnothing$, the object $\sW_x := \sW^{\varnothing}_x$ is the usual Wakimoto sheaf associated with $x$ (with the conventions of~\cite{adr}). In case $J=S$, one can easily see that $\sW_x^S = \Delta_x$ for any $x \in W_\ex$. 

If $\theta \in \Xd_J^0$, it is clear that $\sW^J_{t_\theta} = \sW_{t_\theta}$, and it is not difficult to check that for any $x \in W_\ex$ we have
\begin{equation}
\label{eqn:convolution-Waki}
\sW_{t_\theta} \star^\I \sW^J_x \simeq \sW^J_{t_\theta x}.
\end{equation}
Using~\cite[Lemma~4.1(3)]{adr} one sees that $\sW^J_x$ is a perverse sheaf, supported on $\overline{\I \dot{x} \I/\I}$.

When $p>0$ and $\k$ contains a nontrivial $p$-th root of unity, for $x\in W_\ex^{(A,J),\max}$ we also set
\begin{equation*}
\sW^J_{\chi_A,x}:= \Av^A_{\chi}(q^! \sW^J_x) \quad \in \Shv(\I \backslash \sL G / (\Iu^A, \chi_A)).
\end{equation*}
Note that if $\theta \in \Xd_J^0 \cap \Xd_{J^\comp}^+$ satisfies $t_\theta x \in W_\ex^{+,J} \cap W_\ex^{A,\max}$, then by~\eqref{eqn:Av-chi-D-classical} we have
\[
\sW^J_{\chi_A,x} = \nabla_{t_{-\theta}} \star^\I \Delta_{\chi_A,t_\theta x}.
\]
Here again, $\sW^J_{\chi_A,x}$ is a perverse sheaf supported on the closure of the orbit parame\-trized by $x$.

\begin{lem}
\label{lem:Av-chi-Wak}
Assume that $p>0$ and that $\k$ contains a nontrivial $p$-th root of unity. Let $J \subset S$, and let $A \subset S_\af$ be finitary.
\begin{enumerate}
\item
\label{it:Av-chi-Wak-1}
Let $x \in W_\ex$, and let $w \in W_A$ be the unique element such that $xw \in W_\ex^{(A,J), \max}$. Then we have
\[
\Av^A_{\chi}(q^! \sW_x^J) \simeq \sW_{\chi_A,xw}^J.
\]
\item
\label{it:Av-chi-Wak-2}
Let $x \in W_\ex^{(A,J),\max}$. The object $\Av^{A}_{!}(\sW^{J}_{\chi_A, x})$ is an iterated extension of objects $q^! \sW^J_{xw}$ with $w\in W_A$, each appearing once.
\end{enumerate}
\end{lem}

\begin{proof}
\eqref{it:Av-chi-Wak-1}
One can choose $\theta \in \Xd_J^0 \cap \Xd_{J^\comp}^+$ such that $t_\theta x \in W_\ex^{+,J}$ and $t_\theta xw \in W_\ex^{+,J} \cap W_\ex^{A,\max}$. Then we have $\sW_x^J \simeq \nabla_{t_{-\theta}} \star^\I \Delta_{t_\theta x}$, $\sW^J_{\chi_A,xw} = \nabla_{t_{-\theta}} \star^\I \Delta_{\chi_A,t_\theta xw}$, and the claim follows from~\eqref{eqn:Av-chi-D-classical}.

\eqref{it:Av-chi-Wak-2}
The proof is similar, using a variant of~\cite[Lemma~3.6]{ar-model} (for $\Fl$ instead of $\Gr$, and with left and right multiplications switched).
\end{proof}

The following statement is analogous to~\cite[Lemma~5.5]{adr}, and follows from similar arguments. (See also the discussion in~\cite[\S 3.1.10]{DL23}.)

\begin{lem}
\label{lem:Ras-shift-Wak}
Let $L \subset J \subset S$, and let $X$ be as above.
For any $\theta \in \Xd_J^0$, we have a canonical equivalence
\[
\Ras^{\I^\semi_{L \subset J} \mid \I_L} ( - \langle \theta \rangle) \simeq ((p_L)_* \sW_{t_\theta}) \star^{\I_L} \Ras^{\I^\semi_{L \subset J} \mid \I_L}(-).
\]
\end{lem}


We deduce the following result, which generalizes~\cite[Proposition~5.6]{adr}.

\begin{prop}
\label{prop:Ras-D}
Let $J \subset S$.
\begin{enumerate}
\item
\label{it:Ras-D-1}
For any $x \in W_\ex$ we have an equivalence 
\[
\Ras^{\I_J^{\quar}\mid \I} \left( \Delta^{J,\quar}_x \right) \simeq \sW^J_x.
\]
\item
\label{it:Ras-D-2}
Assume that $p>0$ and that $\k$ contains a nontrivial $p$-th root of unity, and let $A \subset S_\af$ be finitary. For any $x \in W_\ex^{(A,J),\max}$ we have an equivalence
\[
\Ras^{\I_J^{\quar}\mid \I} \left( \Delta^{J,\quar}_{\chi_A,x} \right) \simeq \sW^J_{\chi_A,x}. 
\]
\end{enumerate}
\end{prop}

\begin{proof}
\eqref{it:Ras-D-1}
Using Lemma~\ref{lem:Ras-shift-Wak} together with~\eqref{eqn:D-N-shift} and~\eqref{eqn:convolution-Waki}, we can assume that $x \in W_\ex^{+,J}$. In this case, what we need to prove is that $\Ras^{\I \mid \I_J^{\quar}} ( \Delta_x) \simeq \Delta^{J,\quar}_x$.
If $\theta$ is as in~\eqref{eq:Raskin-formula}, for any $n \geq 0$ we have $\ell(t_{n\theta} x) = \ell(x)+\ell(t_{n\theta})$, so
\[
\Delta_{t_{n\theta}} \star^{\I} \Delta_x \simeq \Delta_{t_{n\theta} x}
\]
by~\cite[Lemma~4.1(1)]{adr}.
The desired claim follows using~\eqref{eq:Raskin-formula}, since for $n \gg 0$ we have
\[
z^{-n\theta}\I z^{n\theta} \dot{x} \I/\I = \I^{\quar}_{J,(n)} \dot{x} \I/\I.
\]

\eqref{it:Ras-D-2}
The proof is similar.
\end{proof}

We now prove a generalization of (the case $L=\varnothing$ of) Lemma~\ref{lem:Ras-shift-Wak}.

\begin{lem}
\label{lem:Ras-convolution-Wak}
For any $\lambda\in \Xd$, we have an equivalence of functors from $\Shv(\I^{\quar}_J\backslash X)$ to $\Shv(\I\backslash X)$, 
\begin{equation*}
\Ras^{\I_J^{\quar}\mid \I} \circ (\sW^{M_J}_{t_\lambda} \star^{\I_{\varnothing \subset J}} -) [-\langle \lambda, 2\rho-2\rho_J \rangle]
\simeq \sW_{t_\lambda} \star^{\I} \Ras^{\I_J^{\quar}\mid \I} (-),
\end{equation*}
where $\sW^{M_J}_{t_\lambda}$ is the Wakimoto sheaf in $\Shv(\I_{\varnothing \subset J}\backslash \sL{M}_J/\I_{\varnothing \subset J})$ attached to $\lambda$. 
\end{lem}

\begin{proof}
Since every coweight is a difference of two dominant coweights, and both sides are compatible with convolution in the appropriate way, invertibility of the Wakimoto sheaves in $\Shv(\I\backslash \sL{G}/\I)$ and $\Shv(\I_{\varnothing \subset J}\backslash \sL{M}_J/\I_{\varnothing \subset J})$ reduces the proof to the case $\lambda$ is dominant.
In this case $\sW_{t_\lambda}=\Delta_{t_\lambda}=\sW_{t_\lambda}^J$, and $\sW^{M_J}_{t_\lambda}$ is the standard object $\Delta^{M_J}_{t_\lambda}$ in $\Shv(\I_{\varnothing \subset J}\backslash \sL{M}_J/\I_{\varnothing \subset J})$ attached to $\lambda$. 
Once we have fixed an element $\theta$ as in~\eqref{eq:Raskin-formula},
by Proposition~\ref{prop:Ras-D} and~\eqref{eq:Raskin-formula} it is enough to show that there is an equivalence 
\[
\colim_{n} \, \Delta^{M_J}_{t_\lambda} \star^{\I_{\varnothing \subset J}} z^{-n\theta} \cdot \Delta_{t_{n\theta}} [n\ell(t_\theta) -\langle \lambda, 2\rho-2\rho_J \rangle] \simeq \Delta^{J,\quar}_{t_\lambda} . 
\]
This follows from the isomorphism
\[
 \I_{\varnothing \subset J}z^\lambda \I_{\varnothing \subset J} \times^{\I_{\varnothing \subset J}} z^{-n\theta} \I z^{n\theta}\I/\I 
\xs \big(\I_{\varnothing \subset J}\ltimes \Ad_{z^{\lambda-n\theta}}(\sL^+U_{J})\big) z^{\lambda} \I/\I
\]
induced by multiplication for $n \gg 0$.
\end{proof}

\subsection{Properties of averaging functors} 
\label{sss:prop-av}

As an application of Raskin's equivalence, we obtain some properties of the averaging functors defined in \S\ref{sss:av-fun}. We assume that $p>0$ and that $\k$ contains a nontrivial $p$-th root of unity.

%
%

\begin{lem}
\label{lem:Av-chi-D}
Let $L \subset J \subset S$, and let $A \subset S_\af$ be finitary.

\begin{enumerate}
\item
\label{it:Av-chi-D-1}
For any $x\in {}^L W_\ex$ we have 
\[
\Av^{A}_{\chi}(q^! (\Delta^{L \subset J, \semi}_{x})) \simeq 
\begin{cases}
\Delta^{L \subset J,\semi}_{\chi_A, xw} & \text{if $w\in W_A$ satisfies $xw\in {}^L W_\ex^{(A,J),\max}$;} \\
0 & \text{if $x W_A \cap {}^L W_\ex^{(A,J),\max} = \varnothing$.} 
\end{cases}  
\]
\item 
\label{it:Av-chi-D-2}
For any $x\in {}^L W_\ex^{(A,J),\max}$, the object $\Av^{A}_{!}(\Delta^{L \subset J,\semi}_{\chi_A,x})$
is an iterated extension of objects $q^! (\Delta^{L \subset J, \semi}_{xw})$ with $w\in W_A$, 
each appearing once.
\end{enumerate}
\end{lem}

\begin{proof}
\eqref{it:Av-chi-D-1}
First, assume that $L=\varnothing$. (In this case, an element $w$ as in the statement always exists.) Then the claim follows from
Lemma~\ref{lem:Av-Ras}, Lemma~\ref{lem:Av-chi-Wak}\eqref{it:Av-chi-Wak-1} and Proposition~\ref{prop:Ras-D}. The general case follows, using Lemma~\ref{lem:Av-quar-semi-D-N}\eqref{it:Av-quar-semi-D-N-1} (applied for the subset $K=\varnothing$), the commutation of $\Av_!^{\I_{L \subset J}^\semi \mid \I_J^\quar}$ with $q^!$ and $\Av^{A}_{\chi}$, and Lemma~\ref{lem:av-oblv-whit}\eqref{it:av-oblv-whit-1}.

\eqref{it:Av-chi-D-2}
The proof is similar, using Lemma~\ref{lem:Av-chi-Wak} for the case $L=\varnothing$, and then Lem\-ma~\ref{lem:av-oblv-whit}\eqref{it:av-oblv-whit-1} and Lemma~\ref{lem:Av-quar-semi-D-N}\eqref{it:Av-quar-semi-D-N-1} (with $K = \varnothing$).
\end{proof}

\section{Perverse t-structure} 
\label{sec:perv}

Let $X$ be either $\Fl_K$ for some finitary $K \subset S_\af$, $\tFl$, or the symbol $\sL G/(\Iu^A,\chi_A)$ for some finitary $A \subset S_\af$ (assuming that $p>0$ and that $\k$ contains a nontrivial $p$-th root of unity).

\subsection{Definition} 

If $L \subset J \subset S$, we define a t-structure, called the \textit{perverse t-structure}, on $\Shv\left(\I_{L\subset J}^{\semi}\backslash X\right)$, such that 
\[
\Shv\left(\I_{L\subset J}^{\semi}\backslash X\right)^{\leq 0}  
\]
is the smallest full subcategory that is stable under extensions and colimits and contains the following collection of objects: 
\begin{itemize}
\item $\left( \Delta^{L\subset J,\semi}_{K,x} : x\in {}^LW_\ex^{(K,J)} \right)$ when $X=\Fl_K$; 
\item $\left( q^!\Delta^{L\subset J,\semi}_x : x\in {}^LW_\ex \right)$ when $X=\tFl$; 
\item $\left( \Delta^{L\subset J,\semi}_{\chi_A,x} : x\in {}^LW_\ex^{(A,J),\max} \right)$ when $X = \sL G/(\Iu^A,\chi_A)$.
\end{itemize}
Note that when $J=S$ (so that $\I_{L \subset S}^\semi=\I_L$) and $X=\Fl_K$, this t-structure is the natural one in the sense that the forgetful functor $\oblv^{\I_L} : \Shv(\I_L \backslash \Fl_K) \to \Shv(\Fl_K)$ is t-exact with respect to our perverse t-structure on the domain and the usual perverse t-structure (of middle perversity) on the codomain.

As in~\cite[\S 6.1]{adr}, the nonnegative part $\Shv\left(\I_{L\subset J}^{\semi}\backslash X\right)^{\geq 0}$ can be described as the right orthogonal to $\Shv\left(\I_{L\subset J}^{\semi}\backslash X\right)^{\leq -1}$, or equivalently in terms of $!$-restriction to strata.

As usual, we will say that an object of $\Shv\left(\I_{L\subset J}^{\semi}\backslash X\right)$ \emph{is perverse} to say that it belongs to the heart of this t-structure.

\subsection{Some t-exactness results}

We continue with $X$ as above.

\begin{lem}
\label{lem:oblv-semi-quar-texact}
Let $L_1\subset L_2\subset J\subset S$. 
The functor 
\[
\oblv^{\I_{L_2 \subset J}^{\semi}|\I_{L_1 \subset J}^{\semi}}:
\Shv\left(\I_{L_2 \subset J}^{\semi}\backslash X\right)\rightarrow \Shv\left(\I_{L_1 \subset J}^{\semi}\backslash X\right)
\]
is t-exact. 
\end{lem}

\begin{proof}
We explain the proof in case $X=\Fl_K$; the other cases are similar. (For $X=\sL G/(\Iu^A,\chi_A)$, the argument uses Lemma~\ref{lem:av-oblv-whit}.)
We first assume that $L_1=\varnothing$. 
To show right t-exactness of our functor, it is enough to show that $\oblv^{\I_{L_2 \subset J}^{\semi}|\I_{J}^{\quar}}(\Delta^{L_2\subset J,\semi}_{K,x})$ is contained in $\Shv\left(\I_{J}^{\quar}\backslash X\right)^{\leq 0}$ for any $x\in {}^{L_2}W_\ex^{(K,J)}$, which follows from Lemma~\ref{lem:Av-quar-semi-D-N}\eqref{it:Av-quar-semi-D-N-2}. 
We next show left t-exactness. If $\sF\in \Shv\left(\I_{L_2 \subset J}^{\semi}\backslash X\right)^{\geq 0}$, then for any $x\in W_\ex^{(K,J)}$ and $n> 0$ we have 
\begin{equation*}
\Hom \left( \Delta^{J,\quar}_{K,x}[n], \oblv^{\I_{L_2 \subset J}^{\semi}|\I_{J}^{\quar}}(\sF) \right) \simeq \Hom \left( \Av_!^{\I_{L_2 \subset J}^{\semi}|\I_{J}^{\quar}}(\Delta^{J,\quar}_{K,x})[n],\sF \right). 
\end{equation*}
By Lemma~\ref{lem:Av-quar-semi-D-N}\eqref{it:Av-quar-semi-D-N-1}, the right-hand side
vanishes, which implies that $\oblv^{\I_{L_2 \subset J}^{\semi}|\I_{J}^{\quar}}(\sF)\in \Shv\left(\I_J^{\quar}\backslash X\right)^{\geq 0}$ as desired.

Now let $L_1$ be general. 
It is clear that we have a natural equivalence of functors 
\[
\oblv^{\I_{L_2 \subset J}^{\semi}|\I_{J}^{\quar}} \simeq 
\oblv^{\I_{L_1 \subset J}^{\semi}|\I_{J}^{\quar}} \circ \oblv^{\I_{L_2 \subset J}^{\semi}|\I_{L_1 \subset J}^{\semi}}. 
\]
For any $L\subset J$ the functor $\oblv^{\I_{L \subset J}^{\semi}|\I_{J}^{\quar}}$ is conservative, and is t-exact by our previous discussion. 
Hence the functor $\oblv^{\I_{L_2 \subset J}^{\semi}|\I_{L_1 \subset J}^{\semi}}$ is t-exact. 
\end{proof}

\begin{lem}
\label{lem:pullback-Av-texact}
Let $L\subset J\subset S$, and let $K\subset S_\af$ and $A \subset S_\af$ be finitary. 
\begin{enumerate}
\item 
\label{it:pullback-q-texact}
The functor $q^!: \Shv\left(\I_{L\subset J}^{\semi}\backslash \Fl\right)\rightarrow \Shv\left(\I_{L\subset J}^{\semi}\backslash \tFl\right)$ is t-exact. 
Moreover, its restriction to the hearts of the perverse t-structures is fully faithful, and its essential image is stable under subquotients. 
\item
\label{it:pullback-p-texact}
The functor $(p_K)^\dag:\Shv\left(\I_{L\subset J}^{\semi}\backslash \Fl_K\right)\rightarrow \Shv\left(\I_{L\subset J}^{\semi}\backslash \Fl\right)$ is t-exact. 
Moreover, its restriction to the hearts of the perverse t-structures is fully faithful, and its essential image is stable under subquotients and extensions. 
\item
\label{it:Av-texact}
If $p>0$ and $\k$ contains a nontrivial $p$-th root of unity,
the functor $\Av^A_{\chi}:\Shv\left(\I_{L\subset J}^{\semi}\backslash \tFl\right)\rightarrow \Shv\left(\I_{L\subset J}^{\semi}\backslash \sL G/(\Iu^A,\chi_A)\right)$ is t-exact. 
\end{enumerate}
\end{lem}

\begin{proof}
The proofs of~\eqref{it:pullback-q-texact} and~\eqref{it:pullback-p-texact} are similar to those of~\cite[Lemmas 6.6 \& 6.8]{adr}. (In particular, the proof of~\eqref{it:pullback-p-texact} uses Lemma~\ref{lem:push-pull-pK-D-N-semi}).  
The statement in~\eqref{it:Av-texact} follows from similar arguments as those used in the proof of Lemma~\ref{lem:oblv-semi-quar-texact}, using Lemma~\ref{lem:Av-chi-D}. 
\end{proof}

\subsection{Perversity of standard sheaves} 

We now establish an analogue of~\cite[Lem\-ma~6.9]{adr}, proving that the standard objects are in the heart of the perverse t-structure in some cases. (This property is not satisfied in general in the other cases; e.g.~in $\Shv(\I_L \backslash \Fl_K)$ with $L$ and $K$ nonempty the orbits are not necessarily affine, so that the standard objects are not perverse in general.)

\begin{lem}
\label{lem:standard-perverse}
Let $L\subset J\subset S$, and let $K\subset S_\af$ and $A \subset S_\af$ be finitary. 
The objects $\Delta^{J,\quar}_{K,x}$ for $x \in W_\ex^{(K,J)}$, $\Delta^{L\subset J,\semi}_{x}$ for $x \in {}^L W_\ex$, and $\Delta^{L\subset J,\semi}_{\chi_A,x}$ for $x \in {}^L W_\ex^{(A,J),\max}$ (when $p>0$ and $\k$ contains a nontrivial $p$-th root of unity) are contained in the heart of the perverse t-structure on the corresponding $\infty$-categories. 
\end{lem}

\begin{proof}
First, we consider the case of the objects $\Delta^{J,\quar}_{K,x}$ ($x \in W_\ex^{(K,J)}$). It is enough to show that these objects belong to the ``$\geq 0$'' part of the t-structure, i.e.~that for any $x,y\in W_\ex^{(K,J)}$ and $n>0$ we have 
\begin{equation*}
\Hom(\Delta^{J,\quar}_{K,x}[n],\Delta^{J,\quar}_{K,y})=0. 
\end{equation*}

Let $x \in W_\ex^{(K,J)}$.
By Proposition~\ref{prop:Ras-D}\eqref{it:Ras-D-1}, we have $\Ras^{\I_{J}^{\quar}\mid \I}(\Delta^{J,\quar}_{x})\simeq \sW^{J}_{x}$. Using Lemma~\ref{lem:push-pull-pK-D-N-semi}\eqref{it:push-pull-pK-D-N-semi-1} and Lemma~\ref{lem:Av-Ras}\eqref{it:Av-Ras-1} we deduce that
\[
\Ras^{\I_{J}^{\quar}\mid \I}(\Delta^{J,\quar}_{K,x})\simeq \Ras^{\I_{J}^{\quar}\mid \I}((p_K)_* \Delta^{J,\quar}_{x}) \simeq (p_K)_* \sW^{J}_{x} \simeq \nabla_{t_{-\theta}} \star^{\I} \Delta_{K,t_\theta x}
\]
where $\theta \in \Xd_J^0 \cap \Xd_{J^\comp}^+$ satisfies $t_\theta x \in W_\ex^{+,J} \cap W_\ex^{K}$. If now we consider another element $y \in W_\ex^{(K,J)}$, we can choose $\theta \in \Xd_J^0 \cap \Xd_{J^\comp}^+$ such that $t_\theta x$ and $t_\theta y$ are both in $W_\ex^{+,J} \cap W_\ex^{K}$, and we obtain that
\begin{multline*}
\Hom(\Delta^{J,\quar}_{K,x}[n],\Delta^{J,\quar}_{K,y}) \simeq \Hom(\nabla_{t_{-\theta}} \star^{\I} \Delta_{K,t_\theta x} [n],\nabla_{t_{-\theta}} \star^{\I} \Delta_{K,t_\theta y}) \\
\simeq \Hom(\Delta_{K,t_\theta x} [n], \Delta_{K,t_\theta y})
\end{multline*}
by invertibility of the object $\nabla_{t_{-\theta}}$. We deduce the desired vanishing, since the objects $\Delta_{K,t_\theta x}$ and $\Delta_{K,t_\theta y}$ are perverse.

Now, we treat the objects $\Delta^{L\subset J,\semi}_{x}$ ($x \in {}^L W_\ex$). Here again, what we have to show is that
\[
\Hom(\Delta^{L \subset J,\semi}_{x}[n],\Delta^{L \subset J,\semi}_{y})=0
\]
for $x,y \in {}^L W_\ex$ and $n>0$.

Let $x \in {}^L W_\ex$. By Lemma~\ref{lem:Av-quar-semi-D-N}\eqref{it:Av-quar-semi-D-N-1}, the commutative diagram obtained from~\eqref{eq:comm-diag-ras-oblv} by passing to left adjoints, and Proposition~\ref{prop:Ras-D}\eqref{it:Ras-D-1},
we have equivalences 
\begin{multline*}
\Ras^{\I_{L\subset J}^{\semi}\mid \I_L}(\Delta^{L\subset J,\semi}_{x}) \simeq \Ras^{\I_{L\subset J}^{\semi}\mid \I_L}(\Av_!^{\I^\semi_{L \subset J} \mid \I_J^\quar}(\Delta^{J,\quar}_{x}))[-\ell(w_L)] \\
\simeq \Av^{\I_L\mid \I}_! (\sW^{J}_{x})[-\ell(w_L)] \simeq \Av^{\I_L\mid \I}_!(\nabla_{t_{-\theta}}) \star^\I \Delta_{t_\theta x} [-\ell(w_L)] 
\end{multline*}
where $\theta \in \Xd_J^0 \cap \Xd_{J^\comp}^+$ satisfies $t_\theta x \in W_\ex^{+,J}$.
Now, using Remark~\ref{rmk:Raskin}\eqref{it:ttheta-minimal} one sees that we have canonical isomorphisms
\[
\I_L \times^{\I} \I z^\theta \I / \I \simeq \I_L z^\theta \I/\I = \I_L z^\theta \I_L / \I \simeq \I_L z^\theta \I_L \times^{\I_L} \I_L/\I,
\]
from which we obtain an equivalence
\begin{equation}
\label{eqn:computation-standard-perverse}
\Av^{\I_L\mid \I}_!(\nabla_{t_{-\theta}}) \star^\I \Delta_{t_\theta x} [-\ell(w_L)] \simeq \nabla^L_{L,t_{-\theta}} \star^{\I_L} \Av^{\I_L\mid \I}_!(\Delta_{t_\theta x}) [-\ell(w_L)].
\end{equation}
Using once again Lemma~\ref{lem:Av-quar-semi-D-N}\eqref{it:Av-quar-semi-D-N-1} (for the subset ``$J$'' equal to $S$) we deduce an equivalence
\[
\Av^{\I_L\mid \I}_!(\nabla_{t_{-\theta}}) \star^\I \Delta_{t_\theta x} [-\ell(w_L)] \simeq \nabla^L_{L,t_{-\theta}} \star^{\I_L} \Delta^L_{t_\theta x}.
\]
(Note that here we have $t_\theta x \in {}^L W_\ex$ by Lemma~\ref{lem:min-length}\eqref{it:min-length-0}.)
If now $y \in {}^L W_\ex$ is another element, and if $\theta \in \Xd_J^0 \cap \Xd_{J^\comp}^+$ is such that $t_\theta x$ and $t_\theta y$ are both in $W_\ex^{+,J}$, then we obtain that 
\[
\Hom(\Delta^{L \subset J,\semi}_{x}[n],\Delta^{L \subset J,\semi}_{y})=\Hom(\nabla^L_{L,t_{-\theta}} \star^{\I_L} \Delta^L_{t_\theta x}[n],\nabla^L_{L,t_{-\theta}} \star^{\I_L} \Delta^L_{t_\theta y}).
\]
As in~\cite[Lemma~3.1.5(i)]{DL23}, the object $\nabla^L_{L,t_{-\theta}}$ is invertible in $\Shv(\I_L \backslash \sL G / \I_L)$, so that
\[
\Hom(\nabla^L_{L,t_{-\theta}} \star^{\I_L} \Delta^L_{t_\theta x}[n],\nabla^L_{L,t_{-\theta}} \star^{\I_L} \Delta^L_{t_\theta y}) \simeq \Hom(\Delta^L_{t_\theta x}[n],\Delta^L_{t_\theta y}).
\]
One then concludes as above using the fact that $\Delta^L_{t_\theta x}$ and $\Delta^L_{t_\theta y}$ are perverse.

Finally, the case of the objects $\Delta^{L\subset J,\semi}_{\chi_A,x}$ ($x \in {}^L W_\ex^{(A,J),\max}$) follows from the preceding case, using Lemma~\ref{lem:pullback-Av-texact}\eqref{it:pullback-q-texact}--\eqref{it:Av-texact} and Lemma~\ref{lem:Av-chi-D}\eqref{it:Av-chi-D-1}.
\end{proof}

\subsection{Finitely generated complexes}

\subsubsection{Restriction of the t-structure to finitely generated complexes} 
\label{sss:restriction-fg}

For $L\subset J\subset S$ and $X$ as above,
let us now consider the small stable subcategory
\[
\Shvfg(\I_L\backslash X) \subset \Shv(\I_L\backslash X)
\]
of ``finitely generated'' complexes, see the discussions in~\cite[\S 2.2 \& \S 2.5]{adr}. (The homotopy category of $\Shvfg(\I_L\backslash X)$ is what is usually called the $\I_L$-equivariant derived category of constructible sheaves on $X$, see e.g.~\cite{ARbook, ab, by}.)
We will denote by 
\[
\Shv_{\omega}\left(\I_{L\subset J}^{\semi}\backslash X\right) \subset \Shv\left(\I_{L\subset J}^{\semi}\backslash X\right)
\]
the image of $\Shvfg(\I_L\backslash X)$ under $\Ras^{\I_L | \I^\semi_{L \subset J}}$. 
We think of this subcategory as consisting of the ``constructible'' objects.
As in~\cite[\S 6.4]{adr}, the objects in $\Shv_{\omega}\left(\I_{L\subset J}^{\semi}\backslash X\right)$ admit the following characterization.

\begin{lem}
\label{lem:constr-complx}
An object $\sF$ of $\Shv\left(\I_{L\subset J}^{\semi}\backslash \Fl_K\right)$ is contained in $\Shv_{\omega}\left(\I_{L\subset J}^{\semi}\backslash \Fl_K\right)$ if and only if 
\begin{enumerate}
\item $\sF$ is supported on a finite union of closures of orbits $\fS^{L\subset J}_{K,x}$; 
\item the set $\{x\in {}^LW_\ex^{(K,J)}\mid (\bs^{L\subset J}_{K,x})^* \sF\neq 0\}$ is finite; 
\item for any $x\in {}^LW_\ex^{(K,J)}$, the image of $(\bs^{L\subset J}_{K,x})^* \sF$ under the equivalence~\eqref{eq:semi-orbt} belongs to $\Shvc \left( \mathrm{Stab}_{M_L}(\dot{x}\I_K/\I_K)\backslash \pt \right)$. 
\end{enumerate}
We have a similar description of $\Shv_{\omega}\left(\I_{L\subset J}^{\semi} \backslash X\right)$ for the other choices of $X$. 
\end{lem}



The following statement generalizes~\cite[Lemma~6.10]{adr}.

\begin{lem}
\label{lem:tstr-fg-objects}
The perverse t-structure on $\Shv\left(\I_{L\subset J}^{\semi} \backslash X\right)$ restricts to a t-structure on $\Shv_\omega\left(\I_{L\subset J}^{\semi} \backslash X\right)$. 
\end{lem}

\begin{proof}
In the cases where Lemma~\ref{lem:standard-perverse} applies and for $\Shv \left( \I^\semi_{L \subset J} \backslash \tFl \right)$, the claim follows from the same arguments as for the proof of~\cite[Lemma~6.10]{adr}. The cases that remain are those of the $\infty$-categories $\Shv\left(\I_{L\subset J}^{\semi} \backslash \Fl_K \right)$ when $L$ and $K$ are nonempty. In this case we have to show that for any $\sF$ in $\Shv_\omega\left(\I_{L\subset J}^{\semi} \backslash \Fl_K \right)$ and $n \in \Z$, the complex $\tau_{<n}(\sF)$ belongs to $\Shv_\omega\left(\I_{L\subset J}^{\semi} \backslash \Fl_K \right)$. This follows from the corresponding claim in $\Shv_\omega\left(\I_{J}^{\quar} \backslash \Fl_K \right)$ and the t-exactness of the functor $\oblv^{\I_{L\subset J}^{\semi} \mid \I_J^\quar}$ proved in Lemma~\ref{lem:oblv-semi-quar-texact}.
\end{proof}

\subsubsection{Description of the t-structure on finitely generated complexes} 

Let $X$ be as above and $J \subset S$.
Thanks to Lem\-ma~\ref{lem:tstr-fg-objects} we have a perverse t-structure on the $\infty$-category $\Shv_\omega\left(\I_{J}^{\quar} \backslash X\right)$.
The following proposition describes this t-structure in terms of the standard perverse t-structure on $\Shv(\I\backslash X)$.

\begin{prop}
\label{prop:tstr-convolution-Wak}
For any $\sF\in \Shvfg(\I\backslash X)$, the following conditions are equivalent:
\begin{enumerate}
\item 
$\Ras^{\I\mid \I_J^{\quar}}(\sF) \in \Shv\left(\I^{\quar}_J\backslash X\right)^{\leq 0}$; 
\item 
for any $\lambda\in \Xd_J^0$,
$\sW_{t_\lambda}\star^\I \sF\in \Shv(\I\backslash X)^{\leq 0}$; 
\item there exists $\nu\in \Xd_J^0 \cap\Xd^{++}_{J^\comp}$ such that $\sW_{t_{n\nu}} \star^\I \sF\in \Shv(\I\backslash X)^{\leq 0}$ for any integer $n >0$.
\end{enumerate}
\end{prop}

\begin{proof}
The proof is similar to that of~\cite[Lemma~6.11]{adr}. 
\end{proof}

\begin{corollary}
\label{cor:Ras-right-texact}
For any $L\subset J_1\subset J_2\subset S$, the equivalence $\Ras^{\I_{L\subset J_1}^{\semi}\mid \I_{L\subset J_2}^{\semi}}$ is right t-exact. 
\end{corollary}

\begin{proof}
By \eqref{eq:comm-diag-ras-oblv} and Lemma~\ref{lem:oblv-semi-quar-texact}, together with the fact that the functor $\oblv^{\I^\semi_{L\subset J}\mid \I^\quar_J}$ is conservative for any $L\subset J \subset S$, it suffices to prove the claim when $L=\varnothing$. 
We explain the proof in the case when $X=\Fl$; the other cases are similar. 
What we have to show is that $\Ras^{\I_{J_1}^{\quar}\mid \I_{J_2}^{\quar}}(\Delta^{J_1,\quar}_x)$ is contained in $\Shv \left( \I_{J_2}^{\quar}\backslash \Fl \right)^{\leq 0}$ for any $x\in W_\ex$. 
Since $\Delta^{J_1,\quar}_x\in \Shv_{\omega} \left( \I_{J_1}^{\quar}\backslash \Fl \right)$, this property follows from the criterion in Proposition~\ref{prop:tstr-convolution-Wak} since $\Xd_{J_2}^0\subset \Xd_{J_1}^0$.
\end{proof}

Let us also note the following consequence of Proposition~\ref{prop:tstr-convolution-Wak}, which will be used in Section~\ref{sec:t-exactness}. Recall, for $L \subset J \subset S$, the right action of $\Shv(\I \backslash \Fl)$ on $\Shv(\I_{L \subset J}^{\semi} \backslash \Fl)$, see~\S\ref{ss:convolution}.

\begin{corollary}
\label{cor:right-t-exactness-convolution-nabla}
Let $L \subset J \subset S$.
For any $w \in W_\ex$, the functor 
\[
(-) \star^\I \nabla_w : \Shv\left(\I^{\semi}_{L \subset J} \backslash \Fl\right) \to \Shv\left(\I^{\semi}_{L \subset J}\backslash \Fl\right)
\]
is right t-exact for the perverse t-structure.
\end{corollary}

\begin{proof}
Since the functor $\oblv^{\I^{\semi}_{L \subset J} \mid \I^{\quar}_J}$ is conservative and t-exact (see Lemma~\ref{lem:oblv-semi-quar-texact}) and commutes with $(-) \star^\I \nabla_w$, it suffices to treat the case $L=\varnothing$. For this it suffices to show the right t-exactness of the functor
\[
(-) \star^\I \nabla_w : \Shv_\omega\left(\I^{\quar}_J\backslash \Fl\right) \to \Shv_\omega\left(\I^{\quar}_J\backslash \Fl\right).
\]
This property follows from Proposition~\ref{prop:tstr-convolution-Wak}, the commutation of $\Ras^{\I\mid \I_J^{\quar}}$ with $(-) \star^\I \nabla_w$, and the standard fact that $(-) \star^\I \nabla_w : \Shvfg(\I\backslash \Fl) \to \Shvfg(\I\backslash \Fl)$ is right t-exact. (We recall that this fact can be proved after interpreting this functor as the composition of exterior product with a perversely shifted constant sheaf with $*$-pushforward under an affine morphism, which is right t-exact by Artin's theorem, see~\cite[Th\'eor\`eme~4.1.1]{bbdg}. See also~\cite[\S 4.1.3]{ARbook} for details.)
\end{proof}

\subsection{Perversity of costandard sheaves}
\label{ss:perv-costandards}

Let $J\subset S$. We will denote by $\unlhd^J$ the preorder on $\Xd$ given by 
\[
\mu \unlhd^J \lambda \quad 
\text{if and only if} \quad 
\lambda-\mu \in \Xd_J^0\cap \Xd_{J^{\comp}}^+. 
\]

Let $x\in W_\ex$, and recall the set $\Xd_{J,x}^{0,+}$ defined in~\S\ref{sss:Wakimoto}.
For $x\in W_\ex$ and $\mu\in \Xd_{J,x}^{0,+}$ we set 
\[
\sW^{J,\vee,\mu}_x:=\nabla_{t_{-\mu}} \star^{\I} \nabla_{t_\mu x} . 
\]
For any $\mu, \lambda\in \Xd_{J,x}^{0,+}$ such that $\mu \unlhd^J \lambda$, there is a canonical morphism 
\[
\sW^{J,\vee,\mu}_x \rightarrow \sW^{J,\vee,\lambda}_x 
\]
defined as follows. Since $\lambda-\mu$ is dominant, using~\cite[Lemma~4.1]{adr} we obtain a canonical isomorphism $\nabla_{t_{-\mu}} \simeq \nabla_{t_{-\lambda}} \star^{\I} \Delta_{t_{\lambda-\mu}}$. On the other hand, since $t_\mu x \in W_\ex^{+,J}$ we have $\ell(t_\lambda x) = \ell(t_{\lambda-\mu}) + \ell(t_\mu x)$, hence a canonical isomorphism $\nabla_{t_\lambda x} \simeq \nabla_{t_{\lambda-\mu}} \star^{\I} \nabla_{t_\mu x}$ (again, by~\cite[Lemma~4.1]{adr}). Then the desired map is the composition
\begin{multline*}
\sW^{J,\vee,\mu}_x 
=\nabla_{t_{-\mu}} \star^{\I} \nabla_{t_\mu x} 
\simeq  \nabla_{t_{-\lambda}} \star^{\I} \Delta_{t_{\lambda-\mu}} \star^{\I} \nabla_{t_\mu x} \\ 
\rightarrow  \nabla_{t_{-\lambda}} \star^{\I} \nabla_{t_{\lambda-\mu}} \star^{\I} \nabla_{t_\mu x}
\simeq \nabla_{t_{-\lambda}} \star^{\I} \nabla_{t_\lambda x} 
=\sW^{J,\vee,\lambda}_x ,
\end{multline*}
where the middle map is induced by the canonical morphism $\Delta_{t_{\lambda-\mu}}\rightarrow \nabla_{t_{\lambda-\mu}}$. 
As in \cite[Lemma 4.5]{adr} one checks that the assignment $\mu \mapsto \sW^{J,\vee,\mu}_x$ induces a functor (in the sense of $\infty$-categories) from $(\Xd_{J,x}^{0,+}, \unlhd^J)$ to $\Shv(\I \backslash \Fl)$, so that we can define
the \emph{dual $J$-Wakimoto object} associated with $x$ as
\[
\sW^{J,\vee}_x:= \colim_{\mu\in \Xd_{J,x}^{0,+}} \sW^{J,\vee,\mu}_x .
\]

If $p>0$ and $\k$ contains a nontrivial $p$-th root of unity, for $A \subset S_\af$ a finitary subset and $x \in W_\ex^{(A,J),\max}$ we also set $\Xd^{0,+}_{J,x,A} = \{ \mu \in \Xd^{0,+}_{J,x} \mid t_\mu x \in W_\ex^{A, \max}\}$, and define
\[
\sW^{J,\vee}_{\chi_A,x} := \colim_{\mu\in \Xd_{J,x,A}^{0,+}} \nabla_{t_{-\mu}} \star^{\I} \nabla_{\chi_A,t_\mu x}.
\]
Using an analogue of~\eqref{eqn:Av-chi-D-classical} for costandard objects, one checks that for all $x \in W_\ex^{(A,J),\max}$ we have
\begin{equation}
\label{eqn:Av-dual-Waki}
\Av^A_\chi(q^! \sW^{J,\vee}_x) \simeq \sW^{J,\vee}_{\chi_A,x}.
\end{equation}

The following is a generalization of \cite[Corollary 6.14]{adr}, which follows from similar arguments. 

\begin{lem}
\phantomsection
\label{lem:cost-dwaki}
\begin{enumerate}
\item
\label{it:cost-dwaki-1}
For any $x\in W_\ex$,
we have an equivalence
\[
\Ras^{\I^\quar_J|\I}(\nabla_x^{J,\quar}) \simeq \sW^{J,\vee}_x. 
\]
\item
\label{it:cost-dwaki-2}
Assume that $p>0$ and $\k$ contains a nontrivial $p$-th root of unity, and let $A \subset S_\af$ be finitary.
For any $x\in W_\ex^{(A,J),\max}$,
we have an equivalence
\[
\Ras^{\I^\quar_J|\I}(\nabla_{\chi_A,x}^{J,\quar}) \simeq \sW^{J,\vee}_{\chi_A,x}. 
\]
\end{enumerate}
\end{lem}

In the setting of Lemma~\ref{lem:cost-dwaki}\eqref{it:cost-dwaki-1}, using Lemma~\ref{lem:Av-Ras} and~\eqref{eqn:Av-dual-Waki}, Lemma~\ref{lem:cost-dwaki} implies that for any $x\in W_\ex^{(A,J),\max}$ we have
\begin{equation}
\label{eqn:Avchi-nabla-quar}
\Av^A_\chi(q^! \nabla_x^{J,\quar}) \simeq \nabla_{\chi_A,x}^{J,\quar}.
\end{equation}

Let $L\subset J\subset S$. 
By Lemma~\ref{lem:cost-dwaki}\eqref{it:cost-dwaki-1}, Lemma~\ref{lem:Av-Ras}\eqref{it:Av-Ras-1}, the commutative diagram obtained from~\eqref{eq:comm-diag-ras-oblv} by passing to right adjoints, and Lemmas~\ref{lem:Av-quar-semi-D-N}\eqref{it:Av-quar-semi-D-N-1} and \ref{lem:push-pull-pK-D-N-semi}\eqref{it:push-pull-pK-D-N-semi-1}, for any finitary $K\subset S_\af$ and $x\in {}^LW^{(K,J)}_\ex$ we have an equivalence 
\begin{equation}
\label{equ:costd-K-colim}
\Ras^{\I^\semi_{L\subset J}\mid \I_L}(\nabla_{K,x}^{L\subset J,\semi}) 
\simeq 
\colim_{\{\mu\in \Xd_{J,x}^{0,+} \mid t_\mu x \in W_\ex^K\}} \big( 
\nabla^L_{t_{-\mu}} \star^{\I} \nabla_{K,t_\mu x} 
\big).
\end{equation}

The following statement generalizes \cite[Proposition~6.12]{adr}. 

\begin{lem}
\label{lem:costandard-perverse}
Let $L\subset J\subset S$ and let $A\subset S_\af$ be finitary. 
The costandard objects $\nabla_x^{L\subset J,\semi}$ for $x \in {}^L W_\ex$ and $\nabla_{\chi_A,x}^{L\subset J,\semi}$ for $x \in {}^L W_\ex^{(A,J),\max}$ (when $p>0$ and $\k$ contains a nontrivial $p$-th root of unity) are contained in the heart of the perverse t-structure on the corresponding $\infty$-categories. 
\end{lem}

\begin{proof}
We show the assertion for the objects $\nabla_x^{L\subset J,\semi}$; the other case follows using Lemma~\ref{lem:pullback-Av-texact}\eqref{it:pullback-q-texact}--\eqref{it:Av-texact}, Lemma~\ref{lem:av-oblv-whit}\eqref{it:av-oblv-whit-1} and~\eqref{eqn:Avchi-nabla-quar}.
It is clear from definitions that our object belongs to $\Shv\left(\I_{L\subset J}^\semi \backslash \Fl \right)^{\geq 0}$, so what we need to show is that it also belongs to $\Shv\left(\I_{L\subset J}^\semi \backslash \Fl \right)^{\leq 0}$. 
Since our t-structure is compatible with colimits, by \eqref{equ:costd-K-colim} it is enough to show that for any $\mu\in \Xd_{J,x}^{0,+}$ we have
\[
\Ras^{\I_L\mid \I^\semi_{L\subset J}}(\nabla^L_{t_{-\mu}} \star^{\I} \nabla_{t_\mu x}) 
 \in
\Shv\left(\I_{L\subset J}^\semi \backslash \Fl \right)^{\leq 0}. 
\]
For this, by Lemma~\ref{lem:oblv-semi-quar-texact} and Proposition~\ref{prop:tstr-convolution-Wak}, it suffices to show that for any $\lambda, \mu \in \Xd_{J,x}^{0,+}$ such that $\mu \unlhd^J \lambda$ we have
\begin{equation}
\label{equ:costd-perv}
\Delta_{t_\lambda} \star^\I \oblv^{\I_L\mid \I}(\nabla^L_{t_{-\mu}}) \star^{\I} \nabla_{t_\mu x}
 \in
\Shv(\I\backslash \Fl)^{\leq 0}. 
\end{equation}
Now we have equivalences 
\begin{multline*}
\Delta_{t_\lambda} \star^\I \oblv^{\I_L\mid \I}(\nabla^L_{t_{-\mu}})
\simeq \Delta_{L,t_\lambda} \star^{\I_L} \nabla^L_{t_{-\mu}} 
\simeq \oblv^{\I_L\mid \I}(\Delta^L_{L,t_\lambda} \star^{\I_L} \nabla^L_{t_{-\mu}}) \\ 
\simeq \oblv^{\I_L\mid \I}(\Delta^L_{L,t_\lambda} \star^{\I_L} \nabla^L_{L,t_{-\mu}} \star^{\I_L} \nabla^L_e) 
\simeq \oblv^{\I_L\mid \I}(\Delta^L_{L,t_{\lambda-\mu}} \star^{\I_L} \nabla^L_e) \\
\simeq \oblv^{\I_L\mid \I}(\Delta^L_{t_{\lambda-\mu}}),
\end{multline*}
where the second equivalence follows from the fact that $\I z^\lambda \I_L/\I_L= \I_L z^\lambda \I_L/\I_L$, the third and the last one follow from the isomorphism induced by multiplication 
\[
\I_L z^\theta \I_L \times^{\I_L} \I_L/\I \xs \I_L z^\theta \I_L/\I
\]
for all $\theta\in \Xd_J^0$,
and the fourth equivalence follows from \cite[Lemma 3.1.5]{DL23} (which also holds in our setting, with the same proof).
Therefore, the left-hand side of~\eqref{equ:costd-perv} identifies with $\oblv^{\I_L\mid \I}(\Delta^L_{t_{\lambda-\mu}}) \star^\I \nabla_{t_\mu x}$.
Now $\oblv^{\I_L\mid \I}(\Delta^L_{t_{\lambda-\mu}})$ is perverse, so that by Corollary~\ref{cor:right-t-exactness-convolution-nabla} the object $\oblv^{\I_L\mid \I}(\Delta^L_{t_{\lambda-\mu}}) \star^\I \nabla_{t_\mu x}$ belongs to $\Shv(\I\backslash \Fl)^{\leq 0}$, which finishes the proof.
\end{proof}

\begin{rmk}
We have not included in Lemma~\ref{lem:costandard-perverse} the case of the parabolic affine flag varieties $\Fl_K$ with $K \neq \varnothing$. In fact, when $J=S$, it is a standard fact that the objects $\nabla^{L}_{K,x}$ are not always perverse if $L$ and $K$ are both nonempty. 
It is also observed in~\cite[\S 6.6]{adr} that the objects $\nabla_{S,x}^{\semi}$ are never perverse if the characteristic of $\k$ is good for $G$. (More specifically, this reference treats the case $x=e$. The general case follows using~\cite[Eqn.~(5.6)]{adr}.) We do not know what happens in other cases.
\end{rmk}

%

\subsection{t-exactness of convolution actions}

\subsubsection{Right convolution with central sheaves}

Consider the abelian category
\[
\Shvfg(\sL^+G \backslash \Gr)^\heartsuit,
\]
i.e.~the heart of the t-structure of middle perversity on $\Shvfg(\sL^+G \backslash \Gr)$. This subcategory is stable under the convolution product $\star^{\sL^+ G}$, and the geometric Satake equivalence~\cite{mv} provides an equivalence of monoidal abelian categories
\begin{equation}
\label{eqn:geom-Satake}
(\Shvfg(\sL^+G \backslash \Gr)^\heartsuit, \star^{\sL^+ G}) \simeq (\Rep^{\heartsuit,\fd}(G^\vee_\k), \otimes) 
\end{equation}
where $G^\vee_\k$ is the split reductive group scheme over $\k$ which is Langlands dual to $G$ and $\Rep^{\heartsuit,\fd}(G^\vee_\k)$ is its abelian category of finite-dimensional representations.
Recall that for any finitary subset $K \subset S_\af$ we have a ``central functor''
\[
\ZCen_K : \Shvfg(\sL^+G \backslash \Gr)^\heartsuit \to \Shvfg(\I_K \backslash \Fl_K)^\heartsuit
\]
defined by nearby cycles along a suitable degeneration of $\Gr$ into $\Fl_K$. The case $K=\varnothing$ is studied in detail in~\cite{ARbook}, and the general case is similar. When $K=S$ we have $\ZCen_S=\id$, and for any $K$ we have $(p_K)_* \circ \ZCen_\varnothing \simeq \oblv^{\I_K | \I} \circ \ZCen_K$.

This functor is ``monoidal'' in the sense that for $\sF,\sG$ in $\Shvfg(\sL^+G \backslash \Gr)^\heartsuit$ we have a canonical identification
\[
\ZCen_K(\sF \star^{\sL^+G} \sG) \simeq \ZCen_K(\sF) \star^{\I_K} \ZCen_K(\sG).
\]
(See~\cite[Theorem~3.4.1]{ARbook} for the case $K=\varnothing$; the general case is similar. Here we write monoidal in quotation marks because $\Shvfg(\I_K \backslash \Fl_K)^\heartsuit$ has no natural monoidal structure.)
Recall also that the central sheaves are ``convolution exact'' in the sense that for any $\sF$ in $\Shvfg(\sL^+G \backslash \Gr)^\heartsuit$ the functor
\[
(-) \star^{\I_K} \ZCen_K(\sF) : \Shv(\Fl_K) \to \Shv(\Fl_K)
\]
of right convolution with $\ZCen_K(\sF)$
is t-exact with respect to the perverse t-structure. (See~\cite[Corollary~3.2.5]{ARbook} for the case $K=\varnothing$; the general case is similar.) 
We will now establish a similar property for semiinfinite sheaves, which generalizes~\cite[Proposition~3.3.15]{DL23}. (That reference treats the case when $L=J$, $K=S$ and $\k$ has characteristic $0$. Note that, even specializing to this case, our proof is different.)

\begin{prop}
\label{prop:right-conv-t-exact}
Let $L \subset J \subset S$, and let $K \subset S_\af$ be finitary. For any $\sF$ in $\Shvfg(\sL^+G \backslash \Gr)^\heartsuit$ the functor
\[
(-) \star^{\I_K} \ZCen_K(\sF) : \Shv \left( \I_{L \subset J}^\semi \backslash \Fl_K \right) \to \Shv \left( \I_{L \subset J}^\semi \backslash \Fl_K \right)
\]
is t-exact with respect to the perverse t-structure.
\end{prop}

\begin{proof}
Using Lemma~\ref{lem:oblv-semi-quar-texact} one reduces the proof to the case $L=\varnothing$.
In this case, first we show right t-exactness. For this we need to show that for any $x \in W_\ex^{(K,J)}$ we have $\Delta^{J,\quar}_{K,x} \star^{\I_K} \ZCen_K(\sF) \in \Shv \left( \I_{J}^\quar \backslash \Fl_K \right)^{\leq 0}$. Since right convolution commutes with Raskin's equivalence, this object is finitely generated, and the desired property is clear from the characterization of the t-structure in Proposition~\ref{prop:tstr-convolution-Wak} and convolution exactness of central sheaves as recalled above.

Left t-exactness follows, since the functor under consideration has left adjoint given by $(-) \star^{\I_K} \ZCen_K(\sF^\vee)$ where $\sF^\vee$ is the dual to $\sF$ in the rigid monoidal category $(\Shvfg(\sL^+G \backslash \Gr)^\heartsuit, \star^{\sL^+ G})$, which is right t-exact by the case treated above.
\end{proof}

\subsubsection{Left convolution with central sheaves}

Now we consider convolution on the left. Recall from~\S\ref{ss:convolution} that we have an action of $\Shv \left( \I_{L \subset J} \backslash \sL M_J / \I_{L \subset J} \right)$ on $\Shv \left( \I_{L \subset J}^\semi \backslash \Fl_K \right)$ by convolution on the left. We have the ``central functor'' associated with the reductive group $M_J$ and the parahoric subgroup defined by $L$, which provides a ``monoidal'' functor
\begin{equation}
\label{eqn:ZCen-MJ}
\ZCen_{L \subset J} : \Shvfg(\sL^+ M_J \backslash \sL M_J / \sL^+ M_J)^\heartsuit \to \Shvfg \left( \I_{L \subset J} \backslash \sL M_J / \I_{L \subset J} \right)^\heartsuit.
\end{equation}
Recall also that the connected components of the affine Grassmannian $\sL M_J / \sL^+ M_J$ are naturally parametrized by $\Xd / \Z \check{\Phi}_J$ where $\check{\Phi}_J$ is the coroot system of $(M_J,T)$, and that the pairing $\langle -,2\rho-2\rho_J \rangle : \Xd \to \Z$ factors through $\Xd / \Z \check{\Phi}_J$.

The following statement generalizes~\cite[Proposition~3.3.20(ii)]{DL23}.

\begin{prop}
\label{prop:exactness-left-conv}
If $\sF \in \Shvfg(\sL^+ M_J \backslash \sL M_J / \sL^+ M_J)^\heartsuit$ is supported on the component corresponding to $\zeta \in \Xd / \Z \check{\Phi}_J$, then the functor
\[
\ZCen_{L \subset J}(\sF) \star^{\I_{L \subset J}} (-)[-\langle \zeta, 2\rho-2\rho_J \rangle] : \Shv \left( \I_{L \subset J}^\semi \backslash \Fl_K \right) \to \Shv \left( \I_{L \subset J}^\semi \backslash \Fl_K \right)
\]
is t-exact.
\end{prop}

\begin{proof}
The proof is an easy adaptation of that of~\cite[Proposition~3.3.20(ii)]{DL23}, and boils down in the end to convolution exactness of central sheaves. Namely, using Lemma~\ref{lem:pullback-Av-texact}\eqref{it:pullback-p-texact} it suffices to treat the case $K=\varnothing$, and using Lemma~\ref{lem:oblv-semi-quar-texact} one can also assume that $L=\varnothing$. By the same argument as in the proof of Proposition~\ref{prop:right-conv-t-exact} it suffices to prove right t-exactness of our functor, which amounts to proving that for all $w \in W_\ex$ we have
\begin{equation}
\label{eqn:proof-exactness-left-conv}
\ZCen_{\varnothing \subset J}(\sF) \star^{\I_{\varnothing \subset J}} \Delta^{J, \quar}_w [-\langle \zeta, 2\rho-2\rho_J \rangle] \in \Shv \left( \I_{J}^\quar \backslash \Fl \right)^{\leq 0}.
\end{equation}

Fix $\eta \in \Xd_J^0 \cap \Xd_{J^\comp}^{++}$, and consider the action of $\Gm$ on $\Fl$ defined by $\eta$. Then the fixed points $\Fl^0$ for this action identify with a disjoint union of copies of the affine flag variety $\sL M_J / \I_{\varnothing \subset J}$ parametrized by the minimal representatives for the cosets in $W_J \backslash W$. (Here the copy parametrized by $x$ is the $\sL M_J$-orbit of $\dot{x} \I/\I$, and will be denoted $\Fl^0_x$.) Correspondingly, the attractor $\Fl^+$ for this action (in the sense discussed e.g.~in~\cite[\S 1.1.2]{ARbook}) admits a decomposition $\Fl^+ = \sqcup_x \Fl^+_x$ where each component is $\I_{J}^\quar$-stable. For any $x$ we have a ``limit at $0$'' morphism $a_x : \Fl^+_x \to \Fl^0_x$ and a canonical section $b_x : \Fl^0_x \to \Fl^+_x$, and the natural morphism $c_x : \Fl^+_x \to \Fl$ is a locally closed immersion. As in~\cite[Lemma~3.3.2]{DL23}, the $!$-pullback functor associated with $a_x$ induces an equivalence
\[
\Shv(\I_{\varnothing \subset J} \backslash \Fl_x^0)  \simto  \Shv \left( \I_{J}^\quar \backslash \Fl^+_x \right)
\]
with inverse given by $(a_x)_!$ and also by $(b_x)^!$.

Let now $w$ be as above, and write $w=yx$ with $y$ in the extended affine Weyl group $W_{J,\ex}$ of $(M_J,T)$ and $x \in W$ minimal in $W_J x$. Write also $y=t_\lambda z$ with $\lambda \in \Xd$ and $z \in W_J$. 
Then using Remark~\ref{rmk:lJ}\eqref{it:lJ-WJext} we see that $\Delta^{J, \quar}_w$ is the image under the functor
\[
(c_x)_! (a_x)^![- \langle \lambda, 2\rho-2\rho_J \rangle - \ell(x)]
\]
of the standard perverse sheaf in $\Shv(\I_{\varnothing \subset J} \backslash \Fl_x^0)$ associated with $y$ (where $\Fl^0_x$ is identified with the affine flag variety of $M_J$ as above). Hence the object in~\eqref{eqn:proof-exactness-left-conv} is the image under $(c_x)_! (a_x)^![-\langle \lambda, 2\rho-2\rho_J \rangle - \ell(x) -\langle \zeta, 2\rho-2\rho_J \rangle]$ of the convolution of $\ZCen_{\varnothing \subset J}(\sF)$ with that standard perverse sheaf. By convolution exactness of central sheaves (now for $M_J$) the latter convolution is perverse. Taking also into account its support, it is therefore a finite extension of nonnegative cohomological shifts of standard perverse sheaves associated with elements of the form $t_\mu u$ with $u \in W_J$ and $\mu \in \lambda+ \zeta \subset \Xd$. Here we have $\langle \lambda, 2\rho-2\rho_J \rangle + \langle \zeta, 2\rho-2\rho_J \rangle = \langle \mu, 2\rho-2\rho_J \rangle$, so these considerations show that the object in~\eqref{eqn:proof-exactness-left-conv} is a finite extension of nonnegative cohomological shifts of standard objects, hence belongs to $\Shv \left( \I_{J}^\quar \backslash \Fl \right)^{\leq 0}$ as desired.
\end{proof}

\section{A t-exactness theorem} 
\label{sec:t-exactness}

\subsection{Statement} 

Let $X$ be either $\Fl_K$ for some finitary $K \subset S_\af$, $\tFl$, or the symbol $\sL G/(\Iu^A,\chi_A)$ for some finitary $A \subset S_\af$ (assuming that $p>0$ and that $\k$ contains a nontrivial $p$-th root of unity). Our goal in this section is to prove the following statement.

\begin{thm}
\label{thm:texactness}
The composition of functors 
\[
\mathsf{F}:
\Shv\left(\I_J^{\semi}\backslash X\right) \xrightarrow{\oblv^{\I_J^{\semi}|\I_J^{\quar}}} \Shv\left(\I_J^{\quar}\backslash X\right) 
\xrightarrow{\mathrm{Ras}^{\I_J^{\quar}\mid \I^\semi}} \Shv\left(\I^{\semi}\backslash X\right)
\] 
is t-exact with respect to the perverse t-structure on both sides. 
\end{thm}

Standard arguments allow to reduce the proof to the case where for any $\alpha\in \Phis$ there exists a corresponding fundamental coweight 
$\check{\varpi}_{\alpha}\in \Xd$, i.e.~an element which satisfies $\langle \check{\varpi}_{\alpha}, \beta \rangle =\delta_{\alpha,\beta}$ for any $\beta\in \Phis$. We will therefore assume that this condition is satisfied, and fix choices for these fundamental coweights. 
By other standard arguments (see e.g.~\cite[\S 6]{bbdg}), one can also assume that $p>0$ and $\k$ contains a nontrivial $p$-th root of unity, which we will do.

\subsection{Reduction to the Iwahori--Whittaker case} 
\label{subsect:redtn-IW}

We write $\mathsf{F}$ as $\mathsf{F}_{D}$ if $X=\sL{G}/D$, where $D=\I_K$, $\Iu$ or the symbol $(\Iu^A,\chi_A)$. 
If $K=S$, then $\Iu^-:=\Iu^S=\Ad_{\dot{w}_0}(\Iu)$ is the pro-unipotent radical of the Iwahori subgroup defined by the Borel subgroup of $G$ which is opposite to $B$ with respect to $T$, and $(\Iu^S,\chi_S)$ is known as the Iwahori--Whittaker condition. 
We will abbreviate $\IW=(\Iu^S,\chi_S)$ and $\Av_{\IW}:=\Av^S_{\chi} \circ q^!$. 

We consider the affine braid group $\mathfrak{B}_\af$ associated with $W_\af$, i.e.~the group generated by symbols $\widetilde{w}$ for $w \in W_\af$, with relations $\widetilde{w} \cdot \widetilde{y} = \widetilde{wy}$ if $\ell(wy)=\ell(w)+\ell(y)$.
We also let $\mathfrak{B}^+_\af \subset \mathfrak{B}_\af$ be the monoid of positive elements, i.e.~the submonoid generated by the elements $\widetilde{w}$ for $w \in W_\af$.
We will denote by $\ell : \mathfrak{B}^+_\af \to \Z_{\geq 0}$ the unique monoid map such that $\ell(\widetilde{w})=\ell(w)$ for any $w \in W_\af$, and will call $\ell(b)$ the length of $b$.

It is a standard fact (which follows e.g.~from~\cite[Lemma~4.1]{adr}) that the assignment $\widetilde{w} \mapsto \nabla_w$ extends to a group morphism from $\mathfrak{B}_\af$ to the group of isomorphism classes of invertible objects in $\Shv(\I\backslash \sL{G}/\I)$. We will denote by $\nabla_b$ the image of $b \in \mathfrak{B}_\af$. In particular, for $w \in W_\af$ we have $\nabla_{\widetilde{w}^{-1}} = \Delta_{w^{-1}}$. Similarly, the assignment $\widetilde{w} \mapsto \Delta_w$ extends to a group morphism from $\mathfrak{B}_\af$ to the group of isomorphism classes of invertible objects in $\Shv(\I\backslash \sL{G}/\I)$, and we will denote by $\Delta_b$ the image of $b \in \mathfrak{B}_\af$.

The following statement is an extension to general coefficients of~\cite[Proposition~3]{BLin}, whose proof (which is similar to that in~\cite{BLin}) will be explained in~\S\ref{ss:pf-lem}.

\begin{lem}
\label{lem:bez-lin}
Let $\sG\in \Shvfg(\I\backslash \Fl)$. 
We have $\sG\in \Shv(\I\backslash \Fl)^{\leq 0}$ if and only if
\[
\Av_\IW(\sG \star^\I \nabla_b)\in \Shv(\I\backslash \sL{G}/\IW)^{\leq 0}
\]
for any $b\in \mathfrak{B}^+_\af$. 
\end{lem}

We deduce the following corollary.

\begin{corollary}
\label{cor:tstr-criterion-IW}
Let $\sG\in \Shv_{\omega}\left(\I^\semi\backslash \Fl\right)$. 
We have $\sG\in \Shv\left(\I^\semi\backslash \Fl\right)^{\leq 0}$ if and only if $\Av_\IW(\sG \star^\I \nabla_b)\in \Shv\left(\I^\semi\backslash \sL{G}/\IW\right)^{\leq 0}$ for any $b\in \mathfrak{B}^+_\af$. 
\end{corollary}

\begin{proof}
Since the equivalence $\Ras^{\I^\semi \mid \I}$ commutes with $\Av_\IW$ and right convolution with any $\nabla_b$, the assertion follows from Proposition~\ref{prop:tstr-convolution-Wak} and Lemma~\ref{lem:bez-lin}. 
\end{proof}

The following lemma will be used to reduce the proof of Theorem~\ref{thm:texactness} to the case $X=\sL{G}/\IW$. 

\begin{lem}
\phantomsection
\label{lem:reduction-texactness}
\begin{enumerate}
\item 
\label{it:reduction-texactness-1}
If $\mathsf{F}_\I$ is right t-exact, then so are $\mathsf{F}_{\Iu}$, $\mathsf{F}_{\I_K}$ and $\mathsf{F}_{(\Iu^A,\chi_A)}$ for any finitary $K,A \subset S_\af$. 
\item 
\label{it:reduction-texactness-2}
If $\mathsf{F}_{\IW}$ is right t-exact, then so is $\mathsf{F}_\I$. 
\end{enumerate}
\end{lem}

\begin{proof}
\eqref{it:reduction-texactness-1}
To show that $\mathsf{F}_{\Iu}$, resp. $\mathsf{F}_{(\Iu^A,\chi_A)}$, is right t-exact, it suffices to show that $\mathsf{F}_{\Iu}(q^!(\Delta_x^{J,\semi}))$, resp. $\mathsf{F}_{(\Iu^A,\chi_A)}(\Delta^{J,\semi}_{\chi_A,x})$, is contained in the ``$\leq 0$'' part of the perverse t-structure for any $x\in {}^JW_\ex$, resp. $x\in {}^JW_\ex^{(A,J),\max}$. 
This follows from Lemma~\ref{lem:Av-Ras}, Lemma~\ref{lem:Av-chi-D}\eqref{it:Av-chi-D-1}, Lemma~\ref{lem:pullback-Av-texact}\eqref{it:pullback-q-texact}--\eqref{it:Av-texact}, and the assumption that $\mathsf{F}_{\I}(\Delta_x^{J,\semi})$ is contained in $\Shv(\I^\semi\backslash \Fl)^{\leq 0}$. 

The assertion for $\mathsf{F}_{\I_K}$ follows from Lemma~\ref{lem:Av-Ras}\eqref{it:Av-Ras-1}, Lemma~\ref{lem:pullback-Av-texact}\eqref{it:pullback-p-texact} and the fact that the functor $(p_K)^\dag$ is conservative. 

\eqref{it:reduction-texactness-2}
Let us assume that $\mathsf{F}_{\IW}$ is known to be right t-exact. To show that $\mathsf{F}_{\I}$ is right t-exact,
it suffices to show that for any $\sF\in \Shv_{\omega}\left(\I_J^\semi\backslash \Fl\right)^{\leq 0}$ (in particular, when $\sF=\Delta_x^{J,\semi}$ for some $x \in {}^J W_\ex$) we have $\mathsf{F}_{\I}(\sF)\in \Shv \left(\I^\semi\backslash \Fl\right)^{\leq 0}$. 
By Corollary~\ref{cor:tstr-criterion-IW}, it is equivalent to show that 
\[
\mathsf{F}_{\IW}(\Av_\IW(\sF\star^\I \nabla_b)) \simeq \Av_\IW(\mathsf{F}_{\I}(\sF) \star^\I \nabla_b) \in \Shv\left(\I^\semi\backslash \sL{G}/\IW\right)^{\leq 0} 
\]
for all $b\in \mathfrak{B}^+_\af$.
Since by our assumption $\mathsf{F}_{\IW}$ is right t-exact, and since by Lem\-ma~\ref{lem:pullback-Av-texact}\eqref{it:pullback-q-texact}--\eqref{it:Av-texact} the functor $\Av_\IW$ is t-exact, it is enough to show that the functor 
\[
(-) \star^\I \nabla_b : \Shv\left(\I^{\semi}_J\backslash \Fl\right) \to \Shv\left(\I^{\semi}_J\backslash \Fl\right)
\]
is right t-exact for any $b \in \mathfrak{B}_\af^+$, which follows from Corollary~\ref{cor:right-t-exactness-convolution-nabla}.
\end{proof}


\subsection{Proof of Theorem~\ref{thm:texactness}} 

For $\lambda\in \Xd_J^{++}$, we abbreviate 
\[
\Delta^{J,\semi}_{\IW,\lambda}:=\Delta^{J,\semi}_{\chi_S,w_Jt_{\lambda}w_0}.
\]
By Lemma~\ref{lem:JSJw},
these objects are exactly the
standard perverse sheaves in the $\infty$-category $\Shv \left( \I^\semi_J\backslash \sL{G}/\IW \right)$. By Remark~\ref{rmk:JSJw-length}, the cohomological shift appearing in the definition of this object is $-\langle \lambda, 2\rho \rangle - \ell(w_0)$.

For $\lambda\in \Xd_J^+$, we will denote by 
\[
\sA^{M_J}_{\lambda} \in \Shv(\sL^+M_J\backslash \sL{M_J}/\sL^+M_J)^\heartsuit 
\]
the 0-th perverse cohomology of the standard complex attached to the $(\sL^+M_J\times \sL^+M_J)$-orbit of $z^\lambda$. (These objects are the perverse sheaves corresponding to the Weyl modules for the reductive group Langlands dual to $M_J$ under the geometric Satake equivalence~\eqref{eqn:geom-Satake} for $M_J$; see e.g.~\cite[\S 1.5.1]{ARbook}.)

We set $\zeta_J=\sum_{\alpha\in \Phis_J} \check{\varpi}_{\alpha}$; then $\Xd_J^{++}=\zeta_J+\Xd_J^{+}$. 

\begin{lem}
\label{lem:D-convolution}
For any $\lambda\in \Xd_J^{+}$, we have
\[
\Delta^{J,\semi}_{\IW,\lambda+\zeta_J} \simeq \sA^{M_J}_{\lambda} \star^{\sL^+M_J} \Delta^{J,\semi}_{\IW,\zeta_J} [-\langle \lambda, 2\rho-2\rho_J \rangle]. 
\]
\end{lem}

\begin{proof}
We consider 
the pro-unipotent radical $\I^-_{\varnothing\subset J,\mathrm{u}}$ of the Iwahori subgroup in $\sL{M_J}$ defined by the Borel subgroup of $M_J$ which is opposite to $B_J$ with respect to $T$, and fix a morphism $\I^-_{\varnothing\subset J,\mathrm{u}} \to \Ga$ as in~\S\ref{sss:Whittaker}. Considering the associated Iwahori--Whittaker condition on sheaves on the ind-scheme $\sL^+ M_J \backslash \sL M_J$ we obtain the $\infty$-category $\Shv(\sL^+M_J\backslash \sL{M_J}/\IW_J)$.

On the other hand, considering the pullback under the induced morphism $\I^-_{\varnothing\subset J,\mathrm{u}}\ltimes \sL{U_J} \to \Ga$ of our fixed Artin--Schreier local system on $\Ga$ we obtain a twisted equivariance condition on sheaves on $\sL G$, which we will denote
$\I^\quar_{J,\sW}$. 
We can form the ``bi-Whittaker" category 
\[
\Shv \left( \I^\quar_{J,\sW}\backslash \sL{G}/\IW \right)
\]
following the same recipe as in~\S\ref{sss:Whittaker}; we will
denote by $\Delta^{(J,\sW),\quar}_{\IW,e}$ the standard sheaf
associated with the $(\I^-_{\varnothing\subset J,\mathrm{u}}\ltimes \sL{U_J})\times \Iu^-$-orbit through $e$, i.e.~the $!$-pushforward under the natural embedding of the dualizing sheaf on this orbit, shifted by $-\ell(w_0)$.

Convolution gives a functor
\[
\star^{\IW_J}: \Shv \left( \sL^+M_J\backslash \sL{M_J}/\IW_J \right) \times \Shv \left( \I^\quar_{J,\sW}\backslash \sL{G}/\IW \right) \rightarrow \Shv \left( \I^\semi_{J}\backslash \sL{G}/\IW \right). 
\] 
For $\mu\in \Xd_J^{++}$, let $\sA^{M_J}_{\sW,\mu}\in \Shv(\sL^+M_J\backslash \sL{M_J}/\IW_J)$ be the standard sheaf associated with the orbit $\sL^+M_J\cdot z^\mu\cdot \I^-_{\varnothing\subset J,\mathrm{u}}$, i.e.~the $!$-pushforward under the natural embedding of the dualizing sheaf on this orbit, shifted by $-\langle \mu, 2\rho_J \rangle$. 
For any $\nu\in \Xd$, multiplication induces an isomorphism 
\[
\sL^+M_J\cdot z^{\nu}\cdot \I^-_{\varnothing\subset J,\mathrm{u}} \times^{\I^-_{\varnothing\subset J,\mathrm{u}}} (\I^-_{\varnothing\subset J,\mathrm{u}}\ltimes \sL{U_J})\cdot \Iu^- \xs 
\I^\semi_J\cdot z^{\nu}\cdot \Iu^-.
\]
Hence we have an equivalence, for any $\mu\in \Xd_J^{++}$, 
\[
\sA^{M_J}_{\sW,\mu}\star^{\IW_J} \Delta^{(J,\sW),\quar}_{\IW,e} [-\langle \mu, 2\rho-2\rho_J \rangle]
\simeq \Delta^{J,\semi}_{\IW,\mu}. 
\]
Now the assertion follows from the fact that
\[
\sA_{\lambda}^{M_J}\star^{\sL^+M_J} \sA_{\sW,\zeta_J}^{M_J}\simeq \sA_{\sW,\lambda+\zeta_J}^{M_J},
\]
for any $\lambda \in \Xd_J^{+}$,
see~\cite[Proof of Theorem 3.9]{BGMRR}.
\end{proof}

\begin{proof}[Proof of Theorem~\ref{thm:texactness}]
By Lemma~\ref{lem:oblv-semi-quar-texact} the functor $\oblv^{\I_J^{\semi}|\I_J^{\quar}}$ is t-exact, and by Corollary~\ref{cor:Ras-right-texact} the functor $\mathrm{Ras}^{\I_J^{\quar}\mid \I^\semi}$ is left t-exact, as the inverse of a right t-exact equivalence. Hence $\mathsf{F}$ is left t-exact. 
By Lemma~\ref{lem:reduction-texactness}, to conclude we therefore only need to show the right t-exactness of $\mathsf{F}_\IW$. 
For this it suffices to show that for any $\lambda\in \Xd_J^+$ we have
\[
\mathsf{F}_\IW \left( \Delta^{J,\semi}_{\IW,\lambda+\zeta_J} \right) \in \Shv \left( \I^{\semi}\backslash \sL{G}/\IW \right)^{\leq 0},
\]
i.e.~(in view of Proposition~\ref{prop:tstr-convolution-Wak}) that for any $\mu \in \Xd$ we have
\[
\sW_{t_\mu} \star^\I \Ras^{\I^\semi \mid \I} \left( \mathsf{F}_\IW \left( \Delta^{J,\semi}_{\IW,\lambda+\zeta_J} \right) \right)
\in \Shv \left( \I \backslash \sL{G}/\IW \right)^{\leq 0}.
\]

By Lemma~\ref{lem:D-convolution}, we have 
\begin{multline*}
\mathsf{F}_\IW \left( \Delta^{J,\semi}_{\IW,\lambda+\zeta_J} \right) \\
\simeq \mathrm{Ras}^{\I_J^{\quar}\mid \I^\semi}\circ \oblv^{\I_J^\semi | \I_{J}^\quar } \left( \sA^{M_J}_{\lambda} \star^{\sL^+M_J} \Delta^{J,\semi}_{\IW,\zeta_J} \right) [-\langle \lambda, 2\rho-2\rho_J \rangle] \\ 
\simeq \mathrm{Ras}^{\I_J^{\quar}\mid \I^\semi} \left( \ZCen_{\varnothing \subset J}(\sA^{M_J}_{\lambda}) \star^{\I_{\varnothing \subset J}} \oblv^{\I_J^\semi | \I_J^\quar} (\Delta^{J,\semi}_{\IW,\zeta_J}) \right) [-\langle \lambda, 2\rho-2\rho_J \rangle]
\end{multline*}
where $\ZCen_{\varnothing \subset J}$ is as in~\eqref{eqn:ZCen-MJ}.
Now, recall the Wakimoto sheaves for $M_J$ considered in Lemma~\ref{lem:Ras-convolution-Wak}.
Since $\ZCen_{\varnothing \subset J}(\sA^{M_J}_{\lambda})$ admits a filtration by objects $\sW^{M_J}_{t_\mu}$ where $\mu + \Z\check{\Phi}_J = \lambda + \Z\check{\Phi}_J$ (so that $\langle \lambda, 2\rho-2\rho_J \rangle = \langle \mu, 2\rho-2\rho_J \rangle$), see e.g.~\cite[Theorem~4.4.5]{ARbook}, by Lemma~\ref{lem:Ras-convolution-Wak} it suffices to show that 
\begin{equation}
\label{eq:FIW-D}
\mathsf{F}_\IW \left( \Delta^{J,\semi}_{\IW,\zeta_J} \right) \in \Shv \left( \I^{\semi}\backslash \sL{G}/\IW \right)^{\leq 0}. 
\end{equation}

We set $w_{\zeta_J}:=w_Jt_{\zeta_J}w_0$. By 
Lemma~\ref{lem:JSJw} and its proof,
this element has the following properties.

\begin{claim}
\label{cl:properties-w-zetaJ}
\begin{enumerate}
\item 
\label{it:w-zetaJ-1}
$w_{\zeta_J}\in {}^JW_\ex^{S,\max} \cap {}^JW_\ex^{(S,J),\max}$; 
\item 
\label{it:w-zetaJ-3}
$w_{\zeta_J}$ is a minimal element in ${}^JW_\ex^{S,\max}$ for the Bruhat order.
\end{enumerate}
\end{claim}


By Proposition~\ref{prop:Ras-D}\eqref{it:Ras-D-2} we have 
\[
\mathrm{Ras}^{\I_J^{\quar}\mid \I} \left( \Delta^{J,\quar}_{\chi_S,w_{\zeta_J}} \right)
\simeq \sW^{J}_{\chi_S,w_{\zeta_J}}.
\]
Hence by Claim~\ref{cl:properties-w-zetaJ}\eqref{it:w-zetaJ-1}, Lemma~\ref{lem:av-oblv-whit}\eqref{it:av-oblv-whit-1}, and the diagram obtained from~\eqref{eq:comm-diag-ras-oblv} by passing to left adjoints,
we have 
\begin{multline*}
\mathrm{Ras}^{\I_J^{\semi}\mid \I_J} \left( \Delta^{J,\semi}_{\IW,\zeta_J} \right) \simeq
\mathrm{Ras}^{\I_J^{\semi}\mid \I_J} \circ \Av^{\I^\semi_J \mid \I^\quar_J}_! \left (\Delta^{J,\quar}_{\chi_S,w_{\zeta_J}} \right) [-\ell(w_J)] \\
\simeq \Av^{\I_J \mid \I}_! \circ \mathrm{Ras}^{\I_J^{\quar}\mid \I} \left (\Delta^{J,\quar}_{\chi_S,w_{\zeta_J}} \right) [-\ell(w_J)]
\simeq \Av^{\I_J \mid \I}_! \left( \sW^{J}_{\chi_S,w_{\zeta_J}} \right) [-\ell(w_J)].
\end{multline*}
This complex is supported on the closure of the orbit parametrized by $w_{\zeta_J}$, and by Lemma~\ref{lem:av-oblv-whit}\eqref{it:av-oblv-whit-1} and the properties of Wakimoto objects explained in~\S\ref{sss:Wakimoto} its restriction to this orbit coincides with that of $\Delta^J_{\chi_S,w_{\zeta_J}}$.
Using Claim~\ref{cl:properties-w-zetaJ}\eqref{it:w-zetaJ-3}, we deduce equivalences
\begin{multline*}
\mathrm{Ras}^{\I_J^{\quar}\mid \I} \circ \oblv^{\I_J^{\semi} \mid \I_J^\quar} \left( \Delta^{J,\semi}_{\IW,\zeta_J} \right) \simeq
\oblv^{\I_J\mid \I} \circ \mathrm{Ras}^{\I_J^{\semi}\mid \I_J} \left( \Delta^{J,\semi}_{\IW,\zeta_J} \right) \\
\simeq \oblv^{\I_J\mid \I} (\Delta^J_{\chi_S,w_{\zeta_J}}) \simeq \oblv^{\I_J\mid \I} (\nabla^J_{\chi_S,w_{\zeta_J}}).
\end{multline*}
Hence, by Proposition~\ref{prop:tstr-convolution-Wak} (applied with ``$J$'' being the empty set), to conclude it suffices to prove that for any $\theta\in \Xd^+$ the complex
\[
\Delta_{t_\theta}\star^{\I}\oblv^{\I_J\mid \I}(\nabla^J_{\chi_S,w_{\zeta_J}})
\simeq ((p_J)_* \Delta_{t_\theta}) \star^{\I_J} \nabla^J_{\chi_S,w_{\zeta_J}}
\simeq \Delta_{J,t_\theta} \star^{\I_J} \nabla^J_{\chi_S,w_{\zeta_J}}
\]
belongs to $\Shv(\I\backslash \sL{G}/\IW)^{\leq 0}$. (Here, the second equivalence follows from the fact that $t_\theta \in W_\ex^J$.)
Now, interpreting the complexes as living on the ind-scheme $\I \backslash \sL G$, one can write the functor
\[
(-) \star^{\I_J} \nabla^J_{\chi_S,w_{\zeta_J}} : \Shv(\I \backslash \sL G / \I_J) \to \Shv(\I \backslash \sL G / \IW)
\]
as exterior product with a perversely shifted local system, followed by $*$-pushforward under an affine morphism, hence this functor is right t-exact by the same considerations as in the proof of Corollary~\ref{cor:right-t-exactness-convolution-nabla}. The desired claim follows, since $\Delta_{J,t_\theta}$ belongs to $\Shv(\I \backslash \sL G / \I_J)^{\leq 0}$. (In fact, this object is well known to be perverse.)
%
%
\end{proof}

\subsection{Proof of Lemma~\ref{lem:bez-lin}} 
\label{ss:pf-lem}

In this subsection, for the reader's convenience we explain the proof of Lemma~\ref{lem:bez-lin} following~\cite{BLin}.

\begin{proof}[Proof of Lemma~\ref{lem:bez-lin}]
The ``only if" part is easy, since the functor $(-) \star^{\I}\nabla_b$ is right t-exact (see Corollary~\ref{cor:right-t-exactness-convolution-nabla}) and $\Av_{\IW}$ is t-exact (see~\cite[Lemma~3.2]{AR23} or Lemma~\ref{lem:pullback-Av-texact}). 
It remains to show the ``if" part. 
Let $n$ be the maximal integer such that $\sG':={}^p\sH^n(\sG)\neq 0$ (which exists since $\sG\in \Shvfg(\I\backslash \Fl)$). 

\textit{Step 1.} 
We claim that  
\begin{equation}
\label{eq:AvIW-convolution-nonvanishing}
\Av_\IW(\sG'\star^{\I}\nabla_b) \neq 0 \quad \text{for some $b\in \mathfrak{B}_\af^+$.}
\end{equation}
The proof will use the following fact (see e.g.~\cite[Proposition~4.4.4]{ARbook}): there exists a finite subset $X \subset W_\ex$ such that for any $x\in W_\ex$ and $y\in W_\af$, 
\[
(\mathbf{s}^{\varnothing \subset S}_{\varnothing,x})^!(\sG'\star^{\I}\nabla_{y}) \neq 0 \quad \Rightarrow \quad x\in X \cdot y.
\]

Now choose $\lambda\in \Z\check{\Phi}$ 
such that $X \cdot t_\lambda$ consists of elements of the form $wt_\mu$, with $w\in W$ and $\mu\in \Xd^{++}$.
If $x\in X \cdot t_\lambda$ satisfies $(\mathbf{s}^{\varnothing \subset S}_{\varnothing,x})^!(\sG'\star^\I \nabla_{t_\lambda})\neq 0$ and $x' \in W_\ex^{S,\max}$ is the unique element such that $x'W = xW$, then we have 
\begin{align*}
\Hom^\bullet(\Delta_{\chi_S,x'}, \Av_\IW(\sG'\star^{\I}\nabla_{{t_\lambda}}))
&\simeq \Hom^\bullet(\Av^S_!(\Delta_{\chi_S,x'}), q^!(\sG'\star^{\I}\nabla_{{t_\lambda}}))\\
&\simeq \Hom^\bullet(q^! \Delta_x, q^!(\sG'\star^{\I}\nabla_{{t_\lambda}}))\neq 0,
\end{align*}
where in the second equivalence we use the fact that $\Av^S_!(\Delta_{\chi_S,x})$ is an iterated extension of objects $q^! \Delta_{xw}$ for $w\in W$, each appearing once (see Lemma~\ref{lem:Av-chi-D}\eqref{it:Av-chi-D-2}), and that $x'W \cap (X \cdot t_\lambda)=\{x\}$. 
We deduce that $\Av_\IW(\sG'\star^{\I}\nabla_{{t_\lambda}}) \neq 0$, and hence \eqref{eq:AvIW-convolution-nonvanishing}. 

\textit{Step 2.} 
We let $b_0$ be an element of minimal length in $\{b\in \mathfrak{B}^+_\af \mid \Av_\IW(\sG'\star^{\I}\nabla_{b})\neq 0\}$. Then we claim that ${}^p\sH^0(\Av_\IW(\sG'\star^{\I}\nabla_{b_0}))\neq 0$. 

In fact, for any simple reflection $s\in S_\af$, we have a
long exact sequence
\begin{equation*}
0\rightarrow \Delta_e \rightarrow \Delta_s \rightarrow \nabla_s \rightarrow \Delta_e \rightarrow 0
\end{equation*} 
in $\Shv(\I\backslash \Fl)^\heartsuit$.
Using these exact sequences one sees
that for any $b\in \mathfrak{B}_\af^+$ 
we have
\[
\sG'\star^{\I}\nabla_{b} \simeq \sG'\star^{\I}\Delta_{b} \quad \text{mod}\quad \langle \sG'\star\nabla_{b'} : b' \in \mathfrak{B}_\af^+ \text{ with } \ell(b')<\ell(b)\rangle
\]
where in the right-hand side we consider the full stable subcategory generated by the given objects, and the notation means an equivalence in the Verdier quotient.
By our assumption on $b_0$, we therefore have $\Av_{\IW}(\sG'\star^{\I}\nabla_{b_0}) \simeq \Av_{\IW}(\sG'\star^{\I}\Delta_{b_0})$. Since $(-) \star^{\I}\nabla_{b_0}$ is right t-exact, see Corollary~\ref{cor:right-t-exactness-convolution-nabla}, $(-)\star^{\I}\Delta_{b_0}$ is left t-exact by similar considerations, and $\Av_{\IW}$ is t-exact, this object belongs to both the ``$\leq 0$'' and the ``$\geq 0$'' parts of the t-structure, hence is perverse. 
We deduce that ${}^p\sH^0(\Av_\IW(\sG'\star^{\I}\nabla_{b_0}))\simeq \Av_\IW(\sG'\star^{\I}\nabla_{b_0})\neq 0$, as desired.

\textit{Step 3.} 
We now conclude. Let $b_0$ be as in Step 2.
By right t-exactness of the $\mathfrak{B}_\af^+$-action and t-exactness of $\Av_\IW$, we have ${}^p\sH^{>n}(\Av_\IW(\sG\star^{\I}\nabla_{b_0}))=0$ and 
\[
{}^p\sH^n(\Av_\IW(\sG\star^{\I}\nabla_{b_0}))\simeq {}^p\sH^0(\Av_\IW(\sG'\star^{\I}\nabla_{b_0})) \neq 0. 
\]
Hence the fact that $\Av_\IW(\sG \star^{\I}\nabla_{b_0})\in \Shv(\I\backslash \sL{G}/\IW)^{\leq 0}$ forces $n\leq 0$, or in other words $\sG\in \Shv(\I\backslash \Fl)^{\leq 0}$. 
\end{proof}

\section{Gaitsgory sheaves} 
\label{sec:Gaitsgory}

\subsection{Construction}

In this subsection we explain the construction of parabolic counterparts of the ``Gaitsgory sheaf'' studied in~\cite{Gaits18} (when $\k$ has characteristic $0$) and~\cite{adr}. A similar construction is also considered (when $\k$ has characteristic $0$) in~\cite{DL23}. (The authors of~\cite{DL23} consider a more general construction; ours corresponds to the case $\mu=0$ in their notation.)

We fix $J \subset S$, and work in the $\infty$-category $\Shv\left( \I^\quar_J \backslash \Gr \right)$. According to our general constructions the $\I^\quar_J$-orbits in $\Gr$ are naturally parametrized by $W_\ex^{(S,J)}$, but this set is in bijection with $W_\ex / W$ (by Lemma~\ref{lem:bi-min}\eqref{it:bi-min-1}), hence also with $\Xd$. For $\nu \in \Xd$ we will therefore set
\[
\mathrm{S}_\nu^J = \I^\quar_J z^\nu \sL^+G/\sL^+G,
\]
and denote by $\sigma^J_\nu$ the embedding of $\mathrm{S}_\nu^J$ in $\Gr$. (So, we have $\mathrm{S}^J_\nu = \fS^{\varnothing \subset J}_{S,x^J_\nu}$ where $x^J_\nu$ is the unique element in $t_\nu W \cap W_{\ex}^{(S,J)}$. When $J=\varnothing$, $\mathrm{S}_\nu^\varnothing$ is denoted $\mathrm{S}_\nu$ in~\cite{adr}, and when $J=S$ we are considering the $\I$-orbits on $\Gr$, denoted $(\Gr_\nu : \nu \in \Xd)$ in~\cite{adr}.)

For $\lambda\in \Xd^+$, we denote by $\sI_*(\lambda)$ the 0-th perverse cohomology of the $*$-pushforward of $\omega_{\sL^+G z^\lambda \sL^+G/\sL^+G}[-\langle2\rho,\lambda\rangle]$ to $\Gr$. (Here $\sI_*(\lambda)$ corresponds to the induced $G^\vee_\k$-module $\nabla(\lambda)$ of highest weight $\lambda$ under the geometric Satake equivalence~\eqref{eqn:geom-Satake}.)
For $\lambda\in \Xd^0_J\cap \Xd^+_{J^\comp}$ 
we set 
\[
\sM^J_\lambda:= \Ras^{\I\mid \I^\quar_J}\big( \sW_{t_{-\lambda}} \star^\I \oblv^{\sL^+G\mid \I}(\sI_*(\lambda)) \big) 
\quad \in\ \Shv\left( \I^\quar_J \backslash \Gr \right). 
\]
If $\lambda,\mu \in \Xd^0_J\cap \Xd^+_{J^\comp}$ satisfy $\mu \unlhd^J \lambda$ the construction explained in~\cite[\S 8.1]{adr} provides a canonical morphism 
\[
\theta^J_{\mu,\lambda}: \sM^J_\mu \rightarrow \sM^J_\lambda,
\]
and this construction provides a functor from $(\Xd^0_J\cap \Xd^+_{J^\comp},\unlhd^J)$ to $\Shv\left( \I^\quar_J \backslash \Gr \right)$. (To justify this one can either use the corresponding fact when $J=\varnothing$ proved in~\cite[\S 8.1]{adr} and apply the equivalence $\Ras^{\I^\semi | \I_J^\quar}$, or repeat the same argument as in~\cite{adr} using Lemma~\ref{lem:MlambdaJ-perverse} below.)

We define the \textit{$J$-Gaitsgory sheaf}
\[
\Gait^{J,\semi}:= \colim_{\mu\in \Xd^0_J\cap \Xd^+_{J^\comp}} \sM^J_\mu. 
\]
Note that for any $\lambda\in \Xd^0_J\cap \Xd^+_{J^\comp}$ we have equivalences 
\begin{equation}
\label{equ:Mlambda-Jsemi}
\begin{aligned}
\sW_{t_{-\lambda}} \star^\I \oblv^{\sL^+G\mid \I}(\sI_*(\lambda))
&\simeq \nabla_{J,t_{-\lambda}} \star^{\I_J} \oblv^{\sL^+G\mid \I_J}(\sI_*(\lambda)) \\ 
&\simeq \oblv^{\I_J \mid \I}( \nabla^J_{J,t_{-\lambda}} \star^{\I_J} \oblv^{\sL^+G\mid \I_J}(\sI_*(\lambda))),
\end{aligned}
\end{equation}
where the first equivalence follows from the equality $(p_J)_* \nabla_{t_{-\lambda}} \simeq \nabla_{J,t_{-\lambda}}$, and the second one from the fact that $\I z^{-\lambda}\I_J/\I_J=\I_J z^{-\lambda}\I_J/\I_J$. 
These equivalences are compatible with the morphisms $\theta_{\mu,\lambda}^J$,
hence $\sM^J_\lambda$ and $\Gait^{J,\semi}$ are $\I^\semi_J$-equivariant, and by abuse of notation we view 
these objects as belonging to $\Shv\left( \I^\semi_J \backslash \Gr \right)$.

\begin{rmk}
\phantomsection
\label{rmk:Gait-functions-P/L}
\begin{enumerate}
\item
\label{it:OP/L-colim}
To understand the definition of $\Gait^{J,\semi}$, one should compare it with the fact that if $P^\vee_{J,\k} \subset G^\vee_\k$ is the (positive) standard parabolic subgroup associated with $J$ and $M^\vee_{J,\k}$ its Levi subgroup containing the canonical maximal torus $T^\vee_\k$ then we have
\[
\mathscr{O}(P^\vee_{J,\k} / M^\vee_{J,\k}) \simeq \colim_{\lambda \in \Xd^0_J\cap \Xd^+_{J^\comp}} \nabla(\lambda) \otimes \k_{P^\vee_{J,\k}}(-\lambda)
\]
as $P^\vee_{J,\k}$-modules, for the natural transition maps in the colimit. In fact one should consider $\Gait^{J,\semi}$ as a geometric incarnation of the $P^\vee_{J,\k}$-module $\mathscr{O}(P^\vee_{J,\k} / M^\vee_{J,\k})$. This idea will be made precise in~\cite{adrs}.
\item
The definition of $\Gait^{J,\semi}$ given above matches the description given (for cha\-racteristic-$0$ coefficients) in~\cite[\S 4.1.9]{DL23}. See~\S\ref{ss:Gaitsgory-alternative} for a description which matches the definition chosen in~\cite{DL23}.
\item
We expect that the object $\Gait^{J,\semi}$ is related to the intersection cohomology complex of Drinfeld's compactification of the stack $\mathrm{Bun}_{P_J}$ in the same way as in the characteristic-$0$ setting treated in~\cite[\S 4.2]{DL23} or in the case $J=\varnothing$ treated in~\cite{adr-Drinfeld}. This will be discussed in future work.
\end{enumerate} 
\end{rmk}

\subsection{Main properties}

\begin{lem}
\label{lem:MlambdaJ-perverse}
For any $\lambda\in \Xd^0_J\cap \Xd^+_{J^\comp}$,
the object $\sM^J_\lambda$ is perverse in $\Shv\left( \I^\quar_J \backslash \Gr \right)$. 
Therefore, we have
\[
\Gait^{J,\semi}\in \Shv\left( \I^\semi_J \backslash \Gr \right)^\heartsuit . 
\]
\end{lem}

\begin{proof}
Using Proposition~\ref{prop:Ras-D}\eqref{it:Ras-D-1} and then Lemma~\ref{lem:push-pull-pK-D-N-semi}\eqref{it:push-pull-pK-D-N-semi-1}
we have equivalences 
\begin{align*}
\Ras^{\I\mid \I^\quar_J}\big( \sW_{t_{-\lambda}} \star^\I \oblv^{\sL^+G\mid \I}(\sI_*(\lambda)) \big)
&\simeq \Delta^{J,\quar}_{t_{-\lambda}} \star^\I \oblv^{\sL^+G\mid \I}(\sI_*(\lambda)) \\ 
&\simeq \Delta^{J,\quar}_{S,t_{-\lambda}} \star^{\sL^+G} \sI_*(\lambda),
\end{align*}
and by Lemma~\ref{lem:standard-perverse} and Proposition~\ref{prop:right-conv-t-exact} this object belongs to $\Shv\left( \I^\quar_J \backslash \Gr \right)^\heartsuit$. The second assertion follows, using compatibility of our t-structures with colimits. (Note that by Lemma~\ref{lem:oblv-semi-quar-texact}, the questions of whether our objects are perverse when considered as $\I_J^\quar$-equivariant or as $\I_J^\semi$-complexes are equivalent.)
\end{proof}

As in~\cite[\S 6.7]{adr} we consider
the $q$-analogue of Kostant's partition function $\sP(\nu,q)\in \Z[q]$, for any $\nu\in \Xd$. 
For $\nu\in \Xd$, any $\I^{\quar}_J$-equivariant perverse sheaf $\sF$ on $\mathrm{S}^J_\nu$ is a direct sum of copies of the appropriate shift of the dualizing complex; in case the multiplicity is finite we will denote it $\rk \, \sF$. Finally, we will denote by ${\check \rho} \in \frac{1}{2}\Xd$ the halfsum of the positive coroots.

\begin{lem}
\phantomsection
\label{lem:costk-Mlambda}
\begin{enumerate}
\item
\label{it:costk-Mlambda-2}
If $\lambda,\mu\in \Xd^0_J\cap \Xd^+_{J^\comp}$ satisfy $\mu \unlhd^J \lambda$, then the morphism $\theta^J_{\mu,\lambda}$ is injective. 
\item
\label{it:costk-Mlambda-1}
The objects $\sM^J_\lambda$ (for $\lambda\in \Xd^0_J\cap \Xd^+_{J^\comp}$) and $\Gait^{J,\semi}$ are supported on $\overline{\mathrm{S}^{J}_{0}}$. 
\item
\label{it:costk-Mlambda-3}
 If $\lambda\in \Xd^0_J\cap \Xd^+_{J^\comp}$ and if $\nu\in \Xd^+_J$ satisfies $\nu+\lambda\in \Xd^+$, then we have 
\[
\rk \, {}^{\mathrm{p}}\sH^0 \big( (\sigma^J_\nu)^!\sM^J_\lambda \big) = \begin{cases}
1 & \text{if $\nu=0$;} \\
0 & \text{otherwise},
\end{cases}
\quad
{}^{\mathrm{p}}\sH^1 \big( (\sigma^J_\nu)^!\sM^J_\lambda \big) = 0.
\]
\item
\label{it:costk-Mlambda-4} 
Assume that $\cha(\k)$ is good for $G$. 
If $\lambda\in \Xd^0_J\cap \Xd^+_{J^\comp}$ and if $\nu\in \Xd^+_J$ satisfies $\nu+\lambda\in \Xd^+$, then we have 
\[
\sum_{n \in \Z} \rk \, {}^{\mathrm{p}}\sH^n \big( (\sigma^J_\nu)^!\sM^J_\lambda \big) \cdot q^n 
= \sum_{w\in W} (-1)^{\ell(w)} \sP(w(\lambda+{\check \rho})- (\lambda+\nu+{\check \rho}), q^2). 
\]
\end{enumerate}
\end{lem}

\begin{proof}
\eqref{it:costk-Mlambda-2}
The image of $\theta^J_{\mu,\lambda}$ under the t-exact and conservative functor $\mathsf{F}_{\sL^+G}$ (see \textsection\ref{subsect:redtn-IW} and Theorem~\ref{thm:texactness}) is $\theta^\varnothing_{\mu,\lambda}$, which is injective by~\cite[Lemma~8.1]{adr}. Hence $\theta^J_{\mu,\lambda}$ is also injective.

\eqref{it:costk-Mlambda-1}
By~\eqref{it:costk-Mlambda-2} it suffices to prove that $\Gait^{J,\semi}$ is supported on $\overline{\mathrm{S}^{J}_{0}}$.
For $\lambda\in \Xd^0_J\cap \Xd^+_{J^\comp}$,
as in~\cite{adr} we will denote by $\overline{\nabla}_\lambda$ the costandard perverse sheaf on $\Gr$ associated with the $\I$-orbit of the image of $z^\lambda$, and
set $\sW^{\vee,\lambda}_S:=\sW_{t_{-\lambda}}\star^\I \overline{\nabla}_\lambda$. 
Since the restriction of $\sI_*(\lambda)$ to this $\I$-orbit is the perversely shifted dualizing sheaf, we have a natural morphism $\sI_*(\lambda)\rightarrow \overline{\nabla}_\lambda$, which induces a morphism 
\begin{equation}
\label{equ:Mlambda-pH0-dwaki-1}
\sM^J_\lambda \rightarrow {}^{\mathrm{p}}\sH^0\circ \Ras^{\I\mid \I^\quar_J}(\sW^{\vee,\lambda}_S) . 
\end{equation}

We claim that the morphism~\eqref{equ:Mlambda-pH0-dwaki-1} is injective. In fact, $\overline{\nabla}_\lambda$ is $\I_J$-equivariant since 
$\I z^{\lambda} \sL^+G/\sL^+G=\I_J z^{\lambda} \sL^+G/\sL^+G$. 
By a similar computation as in~\eqref{equ:Mlambda-Jsemi}, the right-hand side of~\eqref{equ:Mlambda-pH0-dwaki-1} is $\I^\semi_J$-equivariant, and~\eqref{equ:Mlambda-pH0-dwaki-1} is a morphism in $\Shv\left( \I^\semi_J \backslash \Gr \right)$. 
Applying the t-exact functor $\mathsf{F}_{\sL^+G}$ (see \textsection\ref{subsect:redtn-IW} and Theorem~\ref{thm:texactness}) to \eqref{equ:Mlambda-pH0-dwaki-1} and using Lemma~\ref{lem:Ras-shift-Wak}, we obtain a morphism
\begin{equation}
\label{equ:Mlambda-pH0-dwaki-2}
\sM^{\varnothing}_\lambda \rightarrow {}^{\mathrm{p}}\sH^0\circ \Ras^{\I\mid \I^\semi}(\sW^{\vee,\lambda}_S) 
\simeq {}^{\mathrm{p}}\sH^0\circ \Ras^{\I\mid \I^\semi}(\overline{\nabla}_\lambda)\langle -\lambda \rangle
\end{equation}
in $\Shv\left( \I^\semi \backslash \Gr \right)^\heartsuit$.
By \cite[Corollary 6.22]{adr}, the socle of $\sM^{\varnothing}_\lambda$ is the standard $\I^\semi$-equivariant perverse sheaf $\overline{\Delta}^\semi_{0}$ associated with the orbit of the base point. 
The precomposition of \eqref{equ:Mlambda-pH0-dwaki-2} with the embedding
\[
\Ras^{\I\mid \I^\semi}(\overline{\Delta}_{\lambda})\langle -\lambda \rangle \simeq \overline{\Delta}^\semi_{0} \hookrightarrow \sM^{\varnothing}_\lambda
\]
(where $\overline{\Delta}_{\lambda}$ is the standard perverse sheaf associated with the $\I$-orbit of the image of $z^\lambda$, and the equivalence follows from~\cite[Proposition~5.6]{adr})
coincides with the morphism induced by the natural map $\overline{\Delta}_{\lambda} \rightarrow \overline{\nabla}_{\lambda}$, and is therefore nonzero. (Here we use the fact that $\Ras^{\I | \I^\semi}$ is left t-exact by Corollary~\ref{cor:Ras-right-texact}, so that $\Ras^{\I\mid \I^\semi}(\overline{\nabla}_\lambda)$ belongs to $\Shv ( \I^\semi \backslash \Gr )^{\geq 0}$.)
Hence \eqref{equ:Mlambda-pH0-dwaki-2} is an injection, which implies that~\eqref{equ:Mlambda-pH0-dwaki-1} is injective as $\mathsf{F}_{\sL^+G}$ is t-exact and conservative. 

By~\eqref{equ:costd-K-colim} we have 
\[
\nabla^{J,\quar}_{S,e}\simeq \colim_{\mu\in \Xd^0_J\cap \Xd^+_{J^\comp}} \Ras^{\I\mid \I^\quar_J}(\sW^{\vee,\mu}_S). 
\]
(Note that $e \in W_\ex^{(S,J)}$, so that the left-hand side makes sense.)
By our claim above, the colimit of the maps~\eqref{equ:Mlambda-pH0-dwaki-1} gives an injection 
\[
\Gait^{J,\semi} \hookrightarrow {}^{\mathrm{p}}\sH^0(\nabla^{J,\quar}_{S,e}), 
\]
from which~\eqref{it:costk-Mlambda-1} follows. 

\eqref{it:costk-Mlambda-3}--\eqref{it:costk-Mlambda-4}
The proofs are similar to those of the corresponding statements in~\cite[Lemma 6.18]{adr}; we only sketch the argument. The rank of ${}^{\mathrm{p}}\sH^n \big( (\sigma^J_\nu)^!\sM^J_\lambda \big)$ can be computed as the dimension of the morphism space in $\Shv(\I_{J,\mathrm{u}}^\quar \backslash \Gr)$ (where $\I_{J,\mathrm{u}}^\quar$ is the semidirect product of the pro-unipotent radical of $\I_{\varnothing \subset J}$ with $\sL U_J$) between the images of the standard $\I_J^\quar$-equivariant complex $\overline{\Delta}^{J,\quar}_{\nu}$ associated with $\mathrm{S}^J_\nu$ and $\sM^J_\lambda[n]$. Applying a pro-unipotent version of the Raskin equivalence and Lemma~\ref{lem:Ras-shift-Wak}, this space identifies with the space of morphisms in $\Shv(\Iu \backslash \Gr)$ between the images of $\Ras^{\I^\quar_J | \I}(\overline{\Delta}^{J,\quar}_{\lambda+\nu})$ and $\sI_*(\lambda)[n]$. (Here, $\Iu$ is the pro-unipotent radical of $\I$.) Since $\lambda+\nu$ is dominant, $t_{\lambda+\nu}$ belongs to $W_\ex^{(S,J)}$, and by Lemma~\ref{lem:push-pull-pK-D-N-semi}\eqref{it:push-pull-pK-D-N-semi-1} we have $\overline{\Delta}^{J,\quar}_{\lambda+\nu}=(p_S)_* \Delta^{J,\quar}_{t_{\lambda+\nu}}$. Using Proposition~\ref{prop:Ras-D}\eqref{it:Ras-D-1} we deduce that $\Ras^{\I^\quar_J | \I}(\overline{\Delta}^{J,\quar}_{\lambda+\nu})$ is the object denoted $\overline{\sW}_{\lambda+\nu}$ in~\cite{adr}. With this identification, the required dimensions are computed in~\cite{adr} and given as in the statement.
\end{proof}

\begin{rmk}
The restriction on the characteristic in Lemma~\ref{lem:costk-Mlambda}\eqref{it:costk-Mlambda-4} comes from the use of the Mirkovi{\'c}--Vilonen conjecture on costalks of the objects $\sI_*(\lambda)$, which is known to be true for all $\lambda \in \Xd^+$ if and only if $\mathrm{char}(\k)$ is good for $G$ (see~\cite[Proof of Lemma~6.18]{adr} for references). In fact, here we only need this property when $\lambda \in \Xd^0_J \cap \Xd_{J^\comp}^+$, so we might hope to relax this assumption. For instance, if $J=S$ it is easily seen that we do not need any restriction on $\mathrm{char}(\k)$. But we do not know how to implement this idea when $J \neq S$.
\end{rmk}

In the following proposition we consider the object $\nabla^{J,\semi}_{S,e}$. Note first that $e \in {}^J W_\ex^{(S,J)}$, so that this object makes sense. Note also that $\I_J^\semi \sL^+G = \I_J^\quar \sL^+ G$; so this object coincides with the object $\nabla^{J,\quar}_{S,e}$ considered in the proof of Lemma~\ref{lem:costk-Mlambda}.

\begin{prop}
\label{prop:costk-Ga}
\begin{enumerate}
\item
\label{it:costk-Ga-1} 
We have an isomorphism $\Gait^{J,\semi} \simeq {}^{\mathrm{p}}\sH^0(\nabla^{J,\semi}_{S,e})$. 
\item
\label{it:costk-Ga-2} 
Assume $\cha(\k)$ is good for $G$. 
For any $\nu\in \Xd^+_J$, we have 
\begin{equation*}
\sum_{n} \rk\, {}^{\mathrm{p}}\sH^n\big( (\sigma^J_\nu)^!\Gait^{J,\semi} \big) \cdot q^n 
= \sum_{w\in W_J} (-1)^{\ell(w)} \sP(w({\check \rho})-{\check \rho} -\nu, q^2). 
\end{equation*}
\end{enumerate}
\end{prop}

\begin{proof}
Here again the proof is similar to that of the counterpart of this statement in~\cite[\S 8.2]{adr}, so we will be brief. 
The claims follow from the computation of the rank of perverse cohomology objects of $!$-restrictions to strata of $\Gait^{J,\semi}$. In fact, since this object is $\I_J^\semi$-equivariant, and since $\mathrm{S}^J_\nu$ is open in the $\I_J^\semi$-orbit containing it when $\nu \in \Xd^+_J$, it suffices to consider the $!$-restrictions $(\sigma^J_\nu)^!\Gait^{J,\semi}$ when $\nu \in \Xd^+_J$. We compute the rank in this case 
as in the proof of Lemma~\ref{lem:costk-Mlambda}\eqref{it:costk-Mlambda-3}--\eqref{it:costk-Mlambda-4}, as the dimension of a certain morphism space in $\Shv(\I_{J,\mathrm{u}}^\quar \backslash \Gr)$. Since $\Gait^{J,\semi}$ is defined as a colimit and since the images of costandard objects in $\Shv(\I_{J,\mathrm{u}}^\quar \backslash \Gr)$ are compact, these spaces are colimits of their counterparts in $\sM_\lambda^J$, and to conclude we need to show that for fixed $\nu$, the latter spaces are independent of $\lambda$ for $\lambda$ ``large enough.'' 

This fact can be reduced to a statement which follows from the same considerations as in the proof of~\cite[Theorem~8.2]{adr}. Namely, using a ``pro-unipotent'' analogue of $\Ras^{\I_J^\quar | \I^\semi}$, it suffices to show that for any fixed $\sF$ in $\Shv_\omega(\I_{\mathrm{u}}^\semi \backslash \Gr)$ the space of morphisms from $\sF$ to the image of $\sM_\lambda^\varnothing$ is independent of $\lambda$ for $\lambda$ ``large enough.'' (See~\cite[\S 6.4]{adr} for the notation used here.) Any such $\sF$ is a finite extension of cohomological shifts of images of standard perverse sheaves, so one can assume that $\sF$ is such an object. Then one proceeds as in the proof of~\cite[Theorem~8.2]{adr}.
(Compared to that proof, in the present setting, all the elements $w \in W_J$ will contribute to the multiplicities, which also explains the sum over $W_J$ appearing here compared to the setting of~\cite{adr}.)
\end{proof}

\begin{rmk}
Comparing the arguments above with the formula in Remark~\ref{rmk:Gait-functions-P/L}\eqref{it:OP/L-colim}, one sees that for $\nu \in \Xd_J^+$ the integer $\sum_{w\in W_J} (-1)^{\ell(w)} \sP(w({\check \rho})-{\check \rho} -\nu, 1)$ is the dimension of the $T^\vee_\k$-weight space of weight $\nu$ in $\mathscr{O}(P^\vee_{J,\k} / M^\vee_{J,\k})$.
\end{rmk}

\subsection{Alternative description}
\label{ss:Gaitsgory-alternative}

As in~\cite[\S 8.3]{adr} one can describe the object $\Gait^{J,\semi}$ as a different colimit. Namely, generalizing slightly the notation introduced in \S\ref{ss:convolution}, for $\mu \in \Xd$ we denote by $z^\mu \cdot (-)$ the pushforward functor under the automorphism of $\Gr$ given by the left action of $z^\mu \in \sL G(\F)$. It is explained in~\cite[\S 8.3]{adr} (following Gaitsgory) that the assignment $\lambda \mapsto z^{-\lambda} \cdot \sI_*(\lambda)[\langle \lambda,2\rho \rangle]$ defines a functor from $\Xd^+$ (endowed with the restriction of the preorder $\unlhd^\varnothing$) to $\Shv(T \backslash \Gr)$. Here we restrict this functor to $\Xd_J^0 \cap \Xd^+_{J^\comp}$. (Note that the restriction of $\unlhd^\varnothing$ to this subset coincides with the restriction of $\unlhd^J$.)

The following lemma follows from the same considerations as for~\cite[Proposition~8.6]{adr}.

\begin{lem}
\label{lem:Gait-alternative}
There exists a canonical equivalence
\[
\Gait^{J,\semi} \simeq \colim_{\lambda \in \Xd_J^0 \cap \Xd^+_{J^\comp}} z^{-\lambda} \cdot \sI_*(\lambda)[\langle \lambda,2\rho \rangle].
\]
\end{lem}

Consider the affine Grassmannians $\Gr_{M_J}$ and $\Gr_{P_J}$ of $M_J$ and $P_J$, and the diagram
\[
\begin{tikzcd}
\Gr_{M_J} \ar[r, "b"',shift right=0.5ex] & \Gr_{P_J} \ar[r, "c"] \ar[l, "a"', shift right=0.5ex] & \Gr
\end{tikzcd}
\]
where the morphisms are induced by the natural maps $P_J \to G$, $P_J \to M_J$ and $M_J \to P_J$. (This setting is similar to that considered in the proof of Proposition~\ref{prop:exactness-left-conv}, but for $\Gr$ instead of $\Fl$, so we use similar notation.) The connected components in $\Gr_{M_J}$ and $\Gr_{P_J}$ are naturally parametrized by $\Xd / \Z {\check \Phi}_J$; for $\zeta \in \Xd / \Z {\check \Phi}_J$ we will denote by $\Gr_{M_J, \zeta}$ and $\Gr_{P_J,\zeta}$ the corresponding components, and by $a_\zeta$, $b_\zeta$, $c_\zeta$ the restrictions of the maps above to these components. As in~\cite[Lemma~3.3.2]{DL23}, the functors $(a_\zeta)^!$ and $(b_\zeta)^!$ induce mutually inverse equivalences
\begin{equation}
\label{eqn:equiv-Gr-P-M}
\Shv \left( \I_J^\semi \backslash \Gr_{P_J, \zeta} \right) \simeq \Shv \left( \sL^+ M_J \backslash \Gr_{M_J, \zeta} \right).
\end{equation}

Let us now denote by $\Sat$ the equivalence~\eqref{eqn:geom-Satake}, and by $\Sat_{J}$ its counterpart of $M_J$. It is well known that for any $\zeta \in \Xd / \Z {\check \Phi}_J$ the functor
\[
(a_\zeta)_! (c_\zeta)^* [\langle \zeta, 2\rho-2\rho_J \rangle] : \Shvfg \left( \sL^+ G \backslash \Gr \right) \to \Shvfg \left( \sL^+ M_J \backslash \Gr_{M_J,\zeta} \right)
\]
is t-exact, and corresponds under the equivalences $\Sat$ and $\Sat_J$ to the composition of the functor $(-)_{M^\vee_{J,\k}}$ of restriction of the action of $G^\vee_\k$ to $M^\vee_{J,\k}$ with the functor $(-)_{\mathrm{Z}_J=\zeta}$ of projection to the direct summand where the center of $M_{J,\k}^\vee$ acts via the character defined by $\zeta$.

The following statement generalizes~\cite[Proposition~4.1.12]{DL23} and~\cite[Remark~8.3]{adr}. Here, for $\mu \in \Xd^0_J$ we denote by $\k_{M^\vee_{J,\k}}(\mu)$ the $1$-dimensional $M^\vee_{J,\k}$-module determined by $\mu$.

\begin{prop}
For any $\zeta \in \Xd / \Z {\check \Phi}_J$ we have an equivalence
\[
(c_\zeta)^* \Gait^{J,\semi} \simeq (a_\zeta)^! \left( \colim_{\lambda \in \Xd_J^0 \cap \Xd^+_{J^\comp}} \Sat_J \left( (\k_{M^\vee_{J,\k}}(-\lambda) \otimes \nabla(\lambda)_{| M^\vee_{J,\k}})_{\mathrm{Z}_J=\zeta} \right) \right)[- \langle \zeta, 2\rho-2\rho_J \rangle]
\]
where the maps in the colimit are induced by those considered in Remark~\ref{rmk:Gait-functions-P/L}\eqref{it:OP/L-colim}.
\end{prop}

\begin{proof}
The proof is identical to that of~\cite[Proposition~4.1.12]{DL23}. Namely, using the equivalence~\eqref{eqn:equiv-Gr-P-M}, it suffices to prove that
\[
(b_\zeta)^! (c_\zeta)^* \Gait^{J,\semi} \simeq \colim_{\lambda \in \Xd_J^0 \cap \Xd^+_{J^\comp}} \Sat_J \left( (\k_{M^\vee_{J,\k}}(-\lambda) \otimes \nabla(\lambda)_{| M^\vee_{J,\k}})_{\mathrm{Z}_J=\zeta} \right) [- \langle \zeta, 2\rho-2\rho_J \rangle].
\]
Now, using Lemma~\ref{lem:Gait-alternative} we obtain that
\[
(b_\zeta)^! (c_\zeta)^* \Gait^{J,\semi} \simeq \colim_{\lambda \in \Xd_J^0 \cap \Xd^+_{J^\comp}} z^{-\lambda} \cdot (b_{\overline{\lambda}+\zeta})_! (c_{\overline{\lambda}+\zeta})^* \sI_*(\lambda)[\langle \lambda,2\rho \rangle]
\]
where $\overline{\lambda}$ is the image of $\lambda$ in $\Xd / \Z {\check \Phi}_J$, and we still use the notation $z^{-\lambda} \cdot (-)$ for the counterpart of this functor for sheaves on $\Gr_{M_J}$. By the comments above we have
\[
(b_{\overline{\lambda}+\zeta})^! (c_{\overline{\lambda}+\zeta})^* \sI_*(\lambda) \simeq \Sat_J \left( (\nabla(\lambda)_{| M^\vee_{J,\k}})_{\mathrm{Z}_J=\overline{\lambda}+\zeta} \right) [-\langle \overline{\lambda}+\zeta, 2\rho-2\rho_J \rangle],
\]
and we deduce that
\begin{multline*}
z^{-\lambda} \cdot (b_{\overline{\lambda}+\zeta})^! (c_{\overline{\lambda}+\zeta})^* \sI_*(\lambda)[\langle \lambda,2\rho \rangle] \\
\simeq \Sat_J \left( (\k_{M^\vee_{J,\k}}(-\lambda) \otimes \nabla(\lambda)_{| M^\vee_{J,\k}})_{\mathrm{Z}_J=\zeta} \right) [- \langle \zeta, 2\rho-2\rho_J \rangle],
\end{multline*}
which finishes the proof.
\end{proof}

\bibliographystyle{plain} 
\bibliography{MyBibtex}

\end{document}